\documentclass[11pt,reqno]{amsart}
\usepackage{geometry}
\usepackage{graphicx}
\usepackage{amssymb}
\usepackage{caption}
\usepackage{subcaption}
\usepackage{enumitem}
\usepackage{mathtools}
\usepackage[all]{xy}
\usepackage{xcolor}
\usepackage{soul}
\usepackage[hidelinks,colorlinks=true,linkcolor=blue,urlcolor=blue,citecolor=blue]{hyperref}

\usepackage{microtype}

\usepackage{bbold}

\usepackage{MnSymbol}

\usepackage{listings}

\definecolor{backcolour}{rgb}{0.62,0.95,0.99}
\definecolor{commentcolour}{rgb}{0,0.6,0}

\usepackage{mdframed}

\lstdefinestyle{mystyle}{
	backgroundcolor=\color{backcolour},   
	commentstyle=\color{commentcolour},
	basicstyle=\ttfamily\small,               
	captionpos=b,                    
	keepspaces=true,                              
	showtabs=false,                  
	tabsize=5,
}

\usepackage[linesnumbered,ruled,vlined]{algorithm2e}

\SetCommentSty{mycommfont}
\SetKwInput{KwInput}{Input}               
\SetKwInput{KwOutput}{Output}     
\let\oldnl\nl% Store \nl in \oldnl
\newcommand{\nonl}{\renewcommand{\nl}{\let\nl\oldnl}}

\usepackage{framed}

\usepackage{graphicx}

\graphicspath{{figures/}}

\numberwithin{equation}{section}

\newtheorem{theorem}{Theorem}[section]

\newtheorem{remark}{Remark}[section]
\newtheorem*{remark*}{Remark}

\newtheorem{problem}{Problem}

\usepackage[numbers,sort&compress]{natbib}
\bibpunct{[}{]}{;}{n}{,}{,}

\usepackage{dutchcal}

\DeclarePairedDelimiterX{\inner}[2]{\langle}{\rangle}{#1, #2}

\DeclareMathOperator*{\ext}{ext}

\allowdisplaybreaks

\title[Practical Perspectives on Symplectic Accelerated Optimization]{Practical Perspectives on \\ Symplectic Accelerated Optimization}
\author{Valentin Duruisseaux and Melvin Leok}

\begin{document}

\maketitle

\begin{abstract}  
	Geometric numerical integration has recently been exploited to design symplectic accelerated optimization algorithms by simulating the Bregman Lagrangian and Hamiltonian systems from the variational framework introduced in~\citeauthor{WiWiJo16}. In this paper, we discuss practical considerations which can significantly boost the computational performance of these optimization algorithms, and considerably simplify the tuning process. In particular, we investigate how momentum restarting schemes ameliorate computational efficiency and robustness by reducing the undesirable effect of oscillations, and ease the tuning process by making time-adaptivity superfluous. We also discuss how temporal looping helps avoiding instability issues caused by numerical precision, without harming the computational efficiency of the algorithms. Finally, we compare the efficiency and robustness of different geometric integration techniques, and study the effects of the different parameters in the algorithms to inform and simplify tuning in practice. From this paper emerge symplectic accelerated optimization algorithms whose computational efficiency, stability and robustness have been improved, and which are now much simpler to use and tune for practical applications.
\end{abstract}

\hfill 

\section{Introduction}
The field of symplectic optimization grew out of efforts to generalize Nesterov's accelerated gradient method~\cite{Nes83}, which was shown to converge in $\mathcal{O}(1/k^2)$ to the minimum of the convex objective function $f$ and improves on the $\mathcal{O}(1/k)$ convergence rate exhibited by standard gradient descent methods. This $\mathcal{O}(1/k^2)$ convergence rate, referred to as acceleration, was shown in~\cite{Nes04} to be optimal among first-order methods using only information about $\nabla f$ at consecutive iterates. Nesterov's algorithm was shown in~\cite{SuBoCa16} to limit to a second-order ordinary differential equation as the timestep goes to 0, and that $f(x(t))$ converges to its optimal value at a rate of $\mathcal{O}(1/t^2)$ along any trajectory $x(t)$ of this ODE. It was then shown in~\cite{WiWiJo16} that in continuous time, an arbitrary convergence rate $\mathcal{O}(1/t^p)$ can be achieved in normed spaces, by considering flow maps generated by a family of time-dependent Bregman Lagrangian and Hamiltonian systems which is closed under time-rescaling. This lead to the field of symplectic optimization~\cite{Jordan2018}, where symplectic discretizations of the Bregman Hamiltonian flow are used to construct accelerated optimization algorithms.

Lagrangian and Hamiltonian flows can also be described variationally. This, together with the time-rescaling property of this family, were exploited in~\cite{duruisseaux2020adaptive} by using time-adaptive geometric integrators to design efficient explicit algorithms for symplectic accelerated optimization. It was observed that a careful use of adaptivity and symplecticity could result in a significant gain in computational efficiency. There has also been work on deriving accelerated optimization algorithms in the Riemannian manifold setting~\cite{Duruisseaux2022Riemannian,Duruisseaux2022Constrained,Duruisseaux2022Projection,Duruisseaux2022Lagrangian,alimisis2020,Alimisis2021,Alimisis2020-1,Sra2016, Sra2018, Sra2020, Liu2017}. 

While the symplectic optimization approach provides a broad framework for constructing accelerated optimization algorithms, the real-world performance of these methods depends on the choice of numerous parameters. In this paper, we will perform a systematic and comprehensive test of a class of symplectic accelerated optimization algorithms, so as to provide practical guidance on how to achieve good real-world performance with less tuning.

\subsection*{Outline of the paper.} After reviewing the basics of geometric integration in Section~\ref{section: Geometric Mechanics and Geometric Numerical Integration}, we introduce variational accelerated optimization and present how to integrate the corresponding dynamics in Sections~\ref{section: Variational Accelerated Optimization} and~\ref{section: Numerical Methods and Problems of Interest}. We then analyze the oscillatory behavior of these dynamical systems in Section~\ref{section: Controlling the Oscillatory Behavior}, and discuss how their unfavorable effect can be neutralized, in particular via the use of momentum restarting techniques which can dramatically improve computational efficiency and robustness. In Section~\ref{section: Time-Adaptivity in the Momentum Restarted Algorithms}, we show that momentum restarting makes time-adaptivity futile, which allows us to simplify the algorithms. Then, in Sections~\ref{sec: Comparison Integrators} and~\ref{section: Tuning the Algorithms}, we compare the different geometric integrators, and investigate how the computational performance depends on the different parameters, which allows us to reduce the numbers of parameters to tune in practice. In Section~\ref{section: Temporal Looping to Improve Numerical Stability}, we see that temporal looping can avoid instability issues due to numerical precision, and finally in Section~\ref{section: ML Testing}, we test the resulting algorithms on problems of interest to the machine learning community. \\

\section{Geometric Mechanics and Geometric Numerical Integration} \label{section: Geometric Mechanics and Geometric Numerical Integration}

\subsection{Lagrangian and Hamiltonian Mechanics}

\hfill  \\

Given a manifold $\mathcal{Q}$, a \textbf{Lagrangian} is a function $L:T\mathcal{Q}  \rightarrow \mathbb{R}$. The corresponding action integral $\mathcal{S}$ is the functional
	\begin{equation} 
		\mathcal{S} (q) = \int_{0}^{T}{L(q,\dot{q})dt},
	\end{equation}
	over the space of smooth curves $q:[0,T] \rightarrow \mathcal{Q}$. Hamilton's variational principle states that $ \delta \mathcal{S}=0$ where the variation $\delta \mathcal{S}$ is induced by an infinitesimal variation $\delta q$ of the trajectory $q$ that vanishes at the endpoints. Given local coordinates $(q^1, \ldots , q^n)$ on the manifold $\mathcal{Q}$, Hamilton's variational principle can be shown to be equivalent to the \textbf{Euler--Lagrange equations},
	\begin{equation}\label{eq: EL Basic}
		\frac{d}{dt} \left( \frac{\partial L} {\partial \dot{q}^k} \right)= \frac{\partial L }{\partial q^k}, \qquad  \text{for } k=1,\ldots , n. 
	\end{equation}
A Lagrangian $L$ is hyperregular if the Legendre transform $\mathbb{F}L : T\mathcal{Q} \rightarrow T^* \mathcal{Q}$ of $L$, defined fiberwise by $\mathbb{F}L : (q^i,\dot{q}^i) \mapsto \left(q^i,\frac{\partial L}{\partial \dot{q}^i} \right)$, is diffeomorphic. A hyperregular Lagrangian on $T\mathcal{Q}$ induces a \textbf{Hamiltonian} system on $T^* \mathcal{Q}$ via 
	\begin{equation}
		H(q,p) = \langle \mathbb{F}L (q,\dot{q}) , \dot{q} \rangle - L(q,\dot{q})= \sum_{j=1}^{n} {p_j  \dot{q}^j} - L(q,\dot{q}) \bigg|_{p_i=\frac{\partial L}{\partial \dot{q}^i}},
	\end{equation} 
	where $p_i = \frac{\partial L}{\partial \dot{q}^i} \in T^* \mathcal{Q}$ is the conjugate momentum of $q^i$. There is a Hamiltonian variational principle on the Hamiltonian side in momentum phase space which is equivalent to \textbf{Hamilton's equations}, 
	\begin{equation} \label{eq: Hamilton Equations Basic} \dot{p}_k = -\frac{\partial H}{\partial q^k} (p,q),  \qquad  \dot{q}^k = \frac{\partial H}{\partial p_k} (p,q), \qquad \text{for } k=1,\ldots , n,
	\end{equation}
	and these equations are equivalent to the Euler--Lagrange equations \eqref{eq: EL Basic}, provided the Lagrangian is hyperregular. Hamiltonian systems possess a long list of structural invariants and constants of motion, the most important of which are the conservation of the Hamiltonian energy and the conservation of the symplectic 2-form. \\

	\subsection{Symplectic and Variational Integrators} 
	
	\hfill  \\
	
	Symplectic integrators form a class of geometric numerical integrators of interest since, when applied to Hamiltonian systems, they yield discrete approximations of the flow that preserve the symplectic 2-form. The preservation of the symplectic 2-form results in the preservation of many qualitative aspects of the underlying dynamical system. In particular, the numerical solution of a Hamiltonian system obtained using a constant time-step symplectic integrator is exponentially-near to the exact solution of a nearby Hamiltonian system for an exponentially-long time~\citep{Benettin1994, HaLuWa2006}. It explains why symplectic integrators exhibit good energy conservation with essentially no accumulation of errors in time, when applied to Hamiltonian systems, and why symplectic methods are best suited to integrate Hamiltonian systems.  We refer the reader to~\citep{IserlesWhyGNI} for a brief recent overview of geometric numerical integration, and to~\citep{HaLuWa2006,Blanes2017,LeRe2005} for a more comprehensive presentation of structure-preserving integration techniques.

Variational integrators form a class of symplectic integrators, derived by discretizing Hamilton's principle instead of discretizing Hamilton's equations directly. As a result, variational integrators are symplectic, preserve many invariants and momentum maps, and have excellent long-time near-energy preservation~\cite{MaWe2001}. Traditionally, variational integrators have been designed based on the Type~I generating function known as the discrete Lagrangian, $L_d:Q \times Q \rightarrow \mathbb{R}$. The exact discrete Lagrangian that generates the time-$h$ flow of Hamilton's equations can be represented both in a variational form and boundary-value form. The latter is given by \begin{equation}
	L_d^E(q_0,q_1;h)=\int_0^h L(q(t),\dot q(t)) dt, \end{equation} where $q(0)=q_0,$ $q(h)=q_1,$ and $q$ satisfies the Euler--Lagrange equations over the time interval $[0,h]$. A variational integrator is defined by constructing an approximation $L_d:Q \times Q \rightarrow \mathbb{R}$ to $L_d^E$, and then applying the discrete Euler--Lagrange equations,
\begin{equation}
	p_k=-D_1 L_d(q_k, q_{k+1}),\qquad p_{k+1}=D_2 L_d(q_k, q_{k+1}),  \label{IDEL}
\end{equation}
where $D_i$ denotes a partial derivative with respect to the $i$-th argument. %, and these equations implicitly define the integrator $\tilde{F}_{L_d}:(q_k,p_k)\mapsto(q_{k+1},p_{k+1})$. 
The error analysis is greatly simplified via Theorem~2.3.1 of~\cite{MaWe2001}, which states that if a discrete Lagrangian, $L_d:Q\times Q\rightarrow\mathbb{R}$, approximates the exact discrete Lagrangian $L_d^E:Q\times Q\rightarrow\mathbb{R}$ to order $r$, i.e., \begin{equation}
		L_d(q_0, q_1;h)=L_d^E(q_0,q_1;h)+\mathcal{O}(h^{r+1}) ,
\end{equation}
then the discrete Hamiltonian map $\tilde{F}_{L_d}:(q_k,p_k)\mapsto(q_{k+1},p_{k+1})$ defined by \eqref{IDEL} and viewed as a one-step method, has order of accuracy $r$. Many properties of the integrator can be determined by analyzing the associated discrete Lagrangian, as opposed to analyzing the integrator directly. 

Variational integrators have been extended to the framework of Type~II/III generating functions, referred to as discrete Hamiltonians~\cite{LaWe2006, LeZh2011,ScLe2017}. Hamiltonian variational integrators are derived by discretizing Hamilton's phase space principle. The boundary-value formulation of the Type~II generating function of the Hamiltonian flow is given by the exact discrete right Hamiltonian, 
\begin{equation}
	H_d^{+,E}(q_0,p_1;h) =  p_1^\top q_1 - \int_0^h \left[ p(t)^\top \dot{q}(t)-H(q(t), p(t)) \right] dt, \end{equation}  where $(q,p)$ satisfies Hamilton's equations with boundary conditions $q(0)=q_0$ and $p(h)=p_1$. A Type~II Hamiltonian variational integrator is constructed by using an approximate discrete Hamiltonian~$H_d^+$, and applying the discrete right Hamilton's equations
\begin{equation}\label{Discrete Right Eq}
	p_0=D_1H_d^+(q_0,p_1), \qquad q_1=D_2H_d^+(q_0,p_1).
\end{equation}
%which implicitly defines the integrator, $\tilde{F}_{H_d^+}:(q_0,p_0) \mapsto (q_1,p_1)$. 
Theorem~2.3.1 of~\cite{MaWe2001}, which simplifies the error analysis for Lagrangian variational integrators, has an analogue for Hamiltonian variational integrators. Theorem~2.2 in~\cite{ScLe2017} states that if a discrete right Hamiltonian $H^+_d$ approximates the exact discrete right Hamiltonian $H_d^{+,E}$ to order $r$, i.e., 
\begin{equation}
	H^+_d(q_0, p_1;h)=H_d^{+,E}(q_0,p_1;h)+\mathcal{O}(h^{r+1}), \end{equation}  then the discrete right Hamiltonian map $\tilde{F}_{H^+_d}:(q_k,p_k)\mapsto(q_{k+1},p_{k+1})$ defined by \eqref{Discrete Right Eq} and viewed as a one-step method, is order $r$ accurate. Note that discrete left Hamiltonians and corresponding discrete left Hamilton's maps can also be constructed in the Type~III case (see~\cite{LeZh2011,duruisseaux2020adaptive}). \\ 

Examples of variational integrators include Galerkin variational integrators ~\cite{LeZh2011,MaWe2001}, Taylor variational integrators~\cite{ScShLe2017}, and prolongation-collocation variational integrators~\cite{LeSh2011}. In this paper, we will use Taylor variational integrators, where a discrete approximate Lagrangian or Hamiltonian is constructed by approximating the flow map and the trajectory associated with the boundary values using a Taylor method, and approximating the integral by a quadrature rule. The Taylor variational integrator is generated by the implicit discrete Euler--Lagrange equations associated to the discrete Lagrangian or by the Hamilton's equations associated with the discrete Hamiltonian. The construction of Taylor variational integrator is presented in the context of accelerated optimization in~\cite{duruisseaux2020adaptive,Duruisseaux2022Lagrangian}. 
	
	 In many cases, the Type~I and Type~II/III approaches produce equivalent integrators, such as for Taylor variational integrators provided the Lagrangian is hyperregular~\cite{ScShLe2017}. However, Hamiltonian and Lagrangian variational integrators are not always equivalent in practice, even when they are analytically equivalent, as they might still have different numerical properties because of numerical conditioning issues~\cite{ScLe2017}. Even more to the point, Lagrangian variational integrators cannot always be constructed when the underlying Hamiltonian is degenerate, which is the case in the adaptive Hamiltonian framework for accelerated optimization presented in Section~\ref{subsubsec: Time-adaptive Hamiltonian Integrators}.  \\

\section{Variational Framework for Accelerated Optimization}
\label{section: Variational Accelerated Optimization}

\hfill 

\subsection{General Framework}
  
  \hfill \\
  
 Efficient optimization has become one of the major concerns in data analysis. Many machine learning algorithms are designed around the minimization of a loss function or the maximization of a likelihood function. Due to the ever-growing size of data sets and problems, there has been a lot of focus on first-order optimization algorithms because of their low cost per iteration, and many gradient-based optimization methods have been proposed since Cauchy's first gradient descent algorithm~\citep{Cauchy1847}. Nesterov's Accelerated Gradient (NAG) method
  \begin{equation}\label{eq: NesterovUpdate}
  	x_k = y_{k-1} -h\nabla f(y_{k-1}) , \qquad y_k = x_k + \frac{k-1}{k+2} (x_k - x_{k-1}), 
  \end{equation}
was introduced in 1983 in~\citep{Nes83}, and converges in $\mathcal{O}(1/k^2)$ to the minimum of the convex objective function $f$, improving on the $\mathcal{O}(1/k)$ convergence rate exhibited by the standard gradient descent methods.
  This $\mathcal{O}(1/k^2)$ convergence rate was shown in~\citep{Nes04} to be optimal among first-order methods using only information about $\nabla f$ at consecutive iterates. This phenomenon in which an algorithm displays this improved rate of convergence is referred to as acceleration, and other accelerated algorithms have been derived, such as accelerated mirror descent~\citep{Nem1983}, and accelerated cubic-regularized Newton's method~\citep{Nes08}. 

It was shown in~\citep{SuBoCa16} that Nesterov's method limits to a second-order ODE, as the step size goes to 0. The authors also proved that the objective function $f(x(t))$ converges to its optimal value at a rate of $\mathcal{O}(1/t^2)$ along the trajectories of this ODE. It was then shown in~\citep{WiWiJo16} that in continuous time, the convergence rate of $f(x(t))$ can be accelerated to an arbitrary rate $\mathcal{O}(1/t^p)$, by considering flow maps generated by a family of time-dependent Bregman Lagrangian and Hamiltonian systems on normed vector spaces which is closed under time rescaling.  More precisely, in a general space $\mathcal{Q}$, given a convex, continuously differentiable function $h:\mathcal{Q} \rightarrow \mathbb{R}$ such that $\Vert \nabla h(q) \Vert \rightarrow \infty $ as $\Vert q \Vert \rightarrow \infty $, its corresponding Bregman divergence is given by
  \begin{equation}  D_h(x,y) = h(y)-h(x) - \langle \nabla h(x), y-x \rangle .  \end{equation} 
   The \textbf{Bregman Lagrangian and Hamiltonian} are defined as
  \begin{equation} \label{eq: Bregman L General Normed}
  	L_{\alpha,\beta,\gamma}(q,v,t) = e^{\alpha_t + \gamma_t} \left[   D_h(q+e^{-\alpha_t}v , q) - e^{\beta_t} f(q)  \right]  ,
  \end{equation}
  \begin{equation} \label{eq: Bregman H General Normed}
  	H_{\alpha,\beta,\gamma}(q,r,t) = e^{\alpha_t + \gamma_t} \left[   D_{h^*}(\nabla h(q)+e^{-\gamma_t}r , \nabla h(q) ) + e^{\beta_t} f(q)  \right] ,
  \end{equation}
  which are scalar-valued functions of position $q\in \mathcal{Q}$, velocity $v\in \mathbb{R}^d$,  momentum $r\in \mathbb{R}^d$, and time $t$, and are parametrized by smooth functions of time, $\alpha,\beta,\gamma$. Here, the function $h^*:\mathcal{Q}^* \rightarrow \mathbb{R}$ denotes the Legendre transform (or convex dual function) of $h$, defined by $h^*(w) = \sup_{z\in \mathcal{Q}}{\left[ \langle w, z \rangle - h(z) \right]}$. These parameter functions $\alpha,\beta,\gamma$ are said to satisfy the ideal scaling conditions if
  \begin{equation}\label{eq: IdealScaling}
  	\dot{\beta}_t \leq e^{\alpha_t} \qquad \text{and} \qquad \dot{\gamma}_t = e^{\alpha_t}.
  \end{equation}
  If the ideal scaling conditions are satisfied, then Theorem 1.1 in~\citep{WiWiJo16} asserts that
  \begin{equation}
  	f(q(t))-f(q^*) \leq \mathcal{O}(e^{-\beta_t}),
  \end{equation}
along the trajectory $q(t)$ associated with the Bregman Lagrangian~\eqref{eq: Bregman L General Normed} and Bregman Hamiltonian~\eqref{eq: Bregman H General Normed}, where $q^*$ is the desired minimizer of the objective function $f$.

  From now on, we take $h(q) = \frac{1}{2} \langle q,q\rangle $. Assuming that the parameter functions $\alpha,\beta,\gamma$ satisfy the ideal scaling conditions \eqref{eq: IdealScaling}, the Bregman Lagrangian and Hamiltonian become
  \begin{align} 
  		L_{\alpha,\beta,\gamma}(q,v,t) & = \frac{1}{2}e^{ \gamma_t - \alpha_t} \langle v , v\rangle  - e^{\alpha_t + \beta_t +  \gamma_t} f(q) , \\ 
  		H_{\alpha,\beta,\gamma}(q,r,t) &= \frac{1}{2}e^{ \alpha_t -  \gamma_t } \langle r, r\rangle  +  e^{\alpha_t + \beta_t + \gamma_t} f(q) ,
  \end{align}
with corresponding Euler--Lagrange equation given by
  	\begin{equation} \label{eq: General Bregman EL}
  \ddot{q}(t)   +\left(  e^{\alpha_t} -  \dot{\alpha}_t \right) \dot{q}(t) + e^{2\alpha_t+\beta_t }  \nabla f(q(t)) = 0. 
  \end{equation}  

\hfill

\subsection{Polynomial Subfamily}  \label{subsec: Polynomial Subfamily} 

\hfill \\

A subfamily of Bregman dynamics of interest, indexed by a parameter $p>0$, is given by the choice of parameter functions
  \begin{equation} \label{eq: AlphaBetaGamma with p}
  	\alpha_t = \log{p} - \log{t} , \qquad  \beta_t = p \log{t} + \log{C},  \qquad \gamma_t = p \log{t},
  \end{equation}
  where $C>0$ is a constant. These parameter functions satisfy the ideal scaling conditions~\eqref{eq: IdealScaling}, and the corresponding Lagrangian and Hamiltonian are given by
     \begin{align} 
   	L_p(q,v,t) & = \frac{t^{ p+1}}{2p} \langle v , v\rangle  - Cpt^{2p-1} f(q) , \\ 
   	H_p(q,r,t) &= \frac{p}{2t^{p+1}} \langle r, r\rangle  +  Cpt^{2p-1} f(q) ,
   \end{align}
   with corresponding Euler--Lagrange equation given by
   \begin{equation} \label{eq: Polynomial EL}
   	\ddot{q}(t)   + \frac{p+1}{t} \dot{q}(t) + Cp^2 t^{p-2}  \nabla f(q(t)) = 0.
   \end{equation}

From Theorem 1.1 in~\citep{WiWiJo16}, the evolution $q(t)$ resulting from this dynamical system satisfies the convergence rate \begin{equation}  f(q(t))-f(q^*) \leq \mathcal{O}(1/t^p). \end{equation}

 Note that this Bregman subfamily has been exploited extensively in \cite{WiWiJo16,JordanSymplecticOptimization,duruisseaux2020adaptive,Duruisseaux2022Lagrangian}, and that the special case where $p=2$ and $C=1/4$ corresponds to  the limiting continuous differential equation introduced in \cite{SuBoCa16} for Nesterov's Accelerated Gradient method. \\

  \subsection{Exponential Subfamily}  \label{subsec: Exponential Subfamily} 
  
  \hfill \\
  
  Another subfamily of Bregman dynamics of interest, indexed by a parameter $\eta >0$, is given by the choice of parameter functions
  \begin{equation} \label{eq: AlphaBetaGamma with eta}
  	\alpha_t = \log \eta,   \qquad  \beta_t = \eta t + \log C, \qquad  \gamma_t = \eta t,
  \end{equation}
  where $C>0$ is a constant. These parameter functions satisfy the ideal scaling conditions~\eqref{eq: IdealScaling}, and the corresponding Lagrangian and Hamiltonian are given by
  \begin{align} 
  	L^{\eta}(q,v,t) & =  \frac{e^{\eta t}}{2\eta }  \langle v,v\rangle  - C \eta e^{2\eta t} f(q), \\ 
  	H^{\eta}(q,r,t) &=  \frac{\eta}{2e^{\eta t}}  \langle r,r\rangle + C\eta e^{2\eta t} f(q),
  \end{align}
  with corresponding Euler--Lagrange equation given by
  \begin{equation} \label{eq: Exponential EL}
  		\ddot{q}(t) + \eta \dot{q} + C \eta^2 e^{\eta t} \nabla f(q(t)) = 0
  \end{equation}
  
  From Theorem 1.1 in~\citep{WiWiJo16}, the evolution $q(t)$ resulting from this dynamical system satisfies the convergence rate \begin{equation}  f(q(t))-f(q^*) \leq \mathcal{O}\left(e^{-\eta t}\right). \end{equation}

\hfill  \\

\subsection{Geometric Numerical Integration of Time-rescaled Bregman dynamics}

\subsubsection{Time-rescaling Property of the Bregman Family} \label{subsubsec: Time-rescaling Property of the Bregman Family}

\hfill \\

  A very important property of the family of Bregman dynamics is its closure under time dilation:
\begin{theorem}[\citep{WiWiJo16}]\label{Theorem: Time Dilation Normed}
	If the curve $q(t)$ satisfies the Euler--Lagrange equations corresponding to the Bregman Lagrangian $L_{\alpha,\beta,\gamma}$, then the reparametrized curve $y(t) = q(\tau(t))$  satisfies the Euler--Lagrange equations corresponding to the Bregman Lagrangian $L_{\tilde{\alpha},\tilde{\beta},\tilde{\gamma}}$ where 
	\begin{equation}\label{eq: Time Dilation Transformation}
		L_{\tilde{\alpha},\tilde{\beta},\tilde{\gamma}} (q,v,t) = \dot{\tau} (t) L_{\alpha,\beta,\gamma}\left(q,\frac{1}{\dot{\tau} (t)}v,\tau(t)\right) , \quad \tilde{\alpha}_t = \alpha_{\tau(t)} + \log{\dot{\tau} (t)},  \quad \tilde{\beta}_t = \beta_{\tau(t)}, \quad \tilde{\gamma}_t = \gamma_{\tau(t)} .  \end{equation}	
\end{theorem} 

Thus, the entire subfamily of Bregman trajectories can be obtained by speeding up or slowing down along any specific Bregman curve in spacetime. It is natural to exploit this time-rescaling property with carefully chosen variable time-steps in the integrator to transform the time-dependent Bregman Hamiltonian or Lagrangian into a simpler autonomous system in some extended phase-space. This allows the higher-order Bregman dynamics to be integrated in a more computationally efficient fashion by time-rescaling the lower-order Bregman dynamics. This was first achieved in \cite{duruisseaux2020adaptive} with the polynomial subfamily of Section~\ref{subsec: Polynomial Subfamily}, where time-rescaling a solution to the $p$-Bregman Euler--Lagrange equations via $\tau(t) =t^{\mathring{p}/p}$ yielded a solution to the $\mathring{p}$-Bregman Euler--Lagrange equations. We can similarly jump from one solution of Bregman dynamics from the exponential subfamily from Section~\ref{subsec: Exponential Subfamily} to another via $\tau(t) =\frac{\mathring{\eta}}{\eta} t$, or jump from exponential Bregman dynamics to polynomial Bregman dynamics via $\tau(t) =\frac{p}{\eta} \log{t}$, and vice-versa via $\tau(t) = e^{\eta t / p}$.

However, when symplectic integrators were first used in combination with variable time-steps, they performed poorly \citep{CalSan93, GlaDunCan91}. A major advantage of using symplectic integrators on conservative Hamiltonian systems is that they exhibit excellent long-time near-energy preservation~\citep{HaLuWa2006}. Backward error analysis \citep{HaLuWa2006} shows that symplectic integrators can be associated with a modified Hamiltonian in the form of a formal power series in the time-step. Using variable time-steps results in a different modified Hamiltonian at every iteration, which is the source of the poor energy conservation and poor overall performance of these symplectic integrators. Fortunately, there are ways to circumvent this issue, which will allow us to exploit the time-rescaling property of the Bregman dynamics with variable time-step integrators, to transform the time-dependent Bregman dynamics into simpler autonomous systems in an extended space.  \\

\subsubsection{Time-adaptive Hamiltonian Integrators}  \label{subsubsec: Time-adaptive Hamiltonian Integrators}

\hfill \\

On the Hamiltonian side, the Poincar\'e transformation is a way to incorporate variable time-steps in geometric Hamiltonian integrators without losing their nice conservation properties \cite{Zare1975,Ha1997,duruisseaux2020adaptive}. Given a Hamiltonian $H(q,p,t)$, consider a desired transformation of time $t \mapsto \tau$ described by the monitor function $\frac{dt}{d\tau} = g(t).$ The time $t$ shall be referred to as the physical time of the system, while $\tau$ will be referred to as the fictive time. A new Hamiltonian system is constructed using the Poincar\'e transformation,
\begin{equation}
	\bar{H}(\bar{q},\bar{p}) = g(\mathfrak{q}) \left(H(q,\mathfrak{q},p) + \mathfrak{p} \right),
\end{equation}
in the extended phase space defined by 
\begin{equation} \bar{q} = \left[\begin{matrix} q \\ \mathfrak{q} \end{matrix} \right] \in \bar{\mathcal{Q}} \qquad \text{and} \qquad \bar{p} = \left[ \begin{matrix} p \\ \mathfrak{p} \end{matrix} \right], \end{equation} where $\mathfrak{p}$ is the conjugate momentum for $\mathfrak{q}=t$ with
$\mathfrak{p}(0)=-H(q(0),0,p(0))$. Then, using a symplectic integrator with constant time-step in fictive time $\tau$ on the Poincar\'e transformed Hamiltonian, has the effect of integrating the original system with the desired variable time-step in physical time $t$ via the relation $\frac{dt}{d\tau} = g(t)$. Note that this framework can be extended to monitor functions which also depend on position $q$ and momentum $p$, $g = g(q,t,p)$ (see \cite{duruisseaux2020adaptive}), but we will only need $g = g(t)$ in this paper. Also note that the Poincar\'e transformed Hamiltonian can be thought of as coming from a variational principle \cite{Duruisseaux2022Lagrangian}.  \\

Going back to accelerated optimization, and denoting momentum by $r$ to avoid confusion, we can jump from one form of Bregman dynamics to another as follows:
\begin{enumerate}
	\item Polynomial-$p$ to Polynomial-$\mathring{p}$: $\tau(t) =t^{\mathring{p}/p}$, $g(t) = \frac{p}{\mathring{p}} t^{1-\mathring{p}/p}$, yielding the Poincar\'e Hamiltonian 
	\begin{equation}
		\bar{H}_{p \rightarrow \mathring{p}}(\bar{q},\bar{r})  =   \frac{p^2}{2 \mathring{p} \mathfrak{q}^{p+\mathring{p}/p}} \langle r, r\rangle  +  \frac{Cp^2}{\mathring{p}} \mathfrak{q}^{2p-\mathring{p}/p } f(q) + \frac{p}{\mathring{p}} \mathfrak{r}  \mathfrak{q}^{1-\mathring{p}/p} .
	\end{equation}  

\item Exponential-$\eta$ to Exponential-$\mathring{\eta}$: $\tau(t) =\frac{\mathring{\eta}}{\eta} t$, $g(t) = \frac{\eta}{\mathring{\eta}}$, yielding the Poincar\'e Hamiltonian 
	\begin{equation}
	\bar{H}^{\eta \rightarrow \mathring{\eta}}(\bar{q},\bar{r})  =    \frac{\eta^2}{2 \mathring{\eta}e^{\eta \mathfrak{q}}}  \langle r,r\rangle + \frac{C\eta^2}{\mathring{\eta}} e^{2\eta \mathfrak{q}} f(q) + \frac{\eta}{\mathring{\eta}}  \mathfrak{r}.
\end{equation}  

\item Exponential-$\eta$ to Polynomial-$p$: $\tau(t) =\frac{p}{\eta} \log{t}$, $g(t) = \frac{\eta}{p} t$, yielding the Poincar\'e Hamiltonian 
	\begin{equation}
	\bar{H}^{\eta}_{ \rightarrow p}(\bar{q},\bar{r})  = \frac{\mathfrak{q} \eta^2}{2 p e^{\eta \mathfrak{q}}}  \langle r,r\rangle + \frac{C \mathfrak{q}\eta^2}{p} e^{2\eta \mathfrak{q}} f(q) +  \frac{\eta}{p} \mathfrak{q} \mathfrak{r}  .
\end{equation}  

\item Polynomial-$p$ to Exponential-$\eta$: $\tau(t) = e^{\eta t / p}$, $g(t) = \frac{p}{\eta} e^{-\eta t / p}$ yielding the Poincar\'e Hamiltonian
	\begin{equation}
	\bar{H}^{\rightarrow \eta}_{p}(\bar{q},\bar{r})  =    e^{-\frac{\eta }{p} \mathfrak{q}} \left( \frac{p^2 }{2 \eta \mathfrak{q}^{p+1}} \langle r, r\rangle  +  \frac{Cp^2}{\eta} \mathfrak{q}^{2p-1} f(q) + \frac{p}{\eta}\mathfrak{r} \right)  . 
\end{equation} %  \\
\end{enumerate}

\hfill 

\subsubsection{Time-adaptive Lagrangian Integrators}   \label{subsubsec: Time-adaptive Lagrangian Integrators}

\hfill \\

The time-adaptive framework for symplectic integration on the Hamiltonian side presented in the previous section relies on a degenerate Hamiltonian which has no associated Lagrangian description. Thus, we cannot exploit the usual correspondence between Hamiltonian and Lagrangian dynamics, and we follow a different strategy to allow time-adaptivity in Lagrangian integrators. Given a time-dependent Lagrangian $L(q,\dot{q},t)$, consider the extended autonomous Lagrangian
\begin{equation}   \bar{L} (\bar{q}(\tau),\bar{q}'(\tau) )  = \mathfrak{q}'(\tau) L\left(q(\tau), \frac{q'(\tau)}{g(\mathfrak{q}(\tau))}  , \mathfrak{q}(\tau) \right)  - \lambda(\tau) \left[ \mathfrak{q}'(\tau) - g(\mathfrak{q}(\tau)) \right] \label{eq: Extended Lagrangian}, \end{equation}
defined in the extended space $\bar{q} = (q,\mathfrak{q} , \lambda )^\top $ where time is viewed as a position coordinate $\mathfrak{q} = t$, where $\lambda$ is a Lagrange multiplier enforcing the desired time rescaling $	\frac{dt}{d\tau} = g(t)$, and where apostrophes denote derivatives with respect to fictive time~$\tau$. Now, if $\left(\bar{q}(\tau), \bar{q}'(\tau) \right)$ satisfies the Euler--Lagrange equations corresponding to the Lagrangian $\bar{L}$, then its components satisfy $\frac{dt}{d\tau} = g(t)$ and the original Euler--Lagrange equations \cite{Duruisseaux2022Lagrangian}. %\\

A discrete variational formulation of these continuous extended Lagrangian mechanics can be formulated \cite{Duruisseaux2022Lagrangian}, by considering a discrete Lagrangian
\begin{equation}
	L_d(q_k, \mathfrak{q}_k,  q_{k+1} , \mathfrak{q}_{k+1})  \approx  \  \ext_{ \substack{ (q,\mathfrak{q} ) \in C^2 ([\tau_{k} , \tau_{k+1} ] , \mathcal{Q} \times \mathbb{R} ) \\ (q, \mathfrak{q})(\tau_k) = (q_k, \mathfrak{q}_k)  , \text{  } (q, \mathfrak{q})(\tau_{k+1}) = (q_{k+1}, \mathfrak{q}_{k+1}) }  }   \text{  } \int_{\tau_k}^{\tau_{k+1}}{L\left(q, \frac{q'}{g(\mathfrak{q})}  , \mathfrak{q} \right) d\tau },
\end{equation}
where $0 = \tau_0 < \tau_1 < \ldots < \tau_N$ partitions the time interval of interest, and $\{  (q_k, \mathfrak{q}_k) \}_{k=0}^{N} $ is a discrete curve in $\mathcal{Q} \times  \mathbb{R}$ such that $q_k \approx q(\tau_k)$ and $\mathfrak{q}_k \approx \mathfrak{q}(\tau_k)$. Defining the discrete momenta via the discrete Legendre transformations, $p_k = - D_1 L_d(q_k,\mathfrak{q}_k,q_{k+1},\mathfrak{q}_{k+1})$, and using a constant time-step $h$ in fictive time $\tau$, the corresponding discrete extended Euler--Lagrange equations can be written as 
\begin{equation} \label{eq: Discrete Extended EL}\begin{aligned}    p_k  & = - D_1 L_d(q_k,\mathfrak{q}_k,q_{k+1},\mathfrak{q}_{k+1}) ,   \\
	p_{k+1}    & = \frac{g(\mathfrak{q}_k)}{g(\mathfrak{q}_{k+1})}  D_3 L_d(q_k, \mathfrak{q}_k,  q_{k+1} , \mathfrak{q}_{k+1}) , \\ \mathfrak{q}_{k+1} & =  \mathfrak{q}_k + hg(\mathfrak{q}_k) ,  \end{aligned}
\end{equation}
with two additional equations for $\mathfrak{p}_k$ and $\mathfrak{p}_{k+1}$ (see \cite{Duruisseaux2022Lagrangian}). For accelerated optimization, we are not interested in the evolution of $\mathfrak{p}$, and since it does not appear in the updates for the other variables we do not need these equations and omit them here. We can then use one of the monitor functions \begin{equation}
	g(t) = \frac{p}{\mathring{p}} t^{1-\mathring{p}/p}, \qquad g(t) = \frac{\eta}{\mathring{\eta}}, \qquad g(t) = \frac{\eta}{p} t, \qquad g(t) = \frac{p}{\eta} e^{-\eta t / p},
\end{equation} to transform from one type of Bregman Lagrangian to another. 

\hfill

\section{Numerical Methods and Problems of Interest}  \label{section: Numerical Methods and Problems of Interest}

\subsection{Numerical Methods} \label{section: Numerical Methods}

\hfill \\

We now present four different methods to design symplectic integrators for the time-rescaled Bregman Lagrangian and Bregman Hamiltonian systems presented in Section \ref{section: Variational Accelerated Optimization}. Keeping in mind the desired applications in machine learning where problem sizes and data sets are very large, we restrict ourselves to explicit first-order optimization algorithms. Each of these four methods can be used within the four different adaptive approaches presented in Section \ref{subsubsec: Time-rescaling Property of the Bregman Family} (polynomial, exponential, polynomial-to-exponential, and exponential-to-polynomial), and the resulting sixteen algorithms are presented in Appendix~\ref{appendix: List of Time-Adaptive Algorithms}. \\

\subsubsection{\textbf{H}amiltonian \textbf{T}aylor \textbf{V}ariational \textbf{I}ntegrator (\textbf{HTVI})} 

\hfill  \\

 Proceeding as in~\cite[Section 4.4]{duruisseaux2020adaptive} or~\cite{ScShLe2017}, we can derive the Hamiltonian Taylor Variational Integrator (HTVI),
\begin{equation}
		\begin{aligned}
		p_{k+1}&=p_k- hH_q(q_k,p_{k+1}), \\ q_{k+1}&=q_k+h H_p(q_k,p_{k+1}).  
	\end{aligned} 
\end{equation}
These updates recover the Symplectic Euler method \cite{HaLuWa2006}, which is a popular symplectic integrator of order 1.  \\

\subsubsection{\textbf{L}agrangian \textbf{T}aylor \textbf{V}ariational \textbf{I}ntegrator (\textbf{LTVI})}  

\hfill \\

As in \cite{Duruisseaux2022Lagrangian}, we can define a discrete Lagrangian,
\begin{align}
	L_d(\bar{q}_0,\bar{q}_1)  & = h L_{p} \left(q_0 , \frac{q_1 - q_0}{ hg(\mathfrak{q}_0)} , \mathfrak{q}_0 \right) ,
\end{align} 
and the updates for the Lagrangian Taylor Variational Integrator (LTVI) can be obtained from the discrete extended Euler--Lagrange equations \eqref{eq: Discrete Extended EL}. \\

\subsubsection{\textbf{S}t\"ormer-\textbf{V}erlet (\textbf{SV})}  

\hfill \\

A popular symplectic integrator is the St\"ormer--Verlet (SV) method,
%	\begin{equation}
	%		\begin{array}{cc}
		%			\textbf{St\"ormer-Verlet (SV)}: \quad 
		%			\begin{aligned}
			%				p_{n+1/2}&=p_n-\frac{h}{2}H_q(p_{n+1/2},q_{n}) \\ q_{n+1}&=q_n+\frac{h}{2}\left(H_p(p_{n+1/2},q_{n})+H_p(p_{n+1/2},q_{n+1})\right)   \\
			%				p_{n+1}&=p_{n+1/2}-\frac{h}{2}H_q(p_{n+1/2},q_{n+1})
			%		\end{aligned} 	\end{array}  
	%	\end{equation}
\begin{equation}
	\begin{aligned}
		p_{k+1/2}&=p_k-\frac{h}{2}H_q(q_k,p_{k+1/2}) ,\\ q_{k+1}&=q_k+\frac{h}{2}\left[H_p(q_k,p_{k+1/2})+H_p(q_{k+1},p_{k+1/2})\right] ,  \\
		p_{k+1}&=p_{k+1/2}-\frac{h}{2}H_q(q_{k+1},p_{k+1/2}),
	\end{aligned} 
\end{equation}
which is a symmetric symplectic integrator of order 2 (see \cite{HaLuWa2006}). A very detailed description of the St\"ormer--Verlet method, its different interpretations, and its beneficial numerical properties can be found in~\cite{Hairer2003}. Note however that in the polynomial and polynomial-to-exponential frameworks, the update for~$\mathfrak{q}$ in the resulting integrators becomes implicit (see PolySV and PolyToExpoSV in Appendix \ref{appendix: List of Time-Adaptive Algorithms}), which makes these integrators less desirable. For the accelerated optimization application, we will usually be able to combine the first and last updates for the momentum vector $p$ into a single update and save roughly a third of the computational time. This is because St\"ormer--Verlet is conjugate to symplectic Euler. \\

\subsubsection{\textbf{S}ymmetric \textbf{L}eapfrog \textbf{C}omposition of Component Dynamics (\textbf{SLC})}

\hfill \\

 The main idea is to decompose the vector field into its components,
 \begin{align*}
 	\frac{d}{d\tau} & = \frac{dq}{d\tau} \frac{d}{dq}  + \frac{d\mathfrak{q}}{d\tau} \frac{d}{d\mathfrak{q}}  + \frac{dr}{d\tau} \frac{d}{dr}  + \frac{d\mathfrak{r}}{d\tau} \frac{d}{d\mathfrak{r}}   = \frac{\partial H}{\partial r} \frac{d}{dq} + \frac{\partial H}{\partial \mathfrak{r}}  \frac{d}{d\mathfrak{q}}  - \frac{\partial H}{\partial q} \frac{d}{dr}  - \frac{\partial H}{\partial \mathfrak{q}}  \frac{d}{d\mathfrak{r}}  = \mathcal{A} +\mathcal{B}  + \mathcal{C}   + \mathcal{D} ,  \end{align*}
and then combine the component dynamics using a symmetric leapfrog composition
\begin{equation}
	\Phi_h = \exp{\left(\frac{h}{2} \mathcal{D} \right)}  \circ   \exp{\left(\frac{h}{2} \mathcal{C} \right)}  \circ   \exp{\left(\frac{h}{2} \mathcal{B} \right)} \circ  \exp{\left(h\mathcal{A} \right)}  \circ   \exp{\left(\frac{h}{2} \mathcal{B} \right)} \circ   \exp{\left(\frac{h}{2} \mathcal{C} \right)}  \circ \exp{\left(\frac{h}{2} \mathcal{D} \right)}  
\end{equation}
which satisfies $\Phi_h = \exp{\left( h H \right)} + \mathcal{O}(h^3)$ (can be shown using the Baker--Campbell--Hausdorff formula). This strategy is similar to the integrator from \cite[Section 3.3]{JordanSymplecticOptimization} and the Splitting algorithms in \cite{duruisseaux2020adaptive}.

As an example, for 
 	\begin{equation}
 	 \bar{H}_{p \rightarrow \mathring{p}}(\bar{q},\bar{r})  =   \frac{p^2}{2 \mathring{p} \mathfrak{q}^{p+\mathring{p}/p}} \langle r, r\rangle  +  \frac{Cp^2}{\mathring{p}} \mathfrak{q}^{2p-\mathring{p}/p } f(q) + \frac{p}{\mathring{p}} \mathfrak{r}  \mathfrak{q}^{1-\mathring{p}/p},
 \end{equation}  
the components of the vector field are given by
 \begin{equation}
 \mathcal{A} =	\frac{p^2}{\mathring{p} \mathfrak{q}^{p+\mathring{p}/p}}  r  \frac{d}{dq} ,    \quad  \mathcal{B}= \frac{p}{\mathring{p}} \mathfrak{q}^{1-\mathring{p}/p}\frac{d}{d\mathfrak{q}}   , \quad  \mathcal{C} = - \frac{Cp^2}{\mathring{p}} \mathfrak{q}^{2p-\mathring{p}/p } \nabla f(q) \frac{d}{dr} ,  \quad  \mathcal{D} = - \frac{\partial \bar{H}_{p \rightarrow \mathring{p}}}{\partial \mathfrak{q}} \frac{d}{d\mathfrak{r}} .
 \end{equation}
Then, the corresponding component dynamics are given as follows:
\begin{itemize}
	\item $\exp{\left(h \mathcal{A} \right)} $ yields the update $ q \leftarrow  q + h \frac{p^2}{\mathring{p} \mathfrak{q}^{p+\mathring{p}/p}}  r  $ \\
	\item  $\exp{\left(h \mathcal{C} \right)} $ yields the update $ r \leftarrow  r - h \frac{Cp^2}{\mathring{p}} \mathfrak{q}^{2p-\mathring{p}/p } \nabla f(q) $ \\
	\item  $\exp{\left(h \mathcal{B} \right)} $ yields the differential equation $\mathfrak{q}' = \frac{p}{\mathring{p}} \mathfrak{q}^{1-\mathring{p}/p}$, which can be solved exactly to obtain the update $\mathfrak{q} \leftarrow  \left( \mathfrak{q}^{\mathring{p}/p} + h \right)^{p/\mathring{p}} $ 
\end{itemize} 
Note that the updates for  $\exp{\left(h \mathcal{A} \right)} $,  $\exp{\left(h \mathcal{B} \right)} $, and  $\exp{\left(h \mathcal{C} \right)} $ do not involve $\mathfrak{r}$, and in practice, we are not interested in the evolution of $\mathfrak{r}$, so we can simplify the composition into
 \begin{equation} \label{eq: SLC Composition}
	\Phi_h = \exp{\left(\frac{h}{2} \mathcal{C} \right)}   \circ   \exp{\left(\frac{h}{2} \mathcal{B} \right)} \circ  \exp{\left(h\mathcal{A} \right)}  \circ   \exp{\left(\frac{h}{2} \mathcal{B} \right)} \circ  \exp{\left(\frac{h}{2} \mathcal{C} \right)}  ,
\end{equation}
which gives the PolySLC algorithm,
\begin{equation}
	\begin{aligned}
					 r & \leftarrow  r - \frac{Cp^2}{2\mathring{p}} h \mathfrak{q}^{2p-\mathring{p}/p } \nabla f(q) , \\
		\mathfrak{q} & \leftarrow  \left( \mathfrak{q}^{\mathring{p}/p} + \frac{h}{2} \right)^{p/\mathring{p}}, \\
			q &\leftarrow  q + \frac{hp^2}{\mathring{p} \mathfrak{q}^{p+\mathring{p}/p}}  r ,\\
			\mathfrak{q} & \leftarrow  \left( \mathfrak{q}^{\mathring{p}/p} + \frac{h}{2} \right)^{p/\mathring{p}} , \\
					 r & \leftarrow  r - \frac{Cp^2}{2\mathring{p}} h \mathfrak{q}^{2p-\mathring{p}/p } \nabla f(q) .
	\end{aligned} 
\end{equation} 

We have chosen to place $\exp{\left(h \mathcal{A} \right)} $ in the middle of the symmetric composition~\eqref{eq: SLC Composition} so that the gradient $\nabla f$ only needs to be evaluated at the iterates $\{q_k\}_{k\in \mathbb{N}}$ and can also be used without further computations in stopping criteria or momentum restarting schemes. We have also chosen to place $\exp{\left(h \mathcal{C} \right)} $ at the left-hand and right-hand of the symmetric composition~\eqref{eq: SLC Composition} so that in practice, we may combine the first and last updates for the vector $r$ into a single update (instead of only being able to combine the first and last updates for the scalar $\mathfrak{q}$ into a single update), and thus save roughly a third of the computational time. \\

\textbf{Remarks.} It was observed in \cite{duruisseaux2020adaptive} that the symplecticity of the integrator was essential for the efficient, robust, and stable discretization of these variational flows describing accelerated optimization. Therefore, we will not consider non-symplectic methods here. %\\ 

Higher-order explicit symplectic integrators can be derived, leveraging higher-order compositions such as Yoshida splittings \cite{Yoshida1990}, but it was observed in \cite{duruisseaux2020adaptive} that these require much more evaluations of the objective function and its gradient at each step. Thus, the resulting algorithms would not be competitive in terms of computational time and number of gradient evaluations, since the other methods usually converge in a similar number of iterations but only require one gradient evaluation per iteration.  

%\hfill  \\

\subsection{Problems of Interest}

\hfill \\

	A subset $C$ of $\mathbb{R}^d$ is \textbf{convex} if  $\lambda x + (1-\lambda) y \in C$ for any $x,y\in C$ and $\lambda \in [0,1]$. A differentiable function $f:\mathbb{R}^d \rightarrow \mathbb{R}$ is \textbf{convex} if its domain $\text{dom}(f)$ is convex and for any $\lambda \in [0,1]$,  \begin{equation} f\left(\lambda x +(1-\lambda)y\right) \leq \lambda f(x) + (1-\lambda) f(y), \qquad  \forall x,y\in \text{dom}(f),\end{equation}
	or equivalently if \begin{equation} f(y) \geq f(x) + \nabla f(x)^\top (y-x), \qquad  \forall x,y\in \text{dom}(f). \end{equation}  A differentiable function $f:\mathbb{R}^d \rightarrow \mathbb{R}$ is \textbf{strongly convex} if there exists $\mu >0$ such that $f(x) - \mu  \| x\|^2$ is convex, or equivalently if 
	\begin{equation}  f(y) \geq f(x) + \nabla f(x)^\top (y-x) + \mu \| y-x \|^2 ,\qquad  \forall x,y\in \text{dom}(f). \end{equation} 
	
	\hfill 
	
	In our numerical experiments, we will use termination criteria of the form,
	\begin{equation}\label{eq: TerminationCriterion}
		|f(x_k) - f(x_{k-1})|<\delta  \quad \text{and} \quad \Vert  \nabla f(x_k) \Vert < \delta,
	\end{equation}
	for various values of the tolerance $\delta$, and solve the following convex problems. %\\

\begin{problem} \label{Problem: Quartic}
	Minimize the quartic polynomial \begin{equation} f(x) = 1 +  \left[(x-1)^\top  \Sigma (x-1) \right]^2, \quad  \text{where }\Sigma_{ij} =0.9^{|i-j|} \text{ and } x\in \mathbb{R}^{d}. \end{equation} 
%	\begin{equation}  \nabla f(x) = 4 (x-1)^\top  \Sigma (x-1) \Sigma (x-1), \end{equation} 
This convex function achieves its global minimum at $x^*= (1,1,\ldots,1)^\top $.  \\
\end{problem}

\begin{problem}\label{Problem: Polynomial and Log}
	Minimize the convex (not strongly convex) function \begin{equation} f(x_1,x_2) = x_1 +x_2^2 - \ln(x_1x_2) ,\end{equation} which achieves its global minimum at  $x^* = (1,\sqrt{2}/2)^\top $.   \\
\end{problem}

\begin{problem}\label{Problem: x log(x)}
	Minimize the strongly convex function \begin{equation} f(x_1,...,x_d) = \sum_{k=1}^{d}{x_k \log x_k}. \end{equation}
	This function, known as negative entropy, achieves its global minimum at $x ^*= (e^{-1} ,e^{-1}  , \ldots, e^{-1}) ^\top $.   \\
\end{problem}

\begin{problem} \label{Problem: Ill-Conditioned Quadratic}
	Minimize the ill-conditioned strongly convex function  \begin{equation} f(x_1,x_2,x_3) = 1 + 0.01x_1^2 + x_2^2 + 100x_3^2, \end{equation} which achieves its global minimum at $x ^*= (0,0,0)^\top$. 
\end{problem}

\hfill 

\begin{problem} [\textbf{Linear Regression} or \textbf{Least Squares}] \label{Problem: Linear Regression}
	Given a matrix $A \in \mathbb{R}^{m\times n}$ with $m \geq n$ and a vector $b \in \mathbb{R}^m$, consider the problem of finding a vector $x \in \mathbb{R}^n$ such that $\| Ax-b \|_2$ is minimized. The least squares problem has many applications in data-fitting and interpolation. It can be formulated as the minimization of
	\begin{equation}
		f(x) =\frac{1}{2} x^\top A^\top Ax - b^\top Ax ,
	\end{equation}
	with a gradient given by $\nabla f(x) = A^\top Ax -  A^\top b$. 
	A vector $x\in \mathbb{R}^n$ is a solution of the least squares problem if and only if it satisfies the normal equation $A^\top A x = A^\top b$. Furthermore, the least squares problem has a unique solution, given by $x^* = (A^\top A)^{-1} A^\top b$, if and only if $A$ has full rank \cite{Trefethen1997}. 
	
	There are also regularized versions of the least squares problem or linear regression \cite{Boyd2004,Bertsekas2009}, to penalize larger values of~$x$. A common form of regularization is Tikhonov regularization \cite{Tikhonov1963,Tikhonov1977,Philips1962} (or $\ell^2$ regularization), where we minimize the convex function
	\begin{equation}
		f(x) =   \| Ax-b\|_2^2 + \lambda \| x\|_2^2,
	\end{equation}
	for some $\lambda >0$, which has a unique minimizer $x^* = (A^\top A + \lambda b)^{-1} A^\top b$. Another regularized version is the $\ell^1$ penalized linear regression (also known as the Lasso problem \cite{Tibshirani1996}), where we minimize the convex (not strongly convex) function
	\begin{equation}
		f(x) =  \frac{1}{2} \| Ax-b\|_2^2 + \lambda \| x\|_1.
	\end{equation}

\end{problem}

\hfill 

\begin{problem} [\textbf{Logistic Regression} for Binary Classification]  \label{Problem: Binary Classification}
	Given a set of feature vectors $x_1, \ldots , x_m \in \mathbb{R}^n$ and associated labels $y_1, \ldots , y_m \in \{-1,1\}$, we want to find a vector $w \in \mathbb{R}^n $ such that $\emph{sign} ( w^\top x )$ is a good model for $y(x)$. This can be formulated as the problem of minimizing the convex (not strongly convex) function
	\begin{equation} \label{eq: Binary Classification Loss}
		f(w)  = \sum_{i=1}^{m}{\log\left( {1+\exp{\left(- y_i w^\top x_i   \right)}} \right) }.
	\end{equation}
	As for linear regression, there are also regularized versions of logistic regression, such as $\ell^1$ and $\ell^2$ regularized logistic regression:
	\begin{align}
		f(w)   = \sum_{i=1}^{m}{\log\left( {1+\exp{\left(- y_i w^\top x_i   \right)}} \right) } + \lambda \| x\|_1,  \quad  \text{ }
		f(w)   = \sum_{i=1}^{m}{\log\left( {1+\exp{\left(- y_i w^\top x_i   \right)}} \right) } +  \lambda \| x\|_2^2.
	\end{align}
\end{problem}

\hfill 

\begin{problem}[\textbf{Fermat--Weber Location Problem} \cite{Drezner2002,Beck2009,Boltyanski1999}]  \label{Problem: Location}

	Given a set of points $y_1, \ldots ,y_m  \in \mathbb{R}^n$ and associated positive weights $w_1, \ldots , w_m \in \mathbb{R}$, we want to find the location $x\in \mathbb{R}^n$ whose sum of weighted distances from the points $y_1, \ldots ,y_m$ is minimized. In other words, we wish to minimize the convex function
	\begin{equation} \label{eq: Location Loss}
		f(x) = \sum_{j=1}^{m}{w_j \|  x - y_j\|}.
	\end{equation}
The Fermat--Weber location problem is at the heart of Location Theory and has countless applications across many fields of science and engineering.
\end{problem}

\hfill

\begin{remark} A Tikhonov-type regularization can also be achieved by modifying the second-order differential equation instead of adding a penalty to the objective function (see \cite{Attouch2016,Attouch2017,Alecsa2021,Jendoubi2010} for instance). The idea is to add an extra term $\epsilon(t) x(t)$ with $\epsilon(t) \rightarrow 0$ as $t \rightarrow \infty$ to the second-order differential equation of interest:
	\begin{equation}
		\ddot{x}(t)+	\alpha(t) \dot{x}(t)   + \gamma(t) \nabla f (x(t)) + \epsilon(t) x(t) = 0.
	\end{equation} 
This extra term forces the generated trajectory to converge
to a solution of minimal norm. This type of modified differential equation can be generated from a variational framework via Lagrangians and Hamiltonians of the form
	\begin{equation}
		L(x,v,t) = \frac{1}{2} \alpha(t) \langle v , v \rangle  + \epsilon(t) \langle v,x\rangle - \gamma(t) f(x) , 
	\end{equation}
	\begin{equation}
	H(x,p,t) = \frac{1}{2\alpha(t)}  \langle p - \epsilon(t) x  , p - \epsilon(t) x \rangle + \gamma(t) f(x) ,
\end{equation}
	whose Euler--Lagrange equation reads
	\begin{equation}
		\alpha(t)	\ddot{x}(t) +	\dot{\alpha}(t)  \dot{x}(t)   + \gamma(t) \nabla f (x(t)) + \dot{\epsilon}(t) x(t) = 0.
	\end{equation} 

\end{remark}

%\newpage 

\section{Controlling the Oscillatory Behavior}   \label{section: Controlling the Oscillatory Behavior}

\hfill 

The Bregman Euler--Lagrange equation \eqref{eq: General Bregman EL} can be written in the form
\begin{equation}
	\ddot{x}(t)  + d(t) \dot{x}(t) + b(t) \nabla f(x(t)) = 0.
	\end{equation}
The introduction of momentum in the dynamical system causes the solution to this ordinary differential equation to overshoot frequently in its path towards the minimizer of the objective function~$f$, and as a result the solution can be highly oscillatory. Therefore, this differential equation can be thought of as modeling a nonlinear oscillator with damping, and the convergence of the function $f$ to its minimum value is not monotone along Bregman trajectories. This is similar to what was observed for the limiting continuous differential equation for Nesterov's accelerated gradient method~\cite{SuBoCa16,Muehlebach2019} and for most momentum methods. These oscillations are problematic since they can significantly slow down optimization algorithms that are derived from the discretization of these Bregman differential equations. Indeed, to resolve the fast oscillations of the differential equation, the time-step in the discretization has to be reduced sufficiently, which can considerably increase the number of iterations and gradient evaluations needed to achieve convergence. If the time-step is not taken small enough, the momentum in the algorithms can lead to large overshoots which can result in divergence. It would therefore be desirable to have a mechanism to neutralize these oscillations. Fortunately, there are ways to reduce the effect of these oscillations, which we will discuss in the remainder of this section. \\

\subsection{Momentum Restarting}

\hfill \\

Momentum causes the solution to the Bregman Euler--Lagrange equation to overshoot frequently on its path towards the minimizer of $f$. One strategy to control these overshoots and reduce the effect of the resulting oscillations is to use restarting or momentum restarting schemes, previously explored in \cite{Powell1977,Dai2004,O'donoghue2015,Giselsson2015,SuBoCa16,Fercoq2016,Donghwan2018,Fercoq2019,Roulet2020,Renegar2022}. We will consider three different momentum restarting schemes:  \\

\begin{itemize}
	\item \textbf{Function Scheme}: Restart momentum whenever $f(q_k) > f(q_{k-1})$ \\ This scheme restarts the momentum $r$ whenever the function evaluation at the new update moves away from the minimum value, to try to avoid wasting iterations in a bad direction. \\

	\item \textbf{Gradient Scheme}: Restart momentum whenever $\nabla f(q_k) (q_k - q_{k-1}) > 0$ \\ 	This restarts the momentum variable $r$ whenever momentum seems to take the new updates in a bad direction, measured using the gradient at that point. \\

	\item \textbf{Velocity Scheme}: Restart momentum whenever $\| q_{k+1}  - q_k \| < \|  q_k - q_{k-1}  \|$ \\
	This scheme restarts the momentum variable $r$ whenever the norm of the (discrete version of the) velocity~$\| \dot{q} \|$ starts decreasing, to try to maintain a high velocity along the trajectory.
	
\end{itemize}

\hfill  

Note that the quantities needed to implement these restarting schemes are already calculated in the standard versions of the optimization algorithms, and thus there is a negligible difference in the computational costs of each iteration in the restarted and non-restarted schemes. 

We can also require a minimum number of iterations between momentum restarts, to avoid having consecutive restarts that are too close to each other. In practice however, it did not seem to really improve the computational efficiency of the algorithm and could sometimes negatively impact the overall performance. For simplicity, we will not impose a minimum number of iterations between consecutive restarts.

\begin{figure}[!ht] 
	\centering
	\includegraphics[width=1\textwidth]{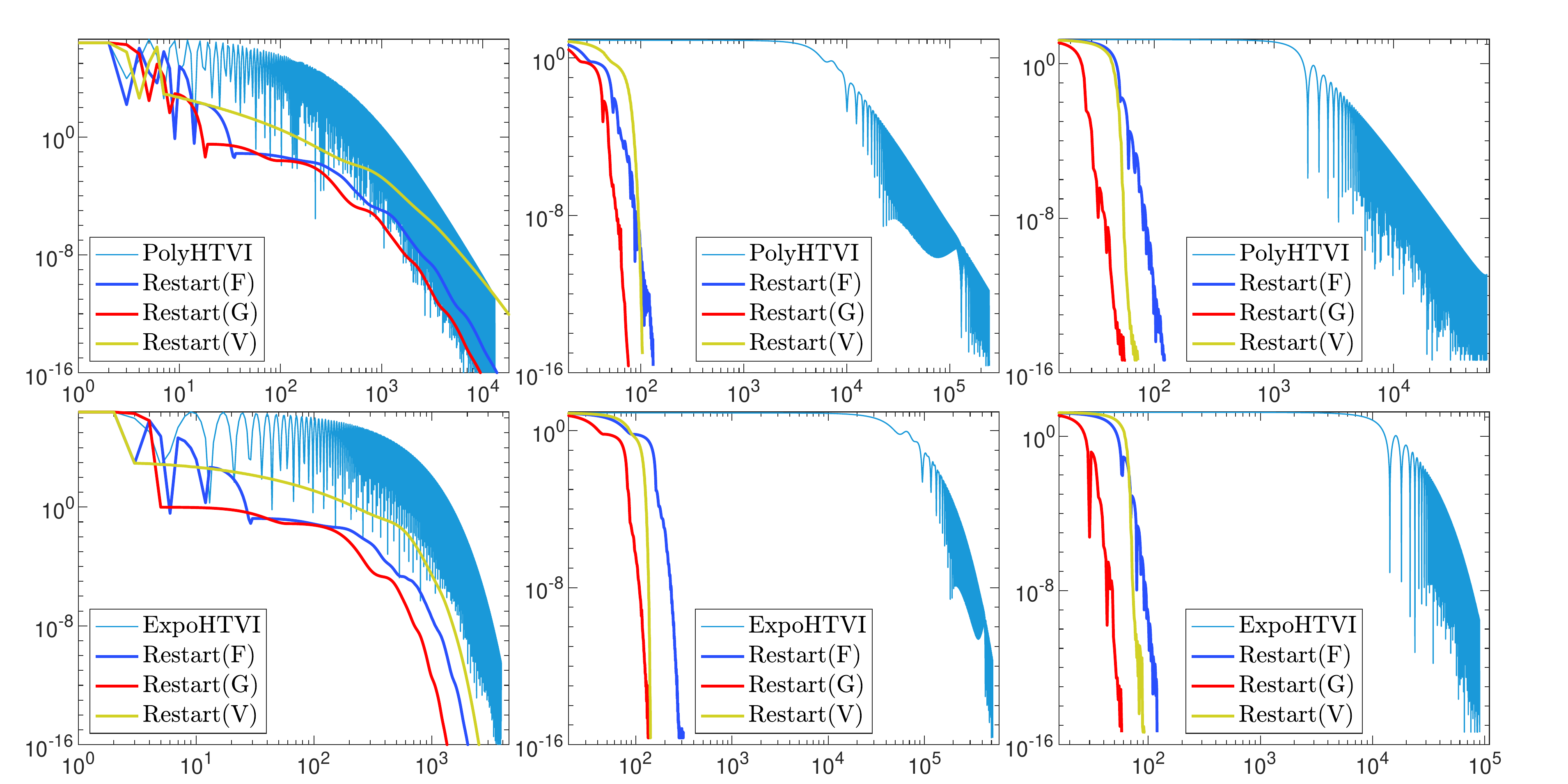} 	\vspace*{-7mm} \caption{ Error vs. Iterations number as the standard PolyHTVI and ExpoHTVI algorithms and their restarted versions (Function~(F), Gradient~(G) and Velocity~(V)) are applied to Problem~\ref{Problem: Quartic} (left), Problem~\ref{Problem: Polynomial and Log} (middle), and Problem~\ref{Problem: x log(x)} (right).}\label{fig: Restarting} 
\end{figure}

In our first numerical experiment, we compared the performance of the standard algorithms to their restarted versions on three different problems for fixed values of all the parameters except the time-step $h$. Figure~\ref{fig: Restarting} shows the  resulting error plots after tuning the value of $h$ optimally. We can clearly see that the restarted versions of the algorithms are much less oscillatory, and they can allow for much larger time-steps leading to significantly faster algorithms, as is the case for Problems~\ref{Problem: Polynomial and Log} and \ref{Problem: x log(x)}. Problem~\ref{Problem: Quartic} is a special instance where larger time-steps cannot be taken in the restarted algorithms, despite their non-oscillatory nature. It should be noted however that although the use of momentum restarting does not lead to significant improvements in computational efficiency, it does not penalize computational efficiency either.

\begin{figure}[!h] 
	\centering
	\includegraphics[width=1\textwidth]{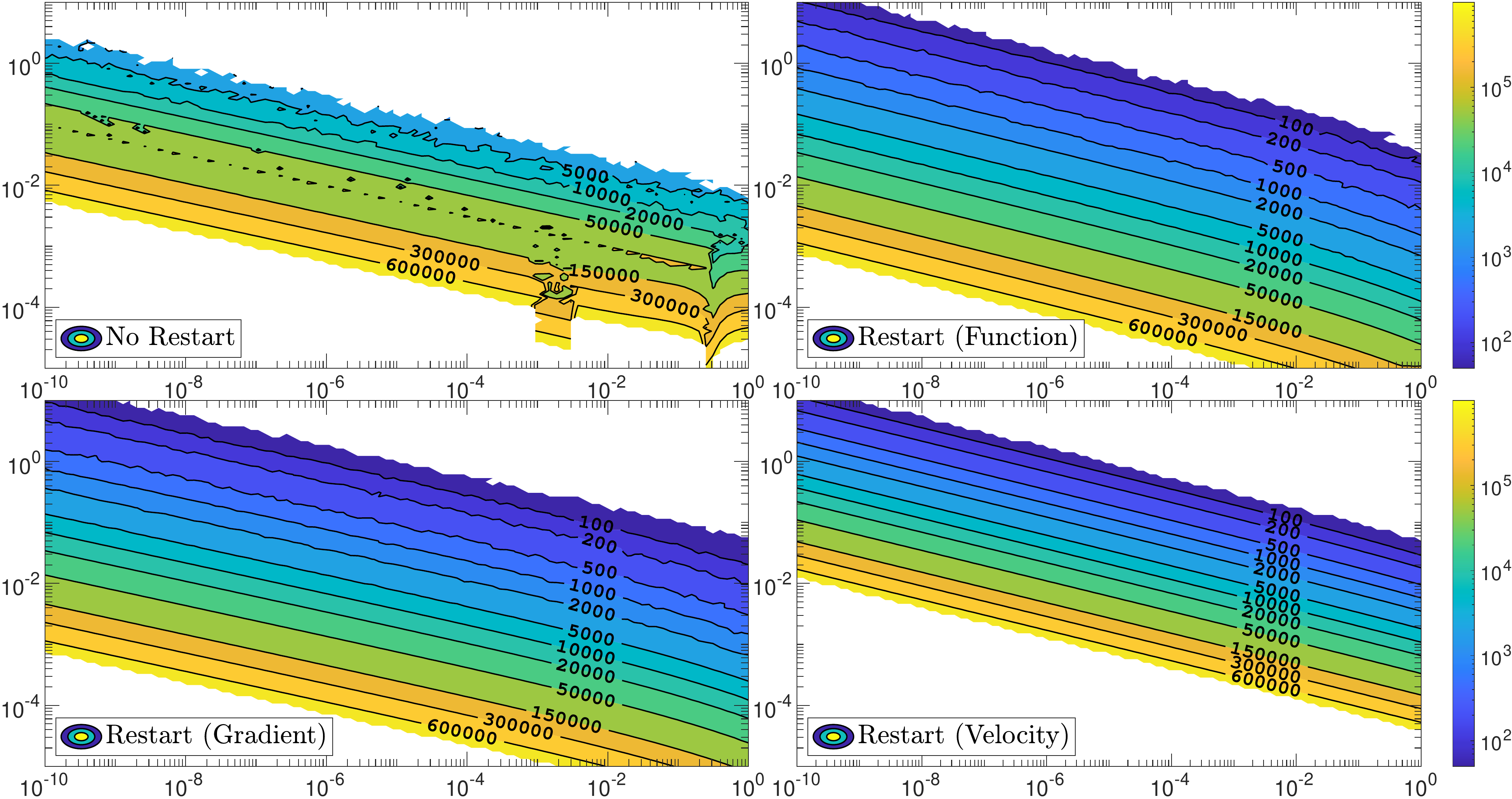} 	\vspace*{-7mm} \caption{ Contour plot of the number of iterations required to achieve convergence with $\delta = 10^{-5}$ in the $(C,h)$-plane, for $p = \mathring{p} = 4$ PolyHTVI applied to Problem~\ref{Problem: Polynomial and Log}. } \label{fig: RestartingRobustnessPoly1}  \vspace*{-1mm}
\end{figure}

\begin{figure}[!h] 
	\centering
	\includegraphics[width=1\textwidth]{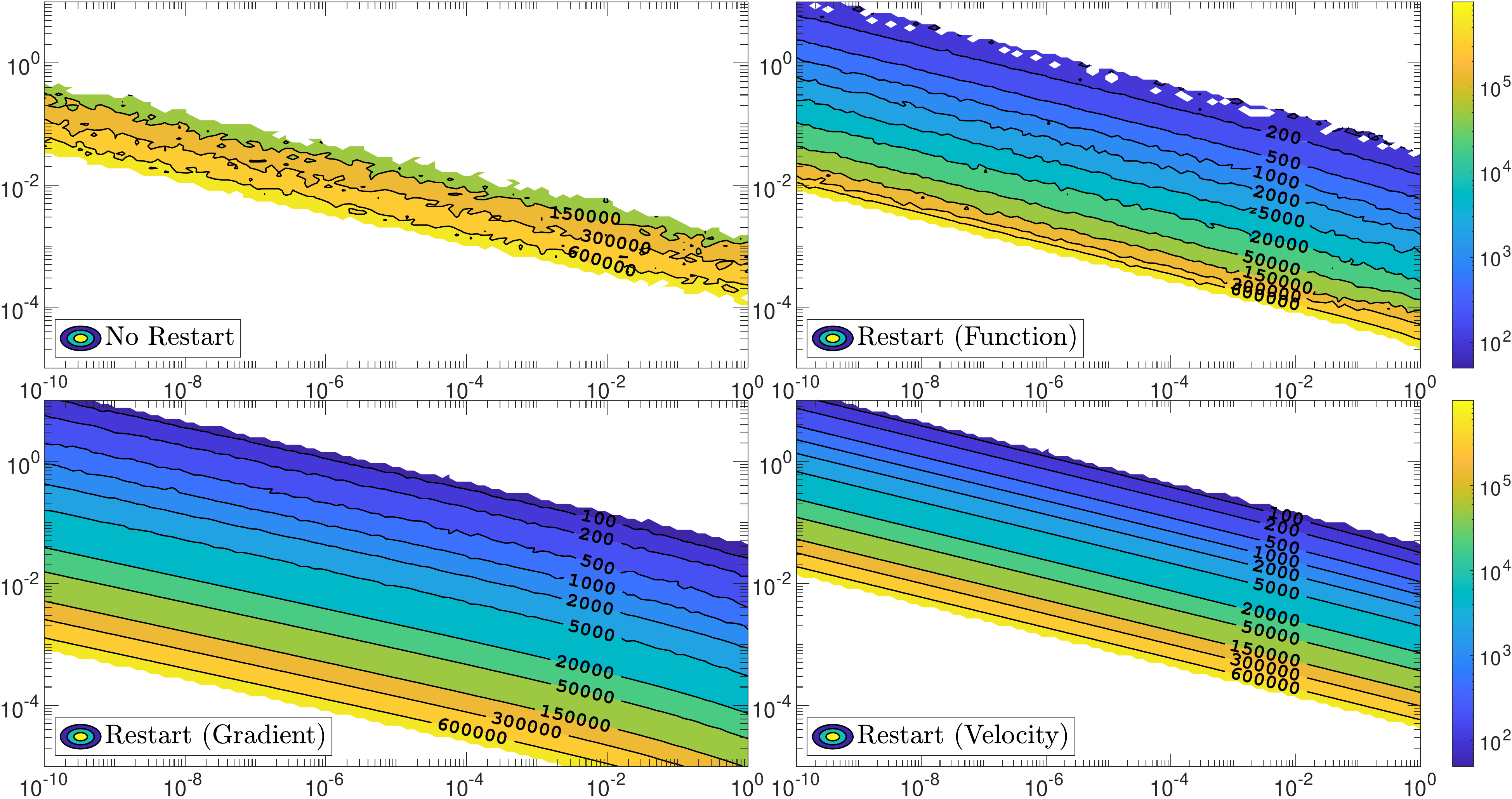} 	\vspace*{-7mm} \caption{ Contour plot of the number of iterations required to achieve convergence with $\delta = 10^{-8}$ in the $(C,h)$-plane, for $p = \mathring{p} = 4$ PolyHTVI applied to Problem~\ref{Problem: Polynomial and Log}. } \label{fig: RestartingRobustnessPoly2} 
\end{figure}

\begin{figure}[!h] 
	\centering
	\includegraphics[width=1\textwidth]{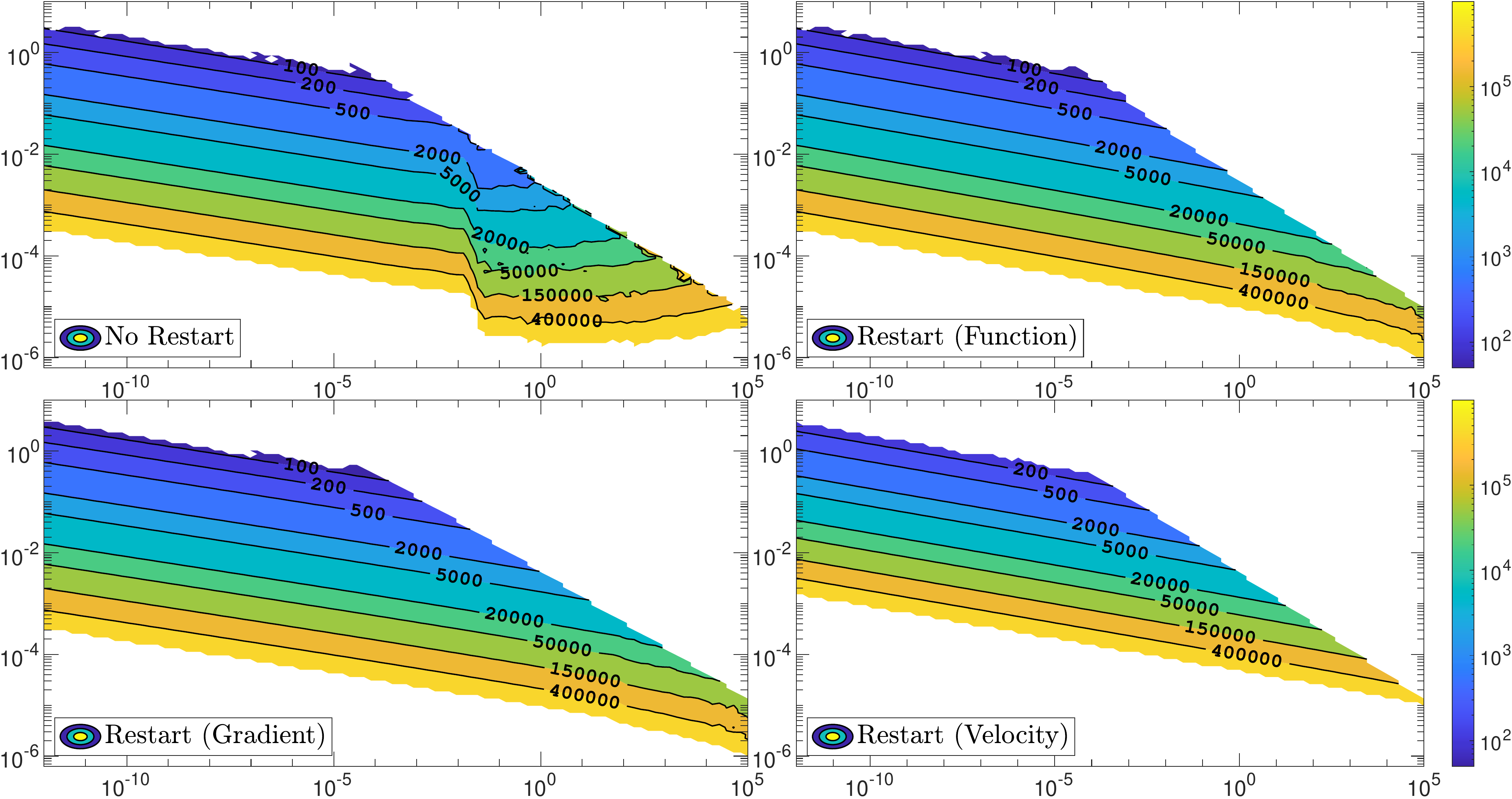} 	\vspace*{-7mm} \caption{Contour plot of the number of iterations required to achieve convergence in the $(C,h)$-plane, for the $p = \mathring{p} = 8$ PolyHTVI algorithm applied to Problem~\ref{Problem: Quartic}. }  \label{fig: RestartingRobustnessPoly3} 
\end{figure}

\begin{figure}[!ht] 
	\centering
	\includegraphics[width=0.96\textwidth]{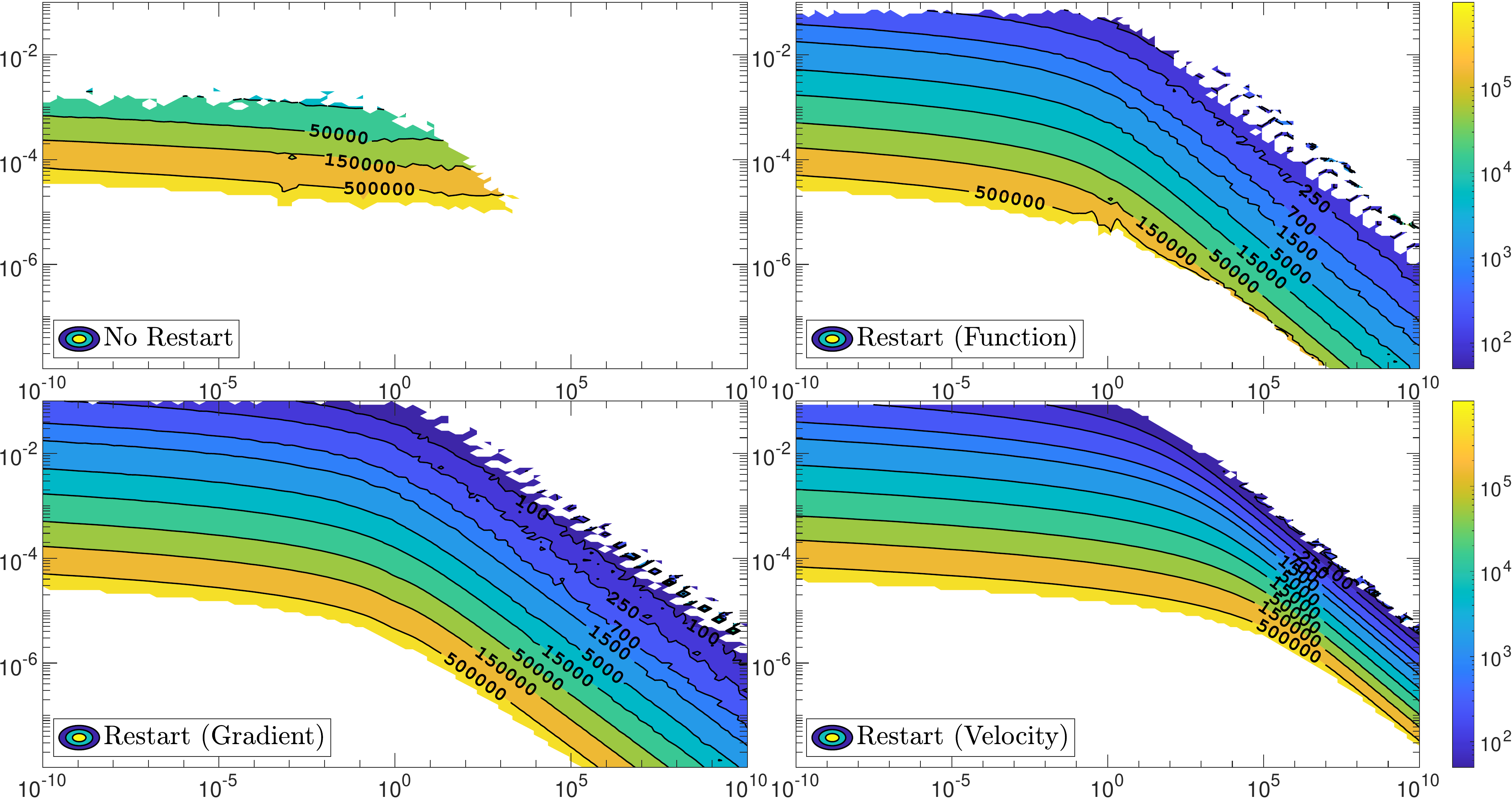} 	\vspace*{-1.5mm} \caption{ Contour plot of the number of iterations required to achieve convergence with $\delta = 10^{-5}$ in the $(C,h)$-plane, for $\eta = \mathring{\eta} = 1$ ExpoLSC applied to Problem~\ref{Problem: Polynomial and Log}. }\label{fig: RestartingRobustnessExpo}
\end{figure}

We have performed additional experiments to obtain a better idea of the benefits of momentum restarting in terms of computational efficiency, robustness and stability. More precisely, we solved optimization problems using the different versions of the algorithms on a $100\times100$ grid with logarithmic spacing in the parameter $(C,h)$-plane, and recorded the number of iterations required to achieve certain convergence criteria. Figures~\ref{fig: RestartingRobustnessPoly1}, \ref{fig: RestartingRobustnessPoly2}, \ref{fig: RestartingRobustnessPoly3}, and~\ref{fig: RestartingRobustnessExpo} display the results as filled contour plots (where the absence of color indicates either divergence or failure of the algorithm to converge in less than $10^6$ iterations). Table~\ref{table: Restarting Fastest Convergence} displays the number of iterations required to converge by each version of the algorithms with its optimal $(C,h)$ pair on the $100\times100$ logarithmically-spaced grid. %\\

Figure~\ref{fig: RestartingRobustnessPoly1} confirms the earlier observation that restarting can significantly reduce the number of iterations needed to converge, and we can also see that the restarted versions of the algorithm are more robust, since the regions of fast convergence are larger than for the standard algorithm. As a result, it is easier to tune the restarted algorithms to achieve fast convergence. Note as well from Figure~\ref{fig: RestartingRobustnessPoly2} that a restarting scheme can significantly improve the stability of the algorithms. Indeed, we can that as the convergence criteria are made stricter going from Figure~\ref{fig: RestartingRobustnessPoly1} to Figure~\ref{fig: RestartingRobustnessPoly2}, the regions of fast convergence have not shrunk as dramatically for the restarted algorithms as for the standard version. Given a converging $(C,h)$ pair for a restarted algorithm in Figure~\ref{fig: RestartingRobustnessPoly1}, the restarted algorithm usually remains convergent for that $(C,h)$ pair with the stricter criteria in Figure~\ref{fig: RestartingRobustnessPoly1} with a slightly increased number of iterations required. This is not true for the standard algorithm where the increase in number of iterations is much more significant, and there is a larger region of initially convergent $(C,h)$ pairs where the standard algorithm diverges when the stricter convergence criteria is imposed. \\


%\newpage 

As was observed earlier in Figure~\ref{fig: Restarting}, we can see from Figure~\ref{fig: RestartingRobustnessPoly3} that momentum restarting does not lead to significant improvements in computational efficiency for Problem~\ref{Problem: Quartic}, but also does not penalize computational efficiency in that case. From Figure~\ref{fig: RestartingRobustnessPoly3}, we see that this observation extends to robustness and stability. since all the different versions of the algorithm share similar convergence regions given the same parameter values and convergence criteria.

All the observations made so far also extend to the other Bregman subfamilies of dynamics and other algorithms, as can be seen, for instance, in Figure~\ref{fig: RestartingRobustnessExpo} for the ExpoSLC algorithm, where momentum restarting leads to significant gains in computational efficiency, robustness and stability. Table~\ref{table: Restarting Fastest Convergence} provides additional data supporting the significant gain in efficiency that can be achieved using momentum restarting.

\begin{table}[!ht]
	\centering
	\resizebox{\textwidth}{!}{
		\begin{tabular}{|c|c|c|c|c|c|c|}
			\hline
			\textbf{Algorithm} & \textbf{Problem}     & $\delta$   & \textbf{No Restart} & \textbf{Function Scheme} & \textbf{Gradient Scheme} & \textbf{Velocity Scheme} \\ \hline
			PolyHTVI  & Problem~\ref{Problem: Quartic} & $10^{-12}$ &      52      &        52            &             52       &      108             \\ \hline
			PolyHTVI  & Problem~\ref{Problem: Polynomial and Log}& $10^{-5}$  &      621       &          39        &    23          &       34        \\ \hline
			PolyHTVI  & Problem~\ref{Problem: Polynomial and Log} & $10^{-8}$  &    15994     &        80      &         51        &      57              \\ \hline
				PolyHTVI  & Problem~\ref{Problem: x log(x)} & $10^{-8}$ &    4121        &                   60 &          11          &       16            \\ \hline
				PolyHTVI  & Problem~\ref{Problem: Ill-Conditioned Quadratic}& $10^{-8}$ &       14723     &          60          &             12       &          14         \\ \hline
				PolyHTVI  & Problem~\ref{Problem: Linear Regression}& $10^{-5}$ &     3917   &     104             &       33            &    38            \\ \hline
					ExpoLSC & Problem~\ref{Problem: Quartic} & $10^{-12}$ &        75    &         68           &            64        &  155                 \\ \hline
			ExpoLSC  & Problem~\ref{Problem: Polynomial and Log} & $10^{-5}$  &      3929    &        50           &     20            &      21         \\ \hline
				ExpoLSC  & Problem~\ref{Problem: Polynomial and Log} & $10^{-8}$  &           204598 &        91          &      27      &       32           \\ \hline
					ExpoLSC & Problem~\ref{Problem: x log(x)} & $10^{-8}$ &       17081     &     66               &         15           &        18           \\ \hline
				ExpoLSC  & Problem~\ref{Problem: Ill-Conditioned Quadratic}& $10^{-8}$ &     58028       &       72             &              10     &          12         \\ \hline
					ExpoLSC & Problem~\ref{Problem: Linear Regression}& $10^{-5}$ &     21440   &      54            &       27            &      41          \\ \hline
	\end{tabular}}
	\caption{Comparison of the fastest convergence achieved by the standard algorithms and the restarting schemes on various problems with different tolerances $\delta$ (displayed in terms of number of iterations required to achieve the termination criteria).} 
	\label{table: Restarting Fastest Convergence}
\end{table}

Overall, all the numerical experiments conducted in this section unequivocally support the use of momentum restarting in the algorithms for accelerated optimization, and it can be seen from Figures~\ref{fig: RestartingRobustnessPoly1}, \ref{fig: RestartingRobustnessPoly2}, \ref{fig: RestartingRobustnessPoly3}, and~\ref{fig: RestartingRobustnessExpo} and Table~\ref{table: Restarting Fastest Convergence} that the gradient-based restarting scheme consistently outperforms the other two restarting schemes in terms of computational efficiency, robustness and stability. Unless stated otherwise, we will now always use momentum restarting based on the gradient scheme in the remainder of this paper, without explicitly stating it every time.  \\

\subsection{The Effect of the Parameter $C$} \label{section: Parameter C}

\hfill \\

The parameter $C$ in the polynomial and exponential subfamilies of Bregman dynamics, presented in Sections~\ref{subsec: Polynomial Subfamily} and \ref{subsec: Exponential Subfamily}, can sometimes provide a simple way to control the oscillatory behavior of the second-order differential equation. From the point of view of perturbation theory, the polynomial and exponential Bregman Euler--Lagrange equations \eqref{eq: Polynomial EL} and \eqref{eq: Exponential EL},
\begin{equation} 
	\ddot{q}(t)   + \frac{p+1}{t} \dot{q}(t) + Cp^2 t^{p-2}  \nabla f(q(t)) = 0, \quad \text{ and } \quad 	\ddot{q}(t) + \eta \dot{q} + C \eta^2 e^{\eta t} \nabla f(q(t)) = 0,
\end{equation}
can be thought of as perturbations of the simpler differential equations,
\begin{equation} 
	\ddot{u}(t)   + \frac{p+1}{t} \dot{u}(t)  = 0, \qquad \text{ and } \qquad 	\ddot{v}(t) + \eta \dot{v}  = 0.
\end{equation}
The solutions to these two unperturbed equations are given by
\begin{equation}
	u(t) = \left(k_1 t^{-p} + k_2 \right) \mathbb{1} ,  \qquad \text{ and } \qquad v(t) = \left( k_3 e^{-\eta t} + k_4 \right) \mathbb{1} ,
\end{equation}
for some constants $k_1,k_2,k_3,k_4$ depending on the initial conditions. They are non-oscillatory, and converge monotonically to a constant vector at the respective rates of $\mathcal{O}(t^{-p})$ and $\mathcal{O}(e^{-\eta t})$. We can thus think of the terms $ Cp^2 t^{p-2}  \nabla f(q(t))$ and  $C \eta^2 e^{\eta t} \nabla f(q(t)) $ as perturbations steering the dynamical system towards the minimizer of the objective function $f$, in an oscillatory fashion. The parameter $C$, which appears in front of these two perturbation terms, should therefore be chosen, in theory, to be small enough to control the oscillations but also large enough to guide the dynamical system towards the minimizer of the objective function. The situation is similar in the ExpoToPoly and PolyToExpo subfamilies of Bregman dynamics. % \\

This perturbation theoretic point of view and the numerical results which will be presented in this section suggest that the parameter~$C$ can play a very important role reducing the effect of oscillations and improving the performance of optimization algorithms. The benefits that tuning the parameter~$C$ can provide have not been sufficiently explored in the literature exploiting the variational framework for accelerated optimization (in \cite{JordanSymplecticOptimization,WiWiJo16,duruisseaux2020adaptive,Jordan2018,Campos2021} for instance), and the resulting dynamical systems were highly-oscillatory and thus required smaller time-steps for their discretizations. As a consequence, the  resulting optimization algorithms might not be as competitive as they could have been. Note as well that the limiting continuous differential equation for Nesterov's Accelerated Gradient introduced in \cite{SuBoCa16} can be thought of as the $p=2$ polynomial Bregman dynamics with $C=1/4$, which results in the highly oscillatory behavior observed in the continuous dynamics associated to most objective functions, and in the numerous discretizations of these dynamics which can be found in the literature. This observation also extends to the Riemannian manifold generalization of this variational framework for accelerated optimization~\cite{Duruisseaux2022Riemannian}, where the constant $C$ might not have been optimally tuned in practice (in \cite{Lee2021,Tao2020,Duruisseaux2022Riemannian,Duruisseaux2022Projection,Duruisseaux2022Constrained,Duruisseaux2022Lagrangian} for instance). % \\

\begin{figure}[!ht] 
	\centering
	\includegraphics[width=1\textwidth]{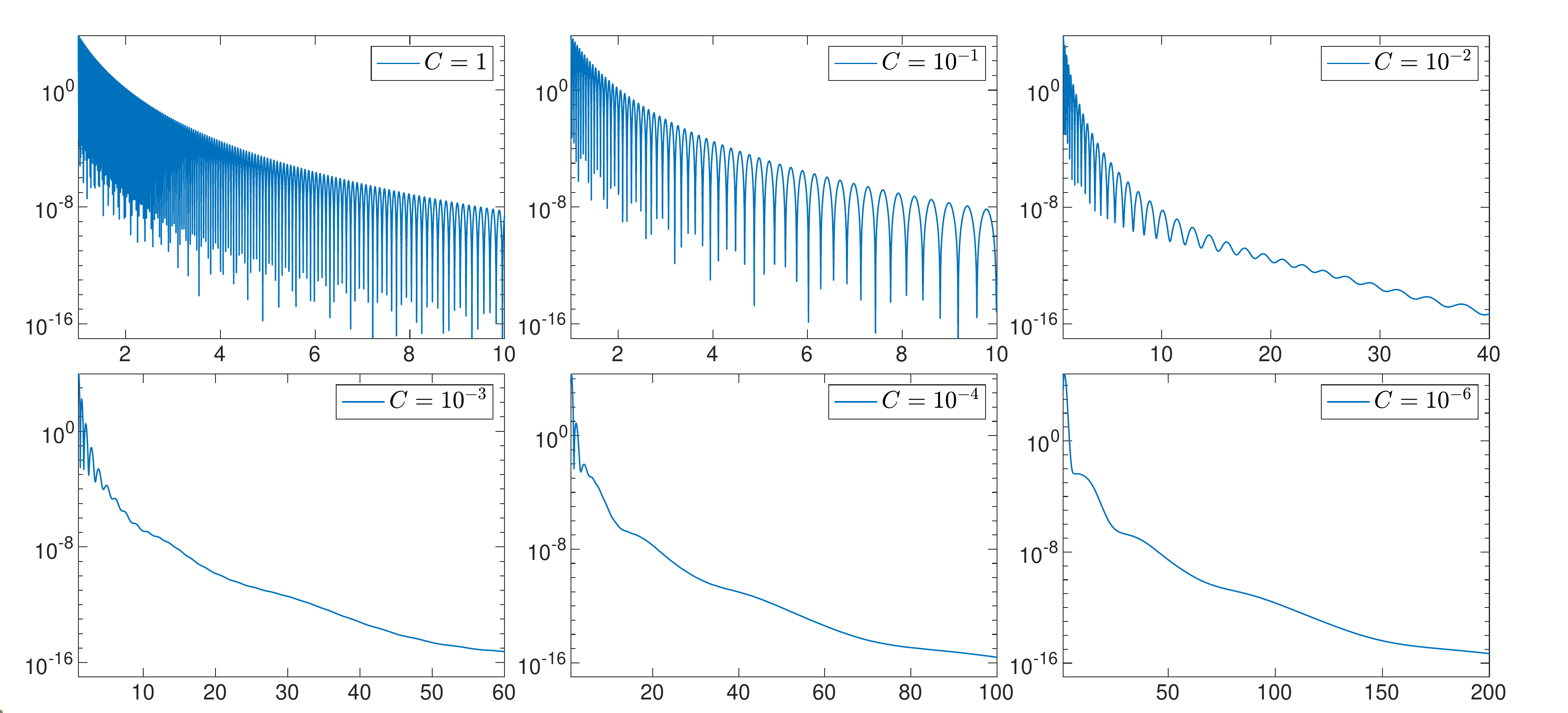} \vspace*{-7mm} \caption{ Error as a function of time $t$ along the $p = \mathring{p} = 6$ polynomial Bregman dynamics for Problem~\ref{Problem: Quartic}, with different values of the constant $C$.}\label{fig: CInvestigationTime1}
\end{figure}

\begin{figure}[!ht] 
	\centering
	\includegraphics[width=1\textwidth]{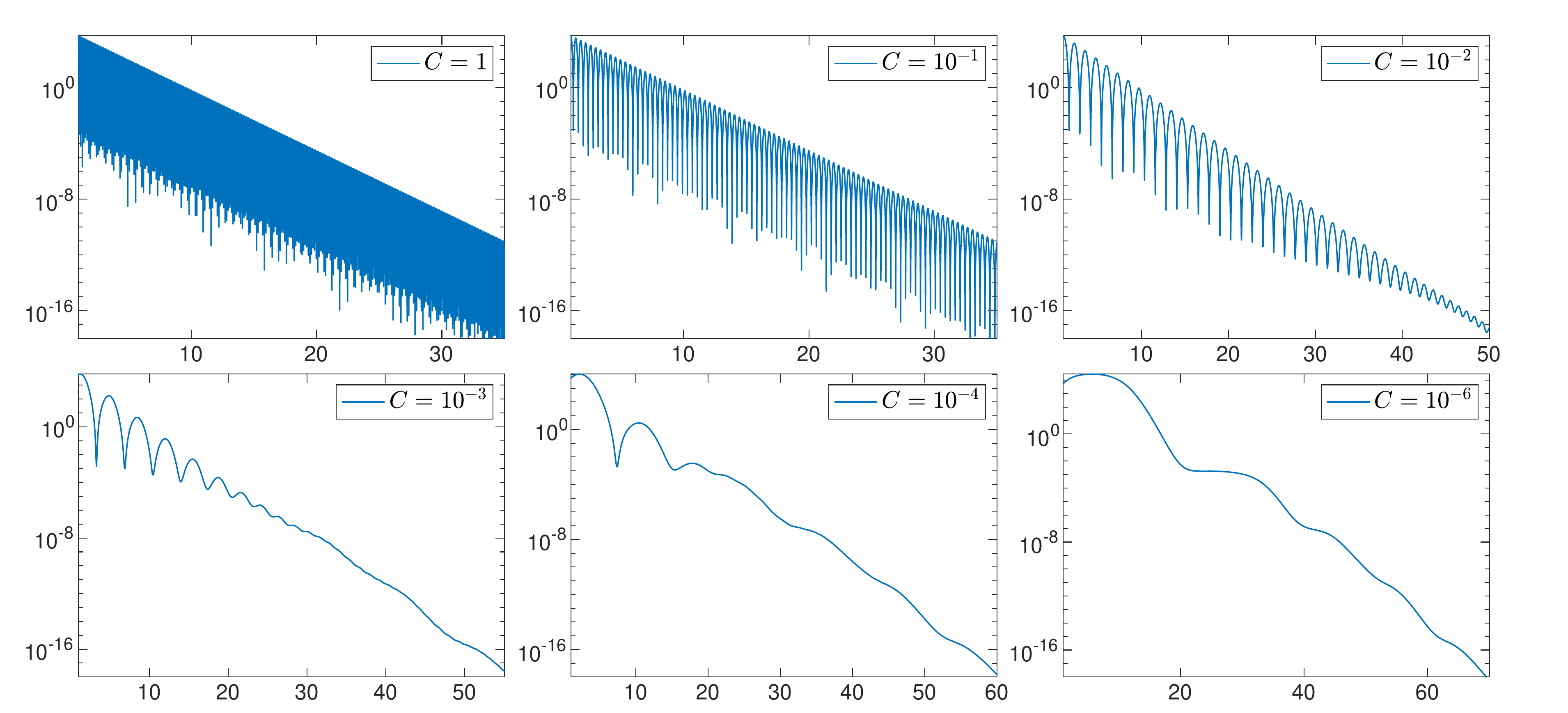} \vspace*{-7mm}  \caption{ Error as a function of time $t$ along the $\eta = \mathring{\eta} = 0.5$ exponential Bregman dynamics for Problem \ref{Problem: Quartic}, with different values of the constant $C$.}\label{fig: CInvestigationTime3}
\end{figure}

\begin{figure}[!ht] 
	\centering
	\includegraphics[width=1\textwidth]{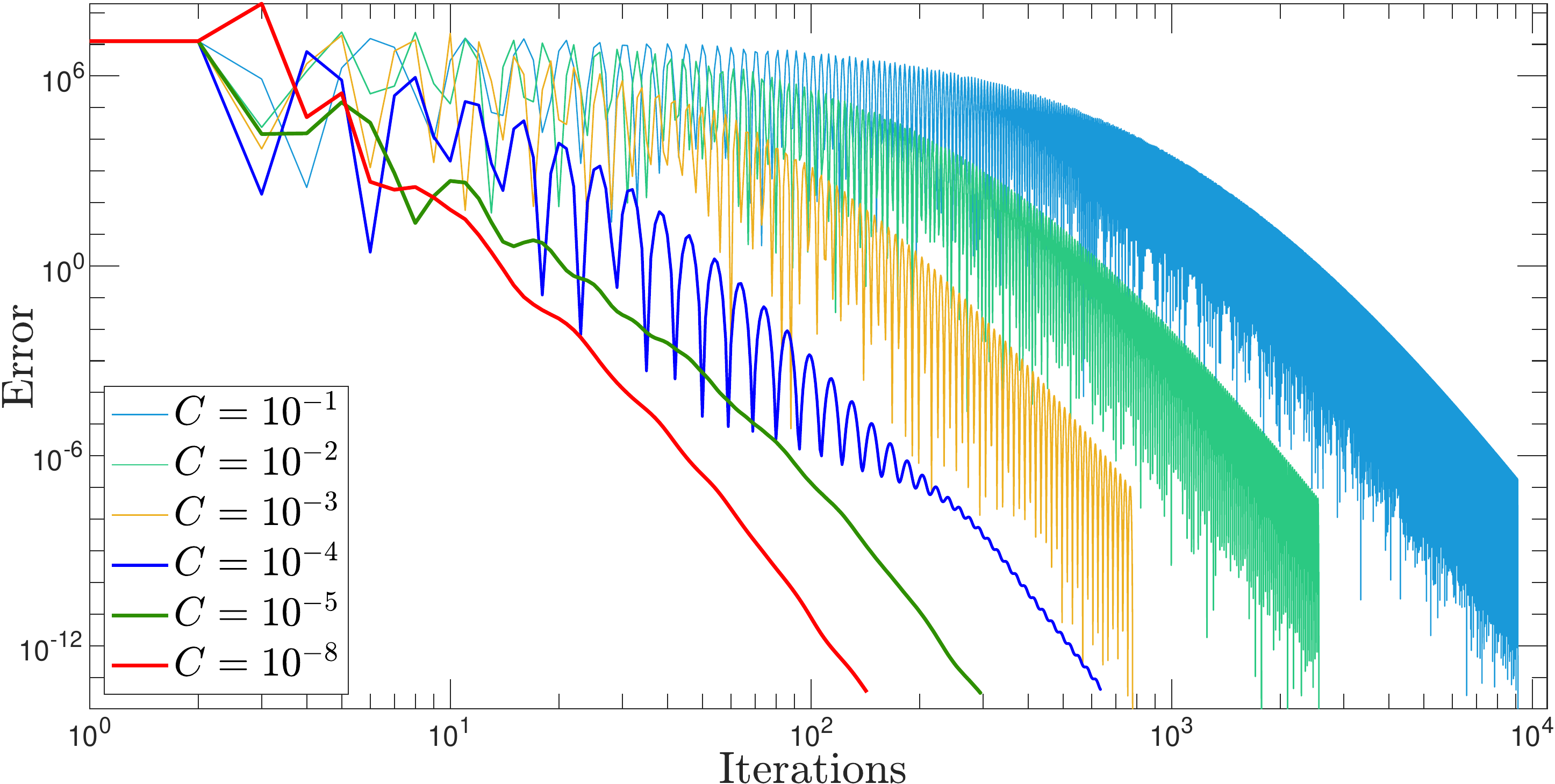} \vspace*{-5mm}  \caption{ Discretization of the polynomial Bregman dynamics using PolyHTVI with different values of $C$ for Problem \ref{Problem: Quartic} (without momentum restarting).   }\label{fig: CInvestigationDiscretization}
\end{figure}

\begin{figure}[!h] 
	\vspace*{-1mm}
	\centering
	\includegraphics[width=1\textwidth]{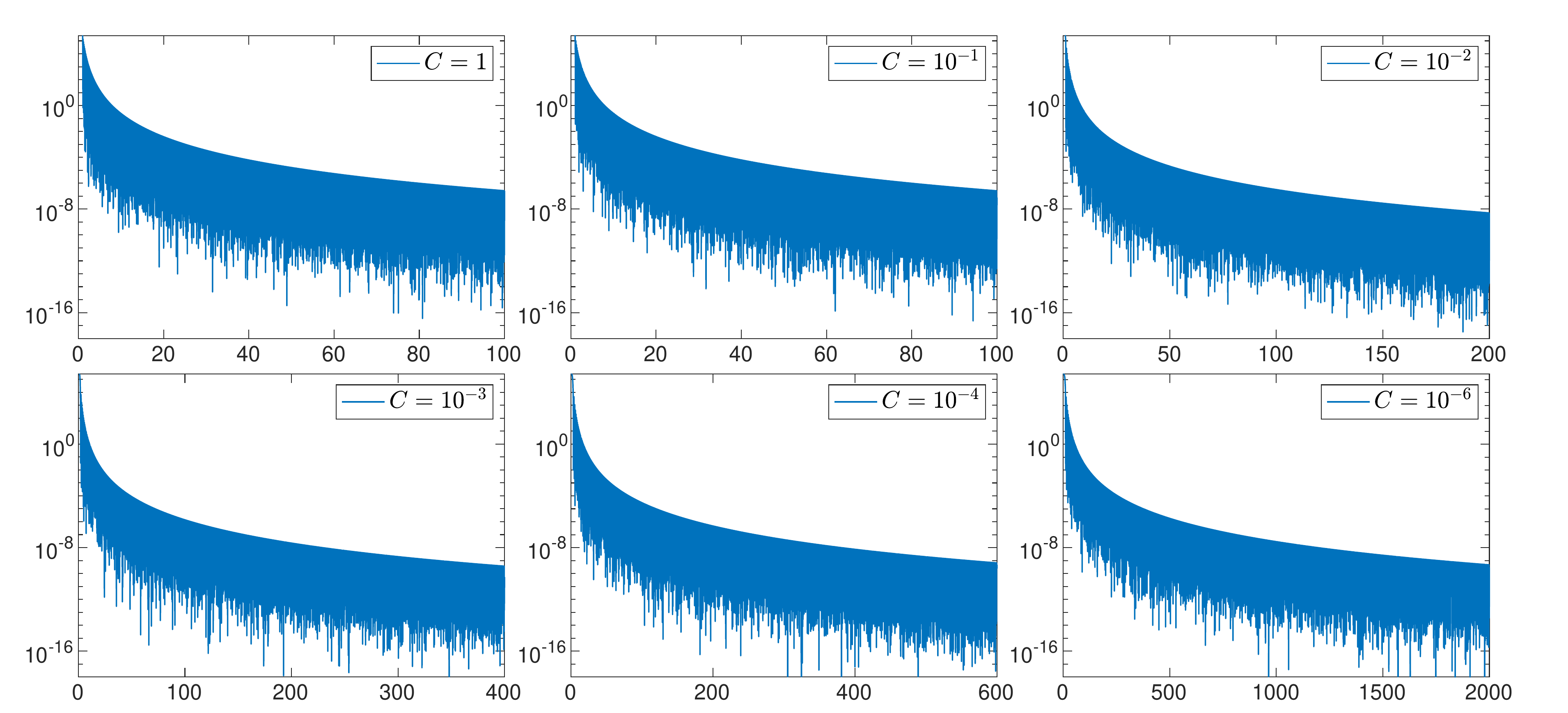}  \vspace*{-7.5mm} 
	\caption{ Error as a function of time $t$ along the $p = \mathring{p} = 6$ polynomial Bregman dynamics for Problem~\ref{Problem: Ill-Conditioned Quadratic}, with different values of the constant $C$. }
	\label{fig: CInvestigationTime2}
\end{figure}

As a first example, Figures~\ref{fig: CInvestigationTime1} and \ref{fig: CInvestigationTime3} display the changes in the polynomial and exponential Bregman dynamics for Problem~\ref{Problem: Quartic} as the parameter $C$ is decreased. The oscillations are clearly neutralized in the continuous Bregman dynamics as $C$ decreases. Although the convergence happens later in time for lower values of $C$, this is usually not an issue since the neutralization of the oscillations allows for larger time-steps when discretizing, as can be seen in Figure~\ref{fig: CInvestigationDiscretization}. This could also be seen in the `No Restart' contour plots presented in Figures~\ref{fig: RestartingRobustnessPoly1}, \ref{fig: RestartingRobustnessPoly2}, \ref{fig: RestartingRobustnessPoly3}, and \ref{fig: RestartingRobustnessExpo}, where lower values of~$C$ allowed for larger time-steps~$h$. Unfortunately, this behavior as $C$ is decreased does not seem to be universal, as can be seen from Figure~\ref{fig: CInvestigationTime2} for Problem~\ref{Problem: Ill-Conditioned Quadratic}. % \\

%\newpage 

%\hfill 

%\hfill 

%\newpage 

\begin{figure}[!h] 
	\vspace*{-4mm}
	\centering
	\includegraphics[width=1\textwidth]{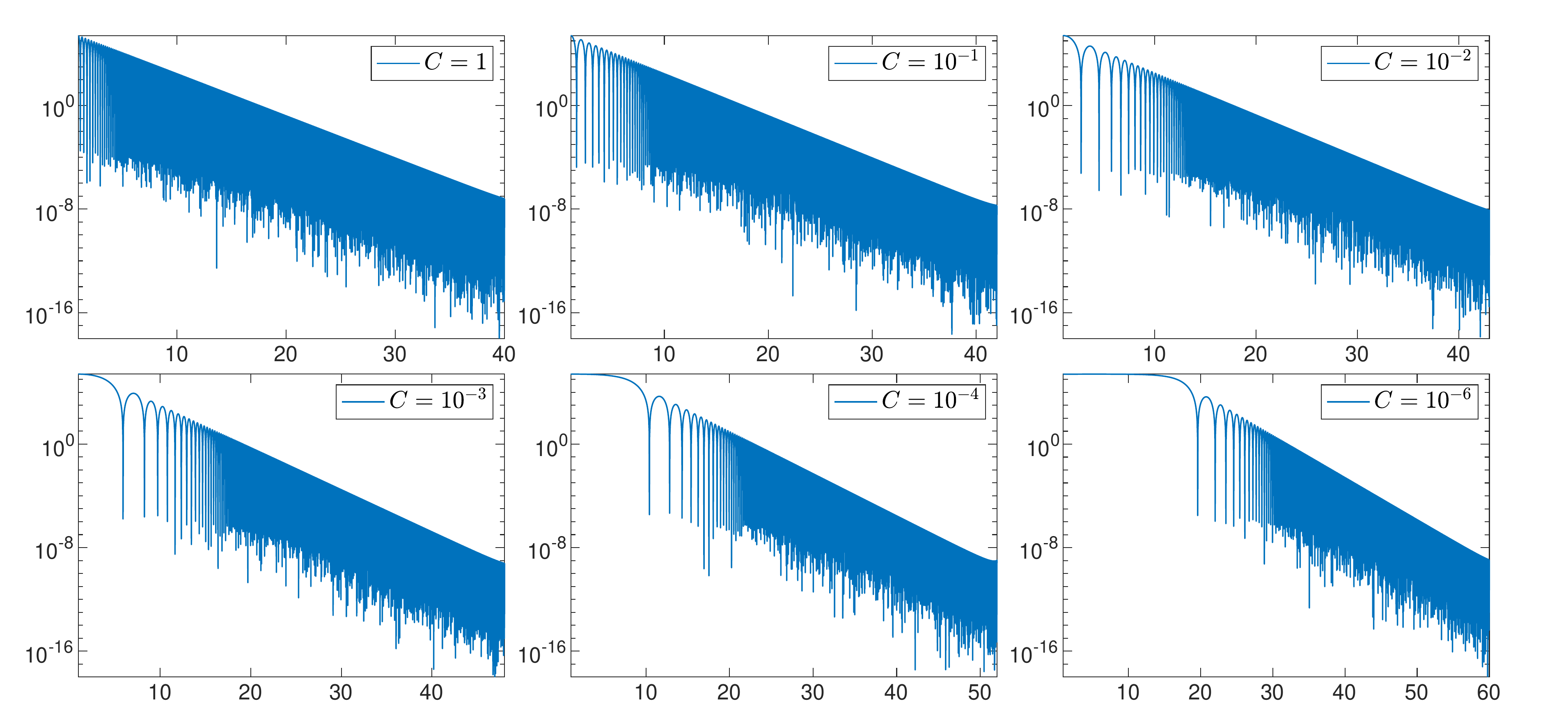} 	\vspace*{-7.5mm}	\caption{ Error as a function of time $t$ along the $\eta = \mathring{\eta} = 0.5$ exponential Bregman dynamics for Problem \ref{Problem: Ill-Conditioned Quadratic}, with different values of the constant $C$.}  \label{fig: CInvestigationTime4} \vspace*{-1mm}
\end{figure}

 We will now try to obtain a better understanding of the dependence on $C$ of the computational efficiency of the optimization algorithms, and of how a good choice of parameter $C$ depends on the other variables. Preliminary experiments showed that the convergence regions are very similar for the different algorithms within the same adaptive family of Bregman dynamics, so we will only use the HTVI algorithms here but the results extend to the other algorithms. We will return to the comparison of the different geometric integrators later in Section~\ref{sec: Comparison Integrators}.

\begin{figure}[!h] 
	\vspace{3mm}
	\centering
	\includegraphics[width=1\textwidth]{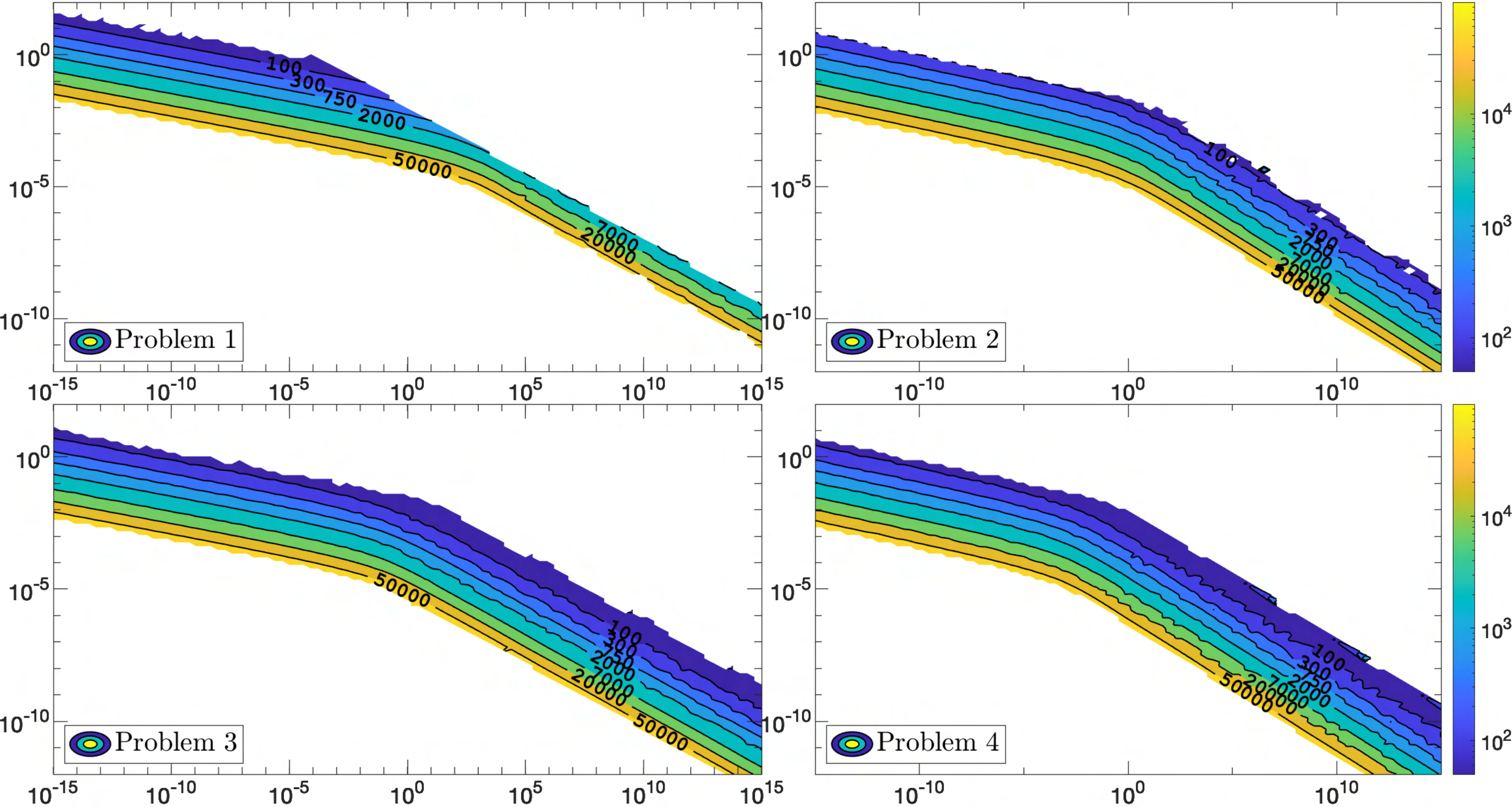} 	\vspace*{-5mm} \caption{ Contour plot of the number of iterations required to achieve convergence ($\delta = 10^{-6}$) in the $(C,h)$-plane, for $p = \mathring{p} = 6$ PolyHTVI applied to Problems~\ref{Problem: Quartic}, \ref{Problem: Polynomial and Log}, \ref{Problem: x log(x)}, \ref{Problem: Ill-Conditioned Quadratic}. } \label{fig: Poly4Problems} 
\end{figure}

\begin{figure}[!h] 
	\vspace{5mm}
	\centering
	\includegraphics[width=1\textwidth]{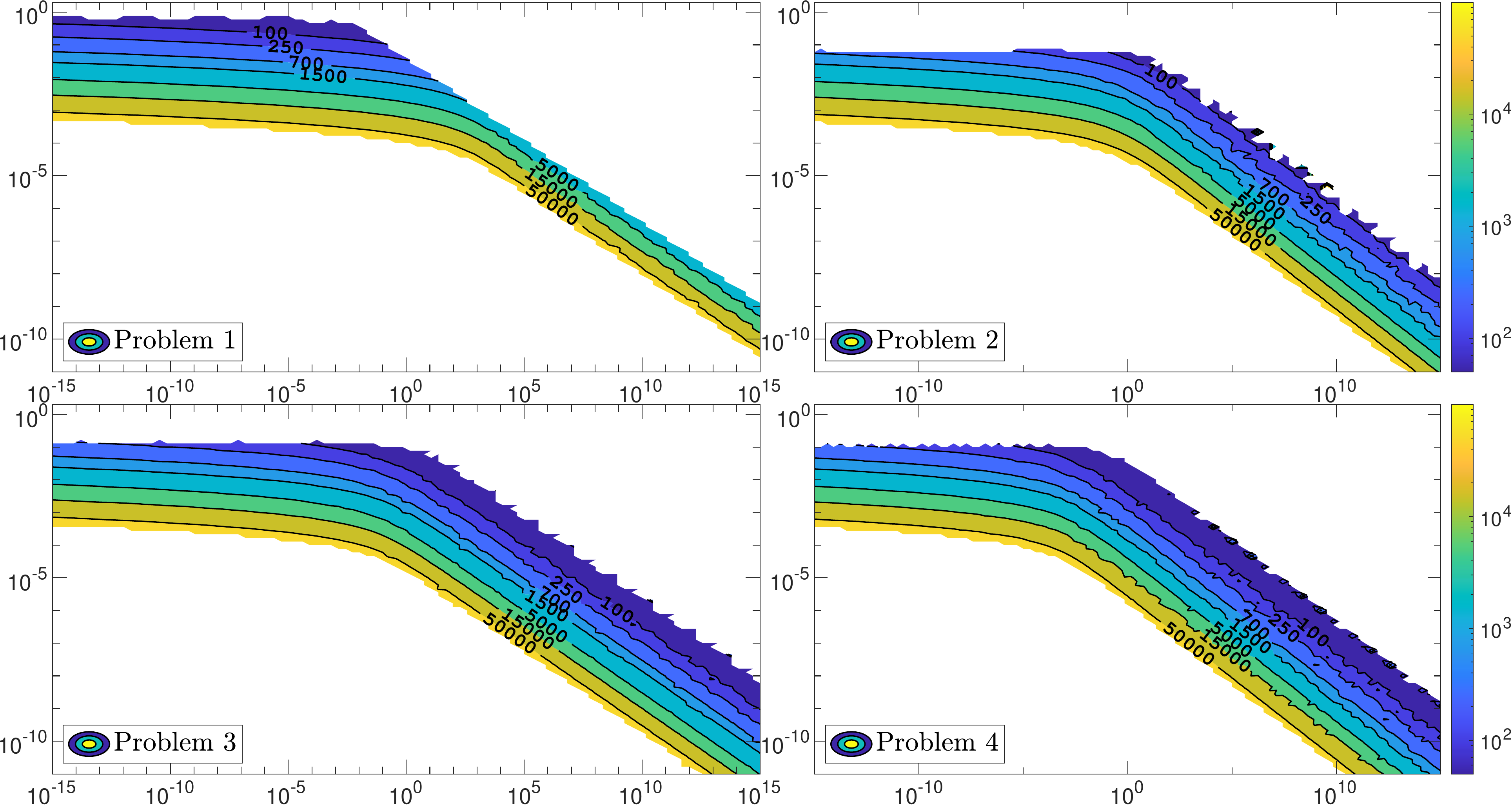} 	\vspace*{-5mm} \caption{ Contour plot of the number of iterations required to achieve convergence ($\delta = 10^{-6}$) in the $(C,h)$-plane, for $\eta = \mathring{\eta} = 6$ ExpoHTVI applied to Problems~\ref{Problem: Quartic}, \ref{Problem: Polynomial and Log}, \ref{Problem: x log(x)}, \ref{Problem: Ill-Conditioned Quadratic}.} \label{fig: Expo4Problems} 
	\vspace{2mm}
\end{figure}

Let us first investigate whether the regions of optimal convergence in the $(C,h)$-plane are problem-dependent. Figures~\ref{fig: Poly4Problems} and \ref{fig: Expo4Problems} display the convergence regions obtained when using the HTVI algorithm for the polynomial and exponential Bregman dynamics on four different objective functions. Unfortunately, these results show that the regions of optimal convergence are problem-dependent and as a result we will not be able to find a single value of $C$ which will achieve almost-optimal performance on all problems.  \\

%\newpage 


\begin{figure}[!h] 
	\vspace{3mm}
	\centering
	\includegraphics[width=1\textwidth]{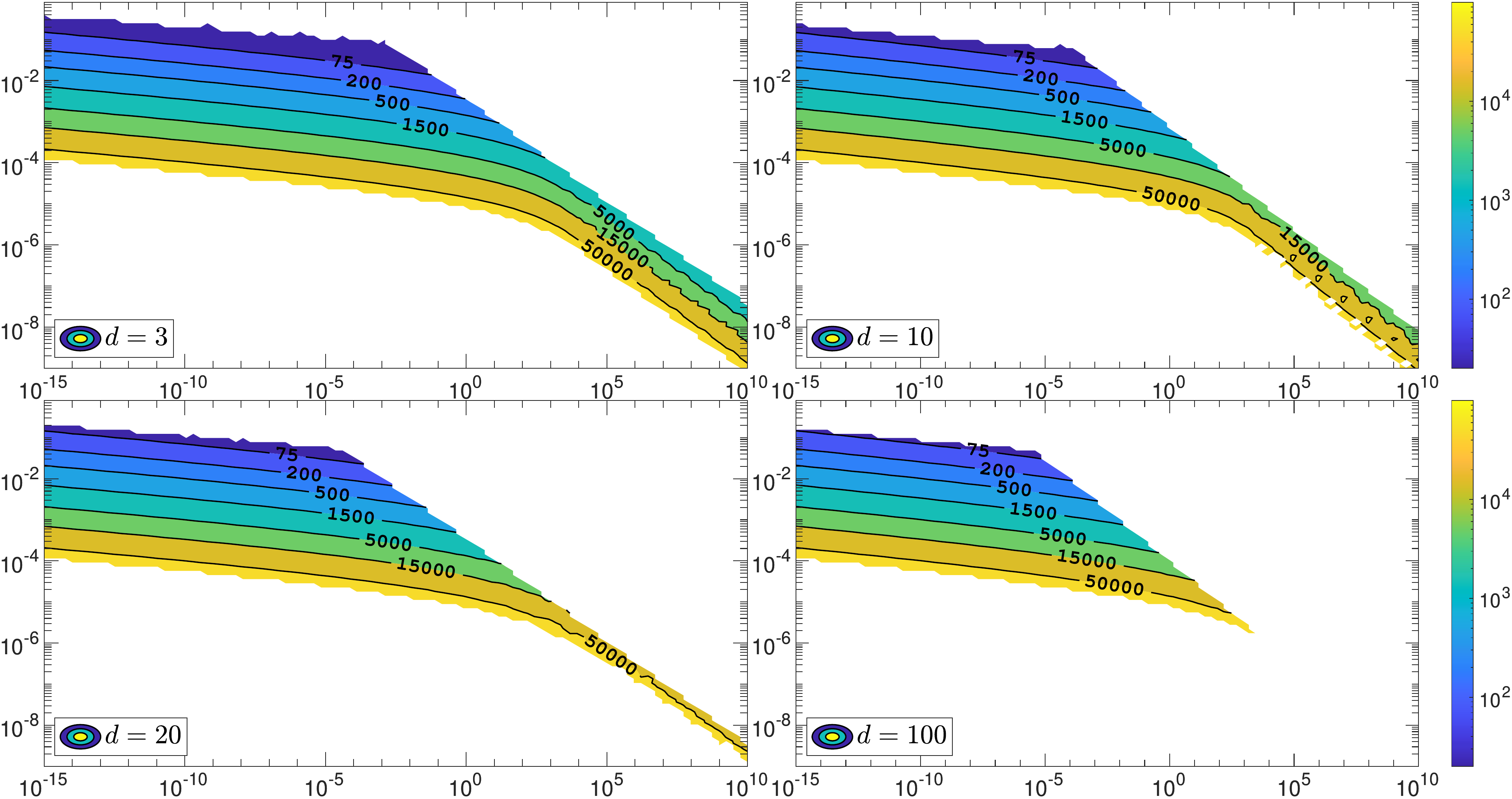} 	\vspace*{-5mm} \caption{ Contour plot of the number of iterations required to achieve convergence ($\delta = 10^{-6}$) in the $(C,h)$-plane, for $p = 6, \mathring{p} = 2$ PolyHTVI applied to Problem~\ref{Problem: Quartic} with different dimensions $d$. } \label{fig: Poly4Dim1} 
\end{figure}

\begin{figure}[!h] 
	\vspace*{5mm}
	\centering
	\includegraphics[width=1\textwidth]{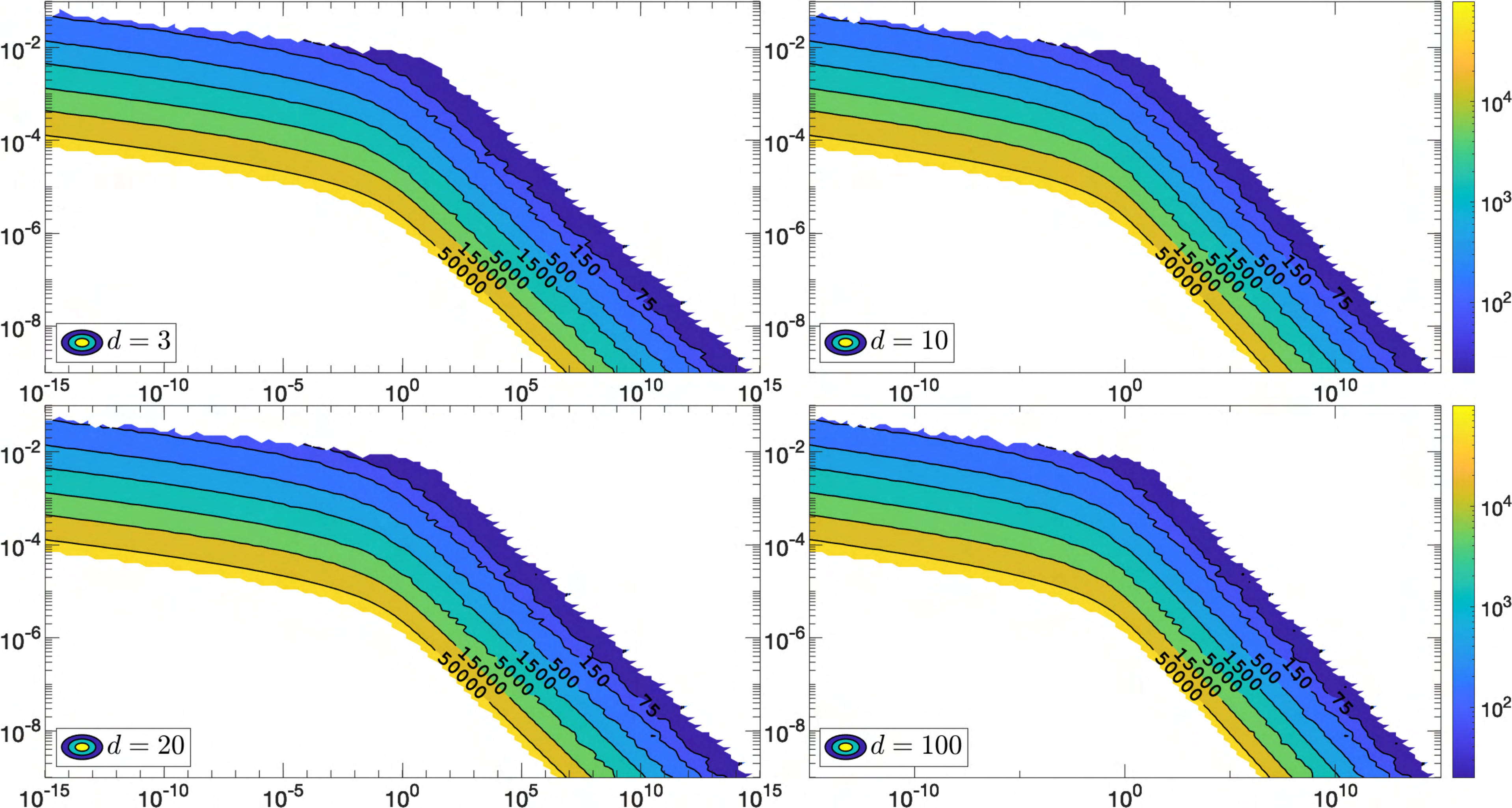} 	\vspace*{-5mm} \caption{ Contour plot of the number of iterations required to achieve convergence ($\delta = 10^{-6}$) in the $(C,h)$-plane, for $p = 6, \mathring{p} = 2$ PolyHTVI applied to Problem~\ref{Problem: x log(x)} with different dimensions $d$. } \label{fig: Poly4Dim2} 
\end{figure}

\begin{figure}[!h] 
	 \vspace{2mm}
	\centering
	\includegraphics[width=1\textwidth]{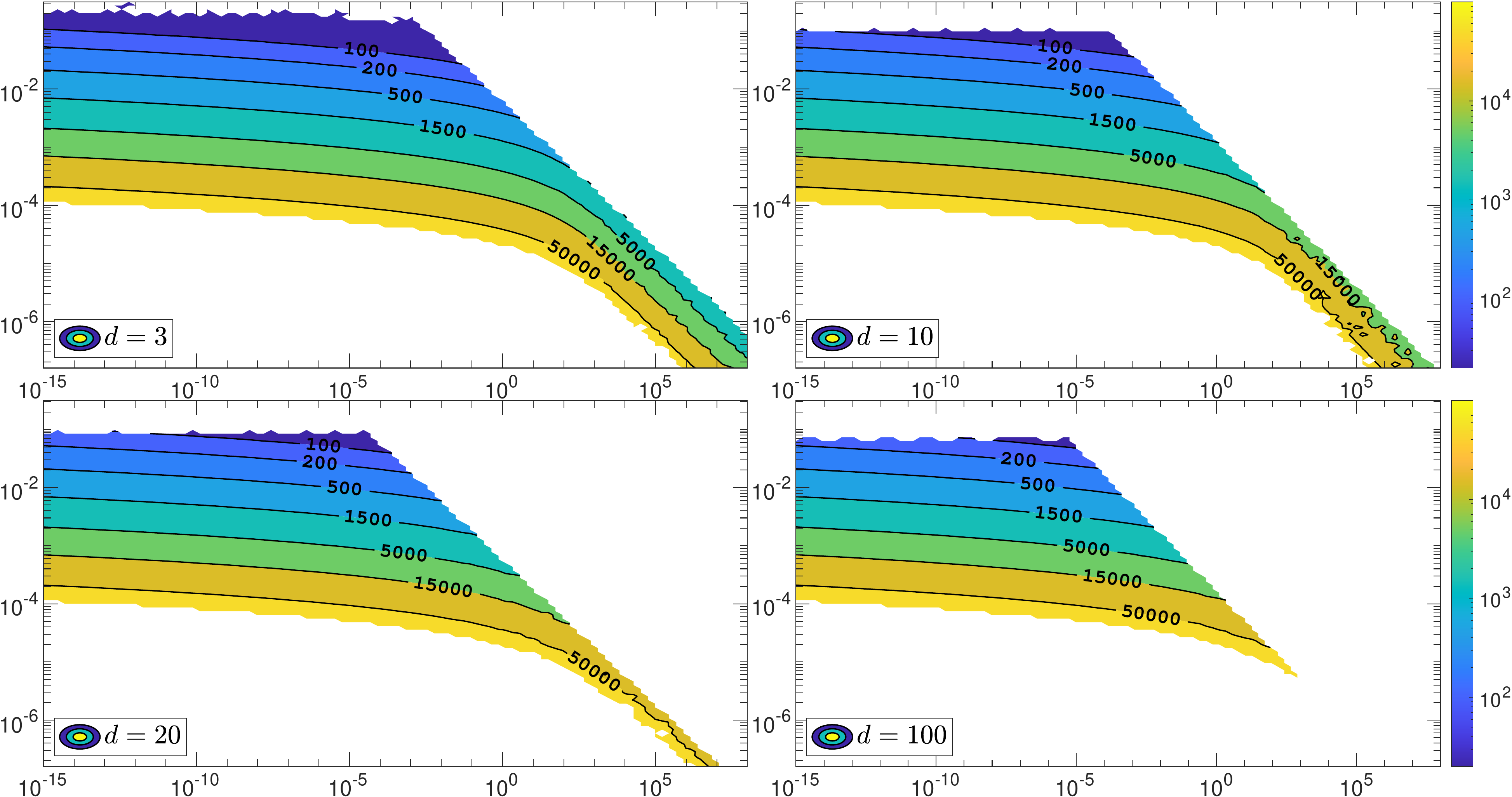} 	\vspace*{-6mm} \caption{ Contour plot of the number of iterations required to achieve convergence ($\delta = 10^{-6}$) in the $(C,h)$-plane, for $\eta = 2, \mathring{\eta} = 1$ ExpoHTVI applied to Problem~\ref{Problem: Quartic} with different dimensions $d$. } \label{fig: Expo4Dim}  \vspace{2mm}
\end{figure}

However, from Figures~\ref{fig: Poly4Dim1}, \ref{fig: Poly4Dim2} and \ref{fig: Expo4Dim}, we can see that for fixed values of $p,\mathring{p},\eta, \mathring{\eta}$, the convergence regions in the $(C,h)$-plane are left almost unchanged as the dimension of the problem is increased from $d=3$ to $d=100$, although the numbers of iterations required increase slightly with $d$. This observation can improve significantly the process of tuning the optimization algorithm for high-dimensional problems by first tuning the algorithm on a similar low-dimensional problem, which could be particularly helpful for certain machine learning applications.  \\


%\newpage 

Note that all the observations made in this section extend to the ExpoToPoly and PolyToExpo subfamilies of time-adaptive Bregman dynamics.  \\

\subsection{Other Approaches to Control Oscillations}

\hfill \\

There are other possible approaches to control the oscillations in second-order nonlinear differential equations. One such method is Hessian-driven damping \cite{Alvarez2002,Attouch2020,Attouch2021,Attouch2022}, where the idea is to add a damping term which involves the Hessian of the objective function, $\beta(t) \nabla^2 f(x(t)) \dot{x}(t)$, to the differential equation of interest:
\begin{equation}
	\ddot{x}(t) + \gamma(t) \dot{x}(t) + \beta(t) \nabla^2 f(x(t)) \dot{x}(t) + b(t) \nabla f(x(t)) = 0.
\end{equation}
The addition of this Hessian-driven damping term appears to neutralize the oscillations in the continuous solution to the differential equation. Furthermore, it was shown using Lyapunov analysis that under suitable assumptions, solutions to the modified equation not only satisfied a similar convergence rate to the minimizer as solutions to the original equation, but also benefited from additional convergence properties for the norm of the gradient $\nabla f$. First-order optimization algorithms were also derived by discretizing the modified differential equation, after rewriting $\nabla^2 f(x(t)) \dot{x}(t)$ as $\frac{d}{dt} \nabla f(x(t))$. Unfortunately, we cannot derive a simple variational formulation for this modified differential equation, so we cannot easily incorporate Hessian-driven damping into our framework which relies on geometric numerical integration of Lagrangian or Hamiltonian systems. % \\

Another possible approach to control oscillations consists in simplifying the Bregman dynamics using local approximations. For instance, one could integrate local linearizations of the Bregman Hamilton's equations, or start from local quadratic Hamiltonian approximations to the Bregman Hamiltonian, or use a local quadratic model for the objective function. We will not consider these methods here because they can suffer from additional numerical stability issues coming from the approximations at play, and it can be very challenging to design a symplectic integrator which preserves the nice properties of the dynamics across all the different local approximations. % \\

A different approach consists in designing a symplectic integrator which can travel faster along the oscillations via larger time-steps. This may be achievable using Spectral or Galerkin variational integrators \cite{LeZh2011,MaWe2001,Hall2015,LeSh2011_sbvi}, which rely on a choice of basis functions that span a good approximation space for the Bregman dynamics (for instance, simulations of the polynomial Bregman dynamics suggest that the error usually follows a trajectory which can be well-approximated using functions of the form $	t^{-\gamma} \cos{(\alpha  t^\beta )} $ or $t^{-\gamma} \sin{(\alpha  t^\beta )} $, where $\gamma$ is the decay rate, $\alpha$ tunes the frequency of oscillations, and $\beta \in (0,1)$ characterizes the slowing down of the oscillation frequency). Due to the oscillatory nature of the dynamical system, it might also be advantageous to use Filon-type \cite{Filon1930} or Levin-type \cite{Levin1982,Levin1996} quadrature rules in the construction of the integrators, since they are designed specifically for highly oscillatory integrals (see \cite{IserlesBook} for a thorough presentation). % \\

Another possibility involves averaging techniques \cite{Sanders2007,Verhulst1996,Farr2009,ScLe2017}. The extended Bregman Hamiltonians or Lagrangians can be split as \begin{equation} H(\bar{q},\bar{r}) = H^A(\bar{q},\bar{r}) + C H^B(\bar{q}), \end{equation} or \begin{equation}  L(\bar{q},\bar{q}' ) = L^A(\bar{q},\bar{q}' ) + C L^B(\bar{q}) ,\end{equation}	where the $A$-component is dominating and can be solved exactly (or efficiently approximated with high accuracy), and the $B$-component generates small perturbations affecting the overall dynamics. One can then hope to integrate the dominating dynamics very accurately with larger time-steps and incorporate the influence of the small perturbations by averaging them out. Unfortunately, although this approach seemed to neutralize the oscillations in the solution in practice, it did not allow the use of larger time-steps, and the resulting algorithm was actually less competitive and robust because of the implicit nature of the update for the momentum $r$. 

\hfill \\

\section{Time-Adaptivity in the Momentum Restarted Algorithms}  \label{section: Time-Adaptivity in the Momentum Restarted Algorithms}

We will first investigate how the optimization algorithms behave as the parameters $\mathring{p}$ and $\mathring{\eta}$ are varied. In \cite{duruisseaux2020adaptive,Duruisseaux2022Riemannian,Duruisseaux2022Constrained,Duruisseaux2022Projection,Duruisseaux2022Lagrangian}, numerical experiments with the polynomial Bregman dynamics suggested that time-adaptivity (i.e. using $\mathring{p} \ne p$) could result in significantly faster optimization algorithms due to the exponentially growing time-steps that they use (instead of constant time-steps for $\mathring{p} = p$). These numerical experiments were however carried with the standard versions of the algorithms and without a careful tuning of the parameter $C$. %\\

\begin{figure}[!h] 
	\vspace{4mm}
	\centering
	\includegraphics[width=1\textwidth]{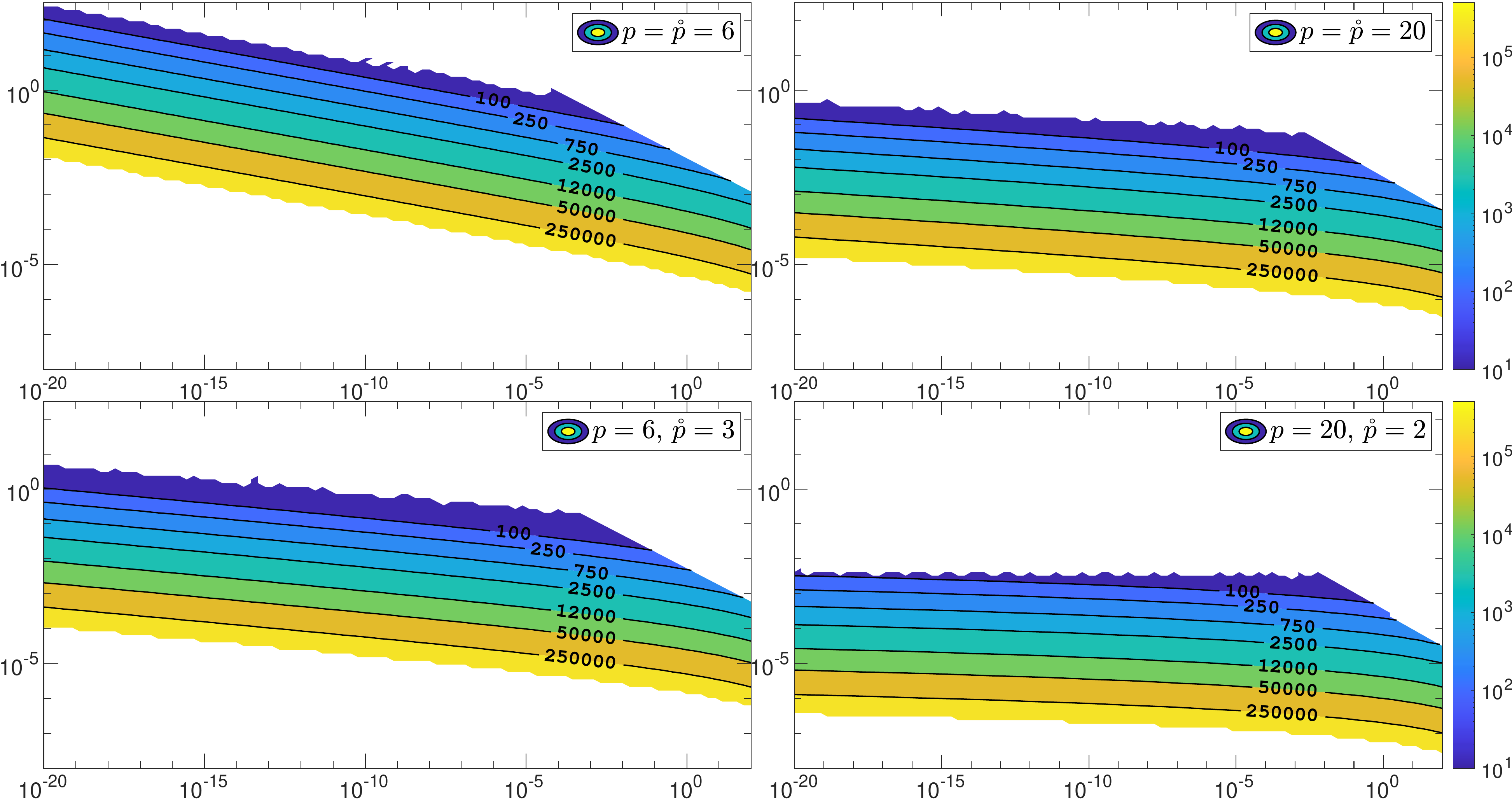} 	\vspace*{-5mm} \caption{ Contour plot of the number of iterations required to achieve convergence ($\delta = 10^{-6}$) in the $(C,h)$-plane, for PolyHTVI applied to Problem~\ref{Problem: Quartic}. } \label{fig: Poly4Param1} 
\end{figure}

\begin{figure}[!h] 
	\vspace{4mm}
	\centering
	\includegraphics[width=1\textwidth]{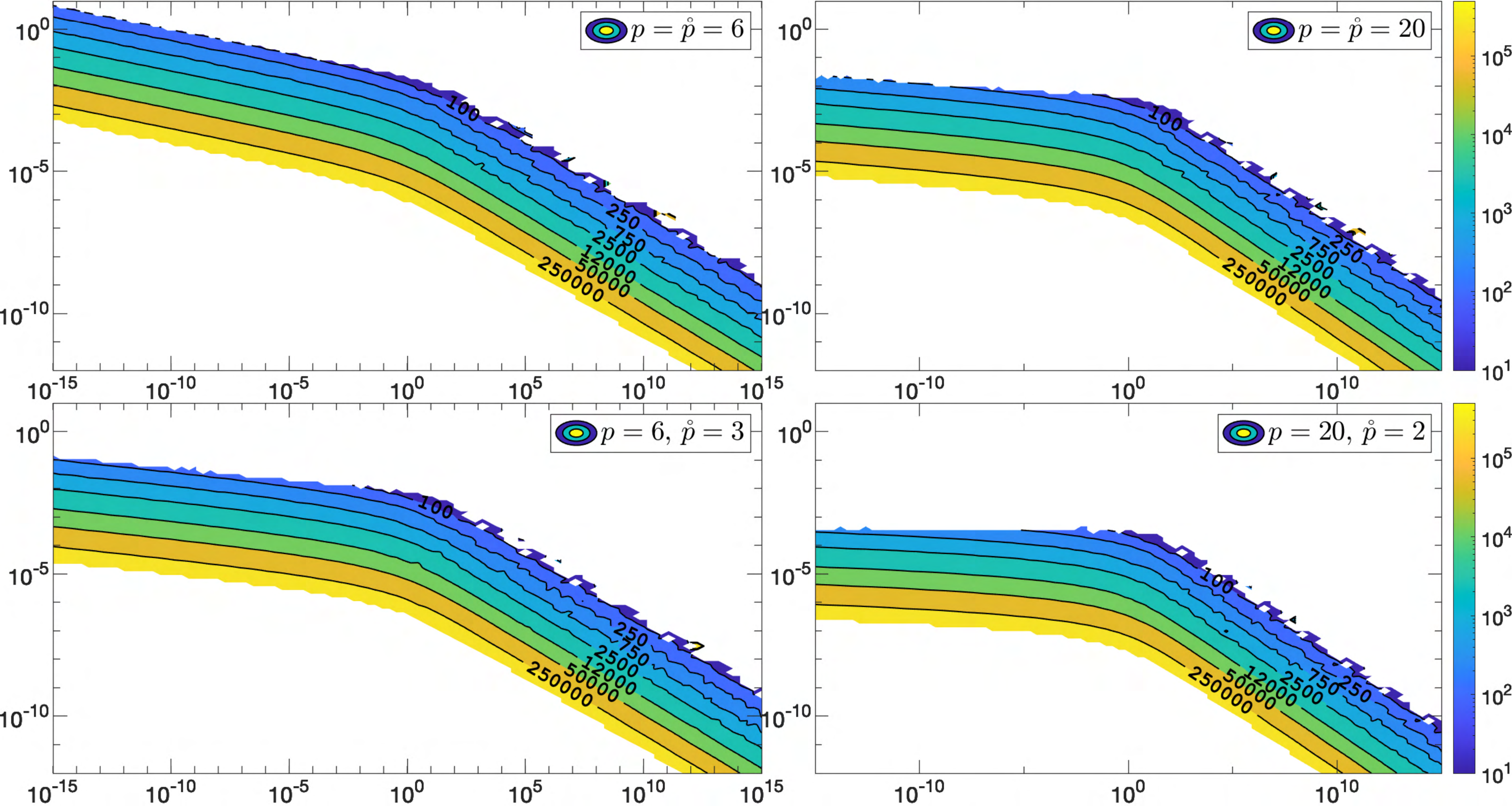} 	\vspace*{-5mm} \caption{ Contour plot of the number of iterations required to achieve convergence ($\delta = 10^{-6}$) in the $(C,h)$-plane, for PolyHTVI applied to Problem~ \ref{Problem: Polynomial and Log}. } \label{fig: Poly4Param2} 
	\vspace{2mm}
\end{figure}

\begin{figure}[!h] 
	\vspace{4mm}
	\centering
	\includegraphics[width=1\textwidth]{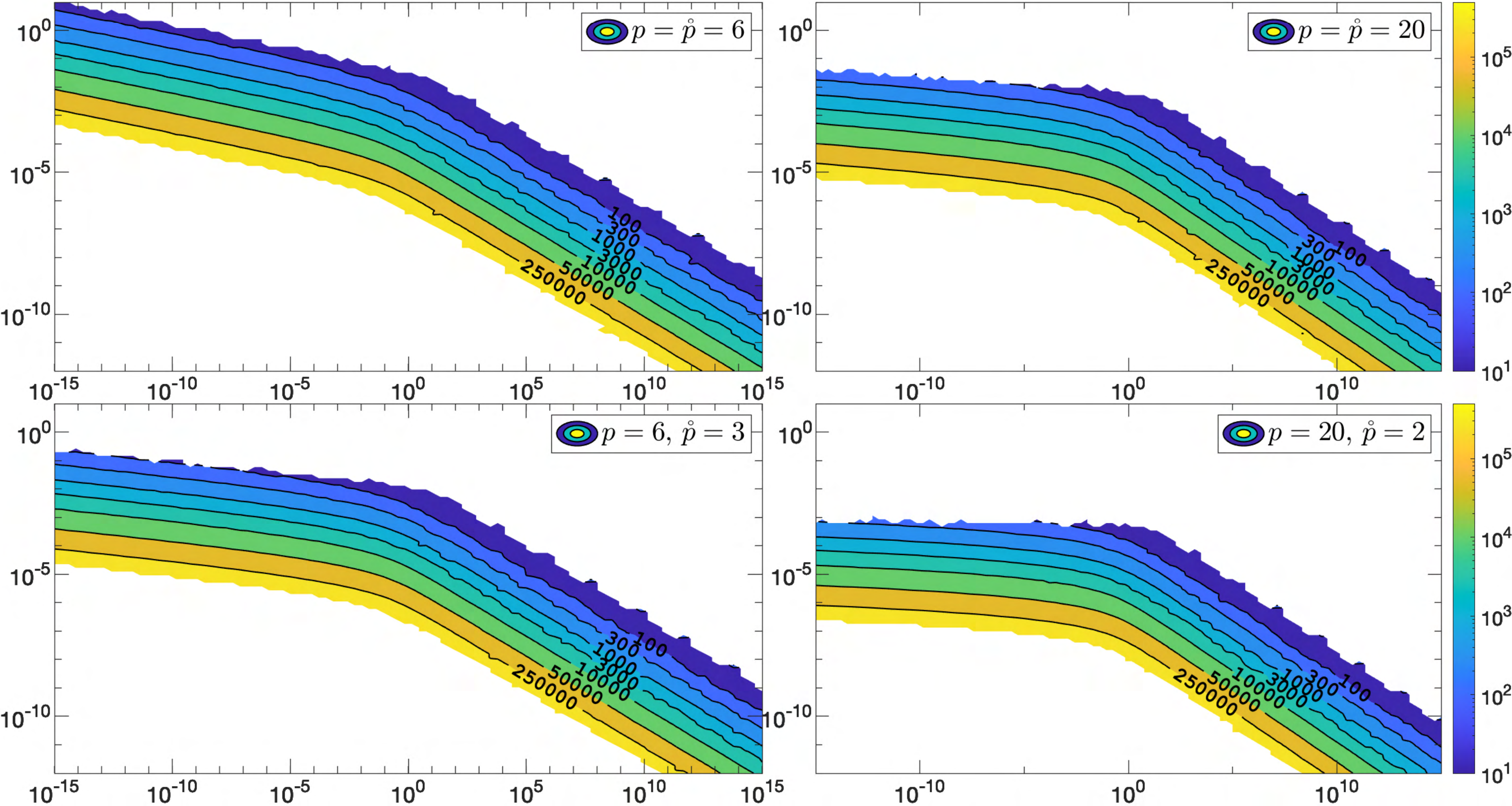} 	\vspace*{-6mm} \caption{ Contour plot of the number of iterations required to achieve convergence ($\delta = 10^{-6}$) in the $(C,h)$-plane, for PolyHTVI applied to Problem~\ref{Problem: x log(x)}. } \label{fig: Poly4Param3} 
\end{figure}

\begin{figure}[!h] 
		\vspace{4mm}
	\centering
	\includegraphics[width=1\textwidth]{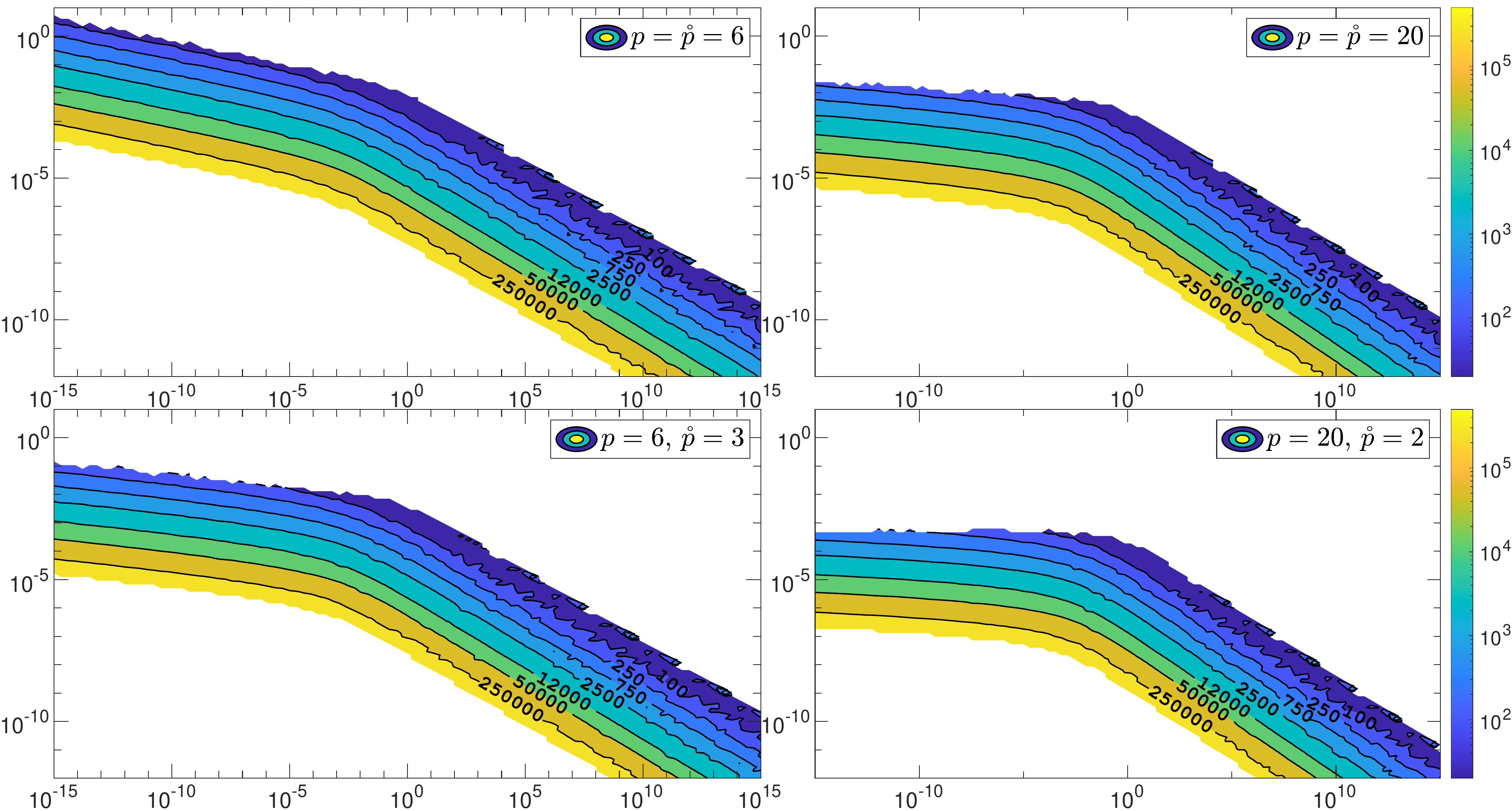} 	\vspace*{-6mm} \caption{Contour plot of the number of iterations required to achieve convergence ($\delta = 10^{-8}$) in the $(C,h)$-plane, for PolyHTVI applied to Problem~\ref{Problem: Ill-Conditioned Quadratic}. } \label{fig: Poly4Param4} 
	\vspace{2mm}
\end{figure}

\begin{figure}[!h] 
	\vspace{4mm}
	\centering
	\includegraphics[width=1\textwidth]{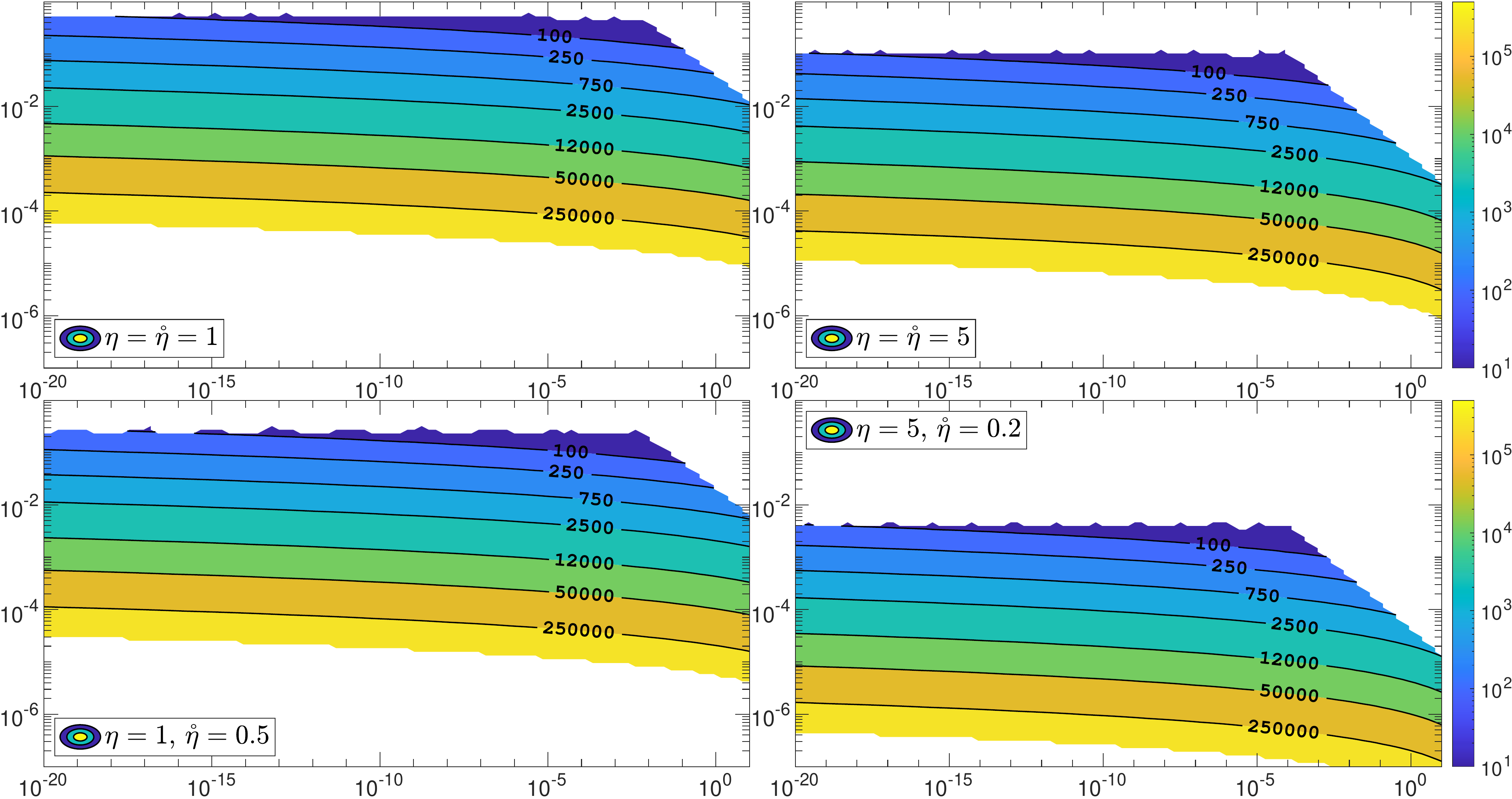} 	\vspace*{-5mm} \caption{ Contour plot of the number of iterations required to achieve convergence ($\delta = 10^{-7}$) in the $(C,h)$-plane, for ExpoHTVI applied to Problem~\ref{Problem: Quartic}. } \label{fig: Expo4Param1} 
\end{figure}

\begin{figure}[!h] 
	\vspace{4mm}
	\centering
	\includegraphics[width=1\textwidth]{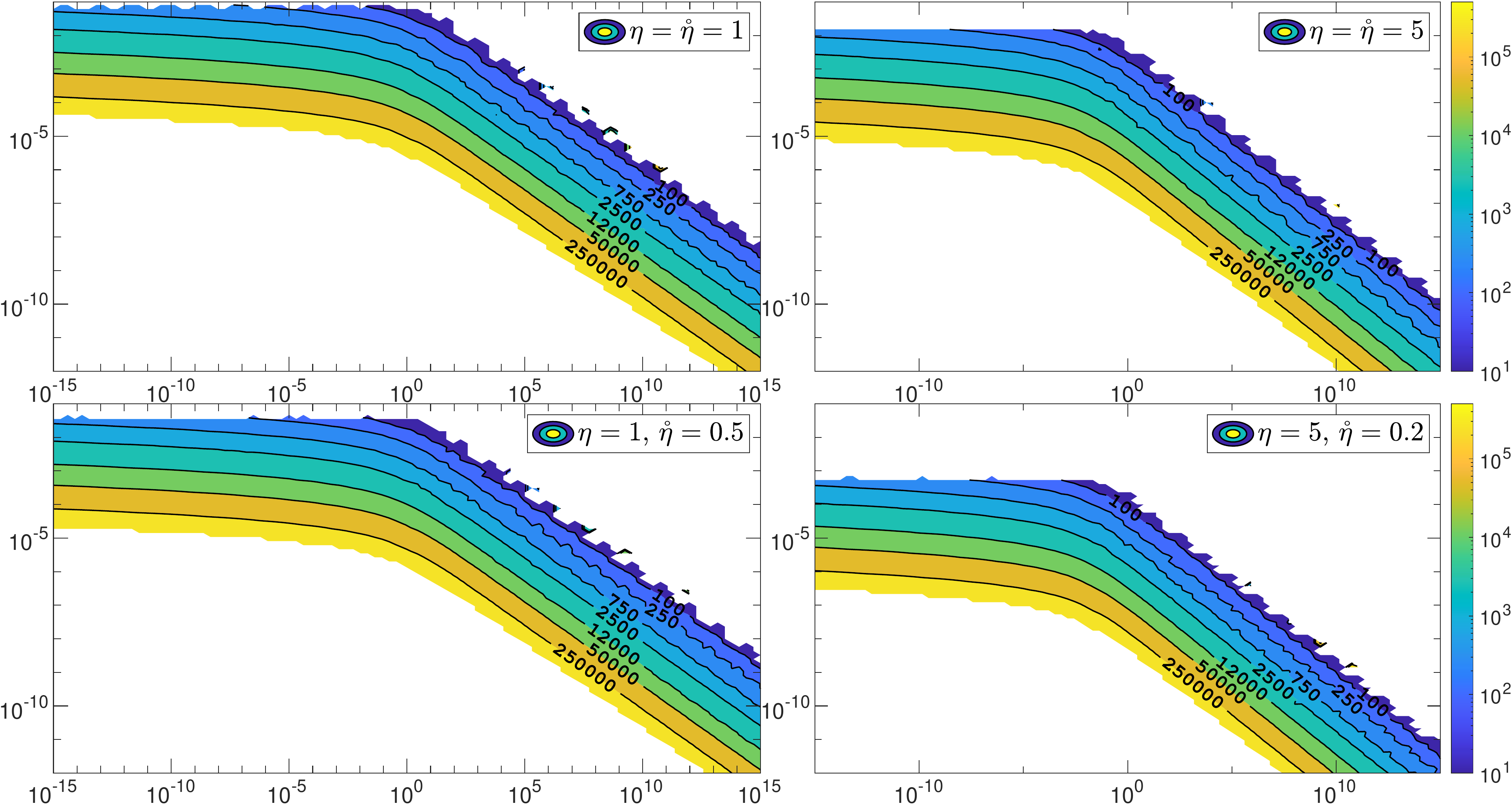} 	\vspace*{-5mm} \caption{ Contour plot of the number of iterations required to achieve convergence ($\delta = 10^{-5}$) in the $(C,h)$-plane, for ExpoHTVI applied to Problem~ \ref{Problem: Polynomial and Log}. } \label{fig: Expo4Param2}  
	\vspace{2mm}
\end{figure}

\begin{figure}[!h] 
	\vspace{4mm}
	\centering
	\includegraphics[width=1\textwidth]{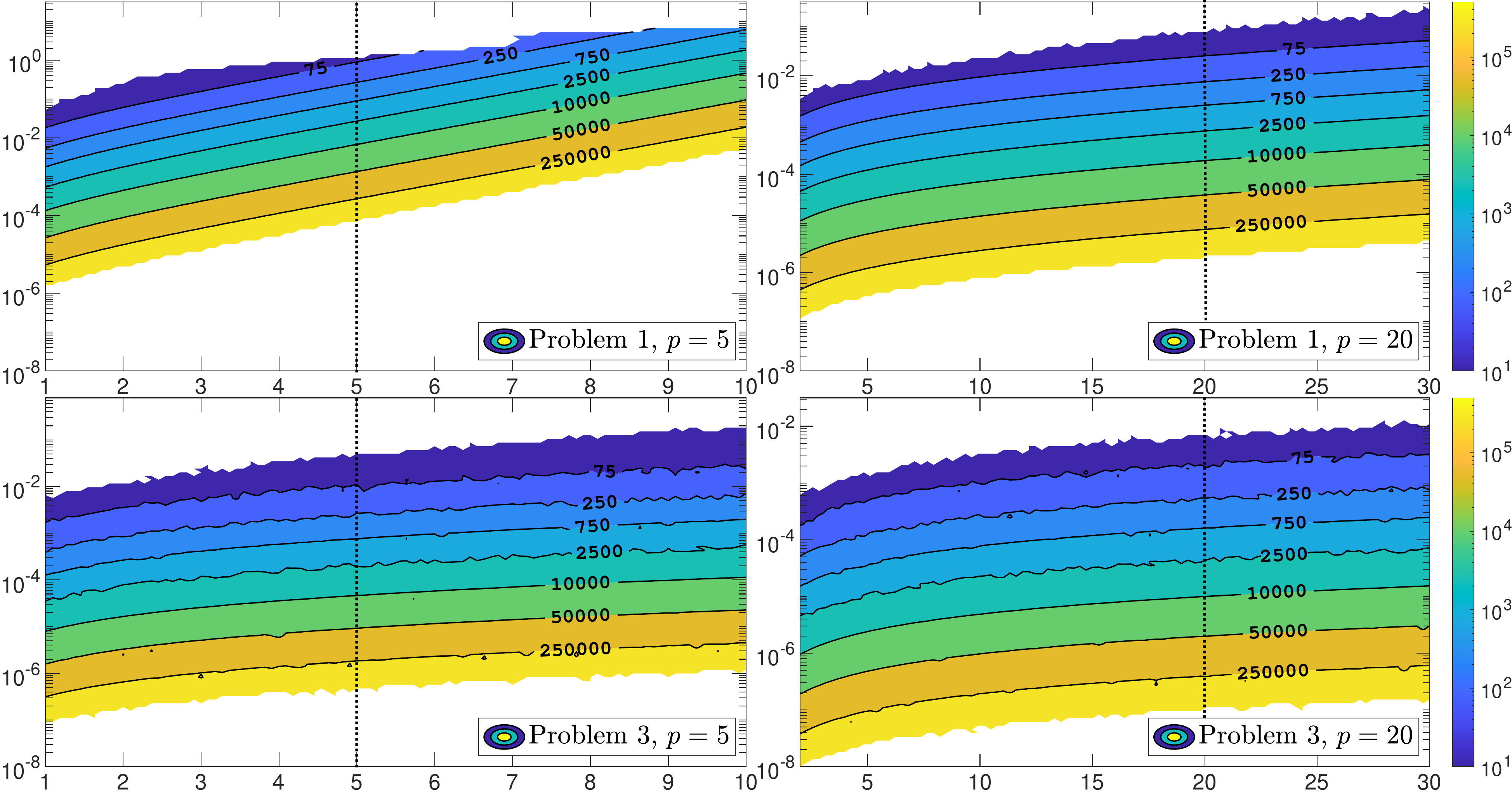} 	\vspace*{-5mm} \caption{Contour plot of the number of iterations required to achieve convergence ($\delta = 10^{-6}$) in the $(\mathring{p},h)$-plane, for PolyHTVI applied to Problems~\ref{Problem: Quartic} (with $C = 10^{-5}$) and \ref{Problem: x log(x)} (with $C = 1$). The dotted line represents the non-adaptive algorithm $p = \mathring{p}$.} \label{fig: AdaptiveTest1} 
\end{figure}

\begin{figure}[!h] 
	\vspace{4mm}
	\centering
	\includegraphics[width=1\textwidth]{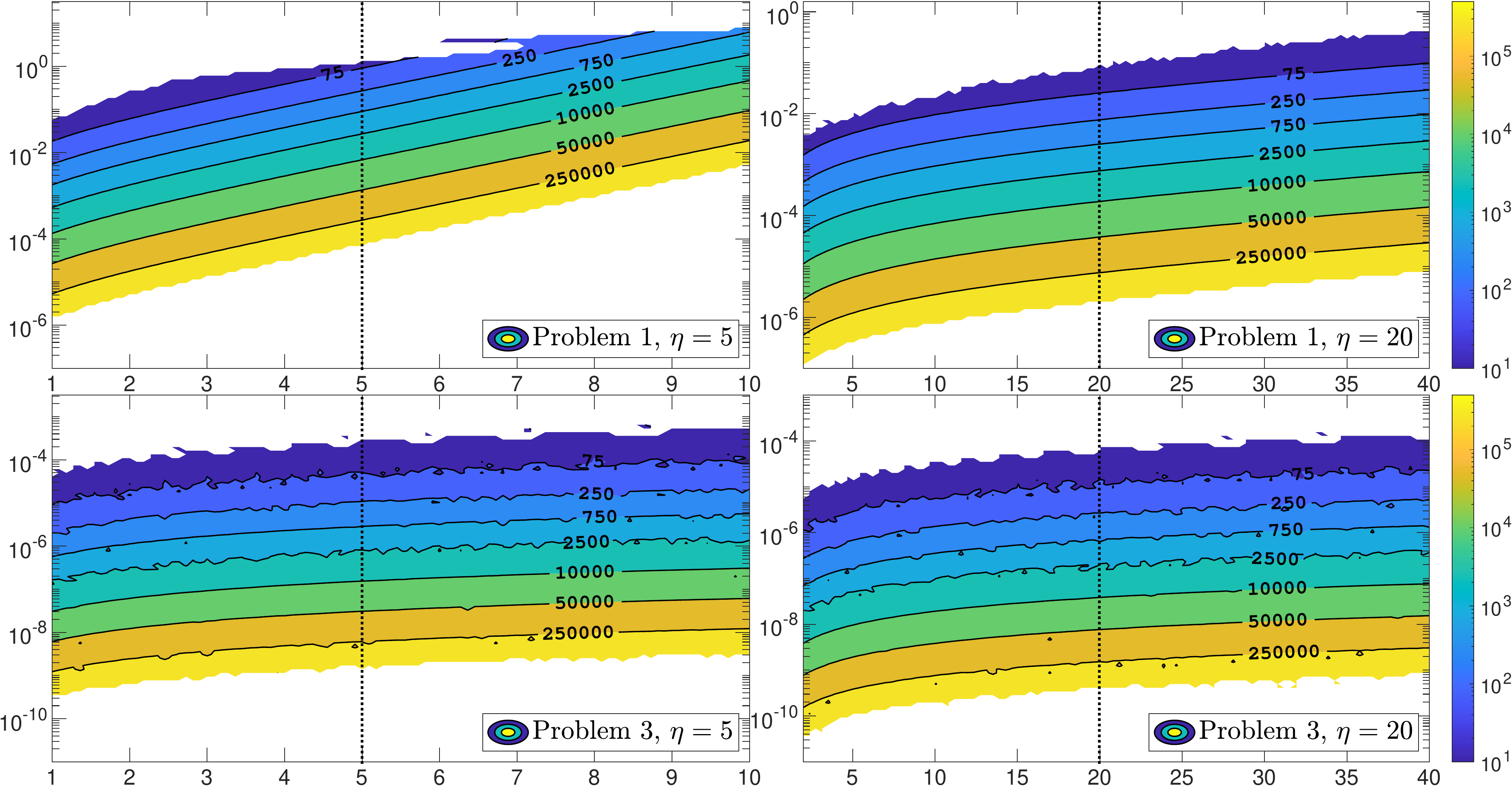} 	\vspace*{-5mm} \caption{Contour plot of the number of iterations required to achieve convergence ($\delta = 10^{-6}$) in the $(\mathring{\eta},h)$-plane, for ExpoHTVI applied to Problems~\ref{Problem: Quartic} (with $C = 10^{-5}$) and \ref{Problem: x log(x)} (with $C = 10^{5}$). The dotted line represents the non-adaptive algorithm $\eta = \mathring{\eta}$.} \label{fig: AdaptiveTest2} 
	\vspace{2mm}
\end{figure}

New numerical experiments carried in this paper suggest that the introduction of momentum restarting schemes in the algorithms enables significantly faster optimization and seems to remove the advantages of the time-adaptive formulation. Indeed, the contour plots in $(C,h)$-space presented in Figures~\ref{fig: Poly4Param1},~\ref{fig: Poly4Param2},~\ref{fig: Poly4Param3},~\ref{fig: Poly4Param4},~\ref{fig: Expo4Param1},~\ref{fig: Expo4Param2} show that the performance and robustness of the PolyHTVI and ExpoHTVI algorithms with momentum restarting is almost unaffected by the introduction of time-adaptivity, regardless of which of Problems~\ref{Problem: Quartic},~\ref{Problem: Polynomial and Log},~\ref{Problem: x log(x)},~\ref{Problem: Ill-Conditioned Quadratic} they are applied to. This is confirmed by the contour plots in $(\mathring{p},h)$-space and $(\mathring{\eta},h)$-space presented in Figures~\ref{fig: AdaptiveTest1} and \ref{fig: AdaptiveTest2} where we can see that for fixed $p$ or $\eta$, the value of $\mathring{p}$ or $\mathring{\eta}$ has very little effect on the performance of the algorithms. 

Overall, the use of time-adaptivity allows for a larger family of algorithms from which one might be able to extract a more efficient algorithm than without time-adaptivity. However, our numerical experiments suggest that with momentum restarting, the benefits time-adaptivity may provide are very limited and are not worth the computational effort of tuning one additional parameter $\mathring{p}$ or $\mathring{\eta}$. For this reason, we will now discard time-adaptivity, and focus on the non-adaptive approaches. More precisely, we will not consider the ExpoToPoly and PolyToExpo Bregman subfamilies anymore, and will only focus on the $p = \mathring{p}$ polynomial and $\eta = \mathring{\eta}$ exponential Bregman subfamilies. 

%\newpage 

%\hfill 

\hfill  \\

%\newpage 

%\newpage 

%\hfill 

\hfill \\
 
%\newpage 

%\newpage 

%\hfill 

\hfill  \\

%\newpage 

%\newpage 

%\hfill 

\hfill  \\

\section{Comparison of Integrators} \label{sec: Comparison Integrators}

Without time-adaptivity, the accelerated optimization algorithms derived in Section~\ref{section: Numerical Methods} can be simplified into the symplectic algorithms presented in Appendix~\ref{appendix: List of Non-Adaptive Algorithms}, and are now all explicit. It should be noted that the HTVI and LTVI algorithms are now equivalent, and will be referred to as LTVI from now on (since the construction of LTVIs extends to Riemannian manifolds, while that of HTVIs does not). Note that we are still using momentum restarting based on the gradient scheme in all our numerical experiments.  \\

\begin{figure}[!hp]
	\centering
	\begin{minipage}[b]{1\textwidth}
		\includegraphics[width=\textwidth]{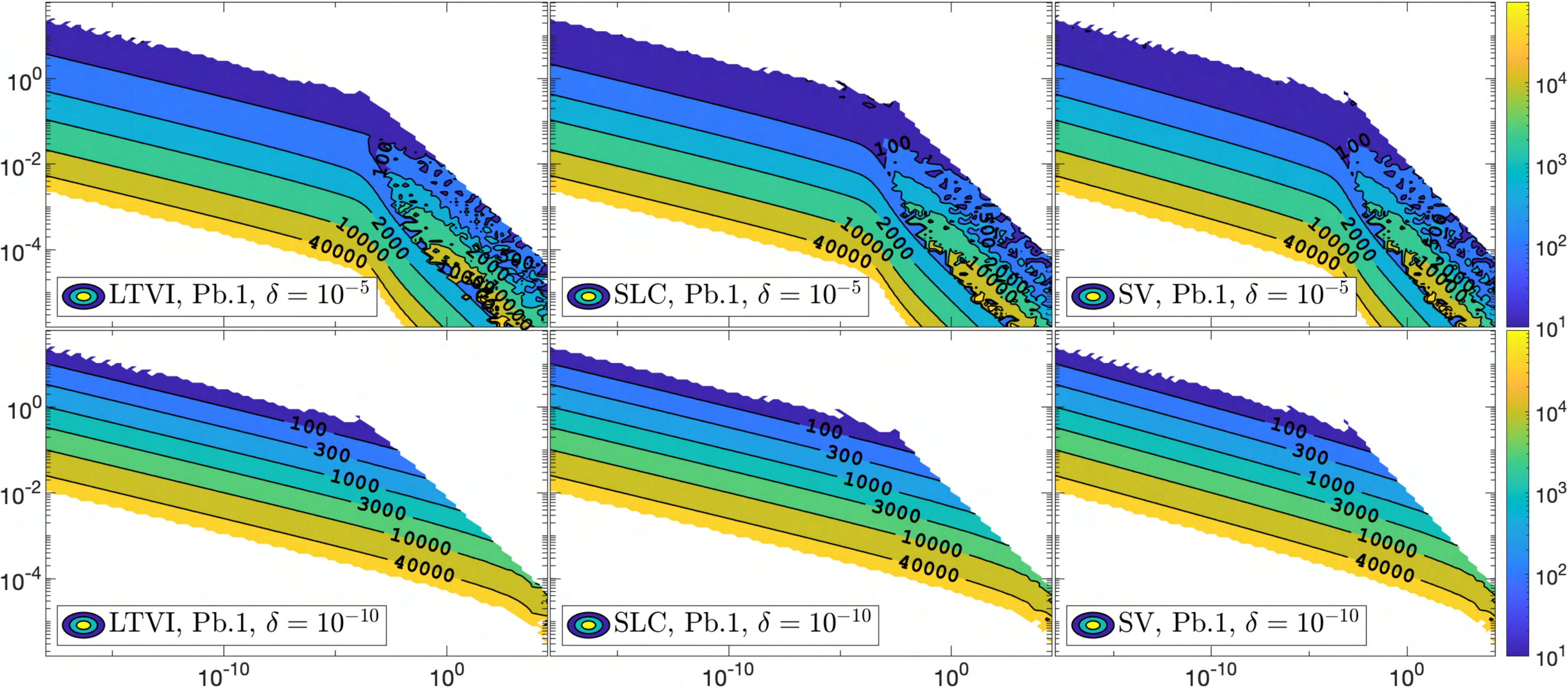} 
	\end{minipage}
	\begin{minipage}[b]{1\textwidth}
		\includegraphics[width=\textwidth]{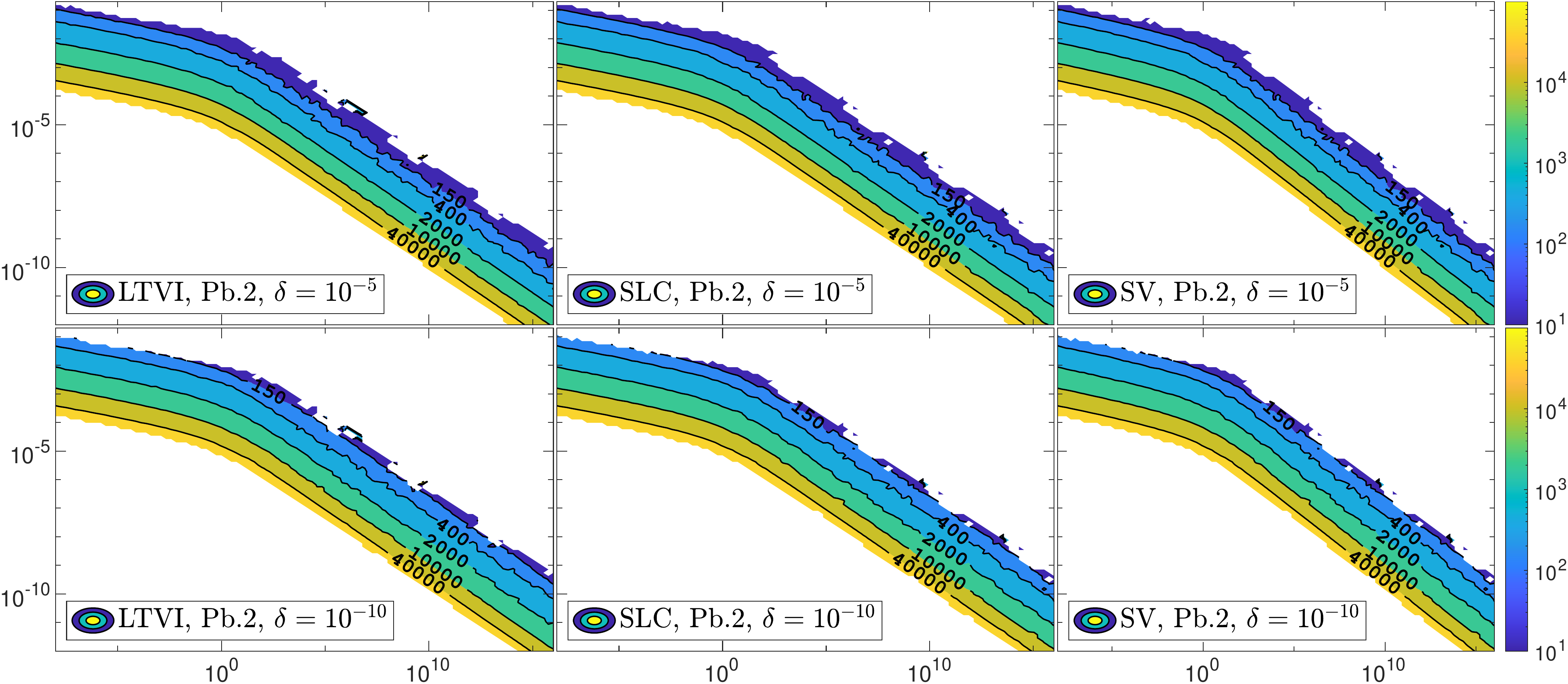} 	
	\end{minipage}
	\begin{minipage}[b]{1\textwidth}
		\includegraphics[width=\textwidth]{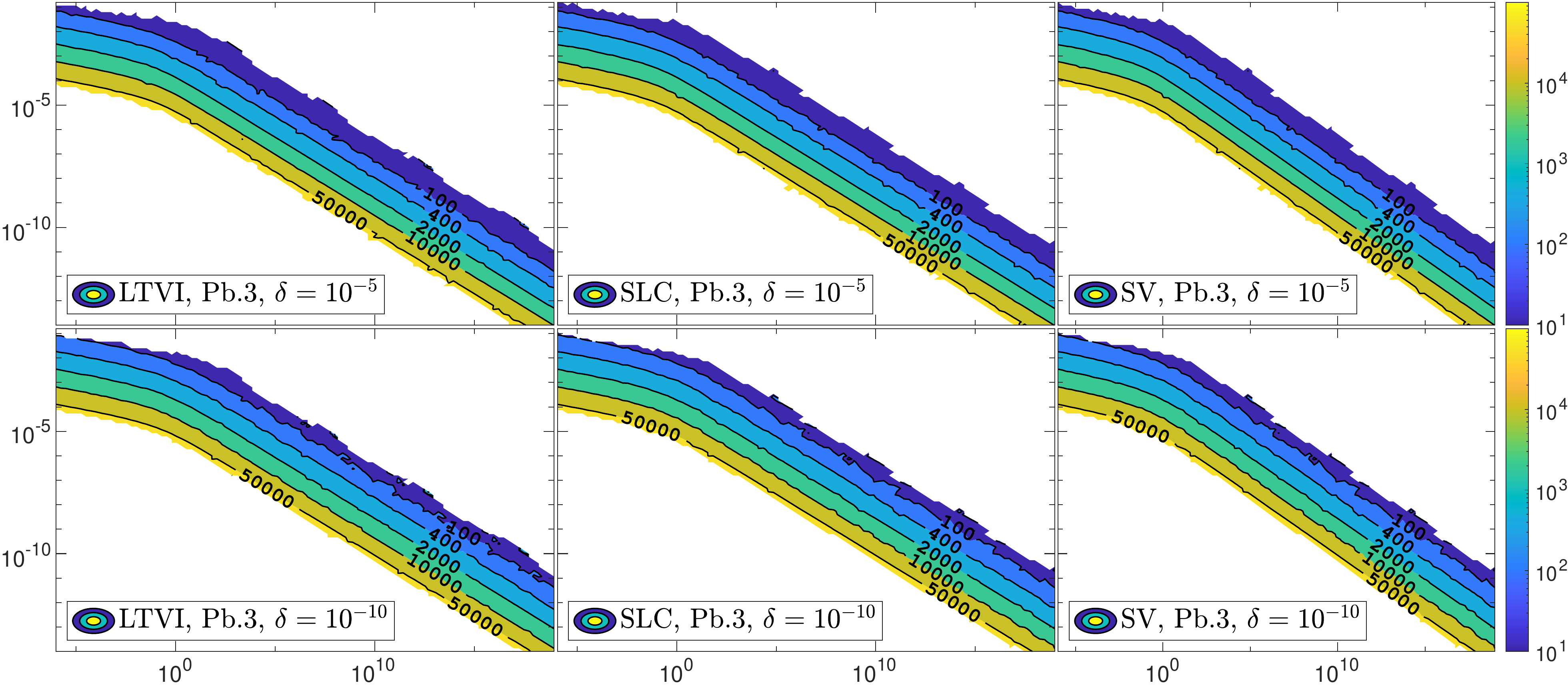}
	\end{minipage} 
	\caption{Convergence regions in the $(C,h)$-plane for PolyLTVI, PolySLC and PolySV applied to Problems~\ref{Problem: Quartic},~\ref{Problem: Polynomial and Log},~\ref{Problem: x log(x)} with $p = 8$.  \label{fig: Ch_Robustness_Poly} }
\end{figure}

\begin{figure}[!h]
	\vspace*{-2mm}
	\centering
	\begin{minipage}[b]{0.925\textwidth}
		\includegraphics[width=\textwidth]{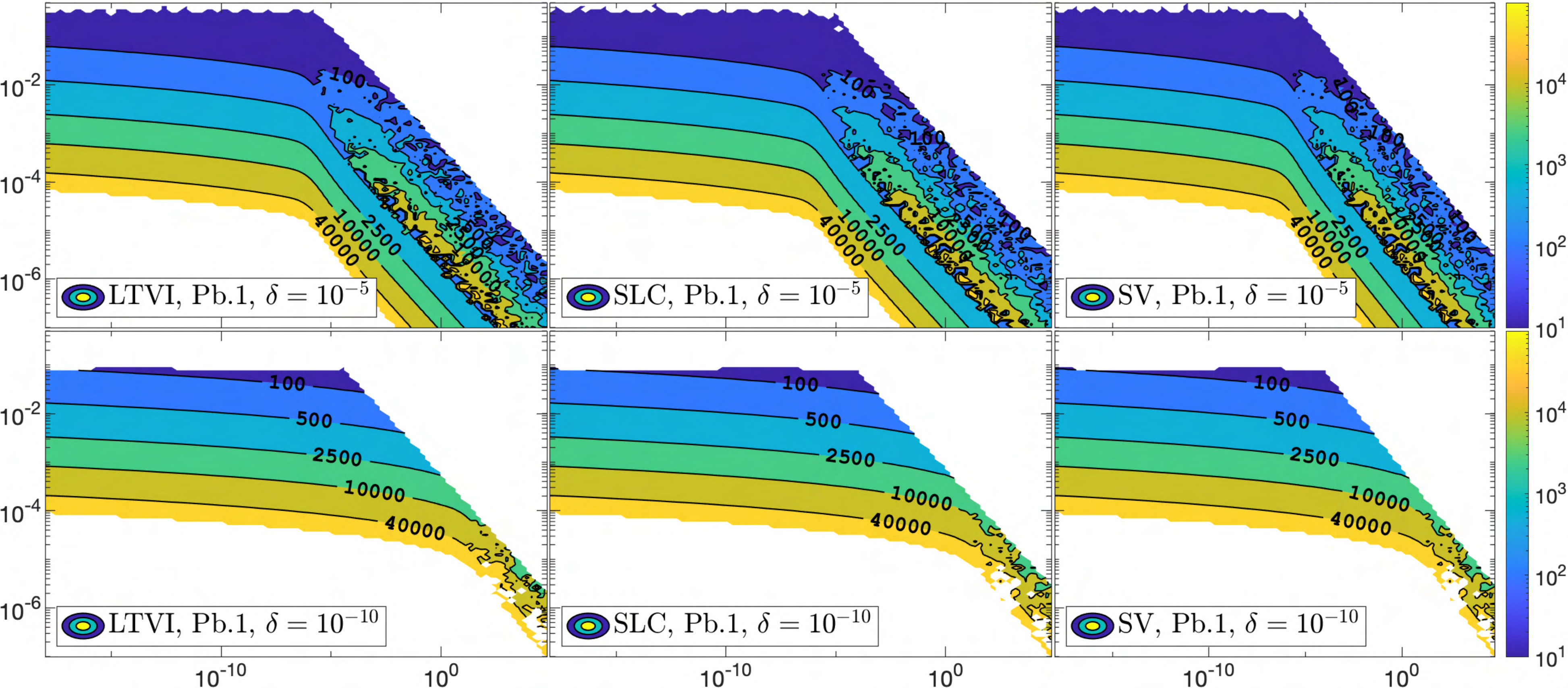} 
	\end{minipage}
	\begin{minipage}[b]{0.925\textwidth}
		\includegraphics[width=\textwidth]{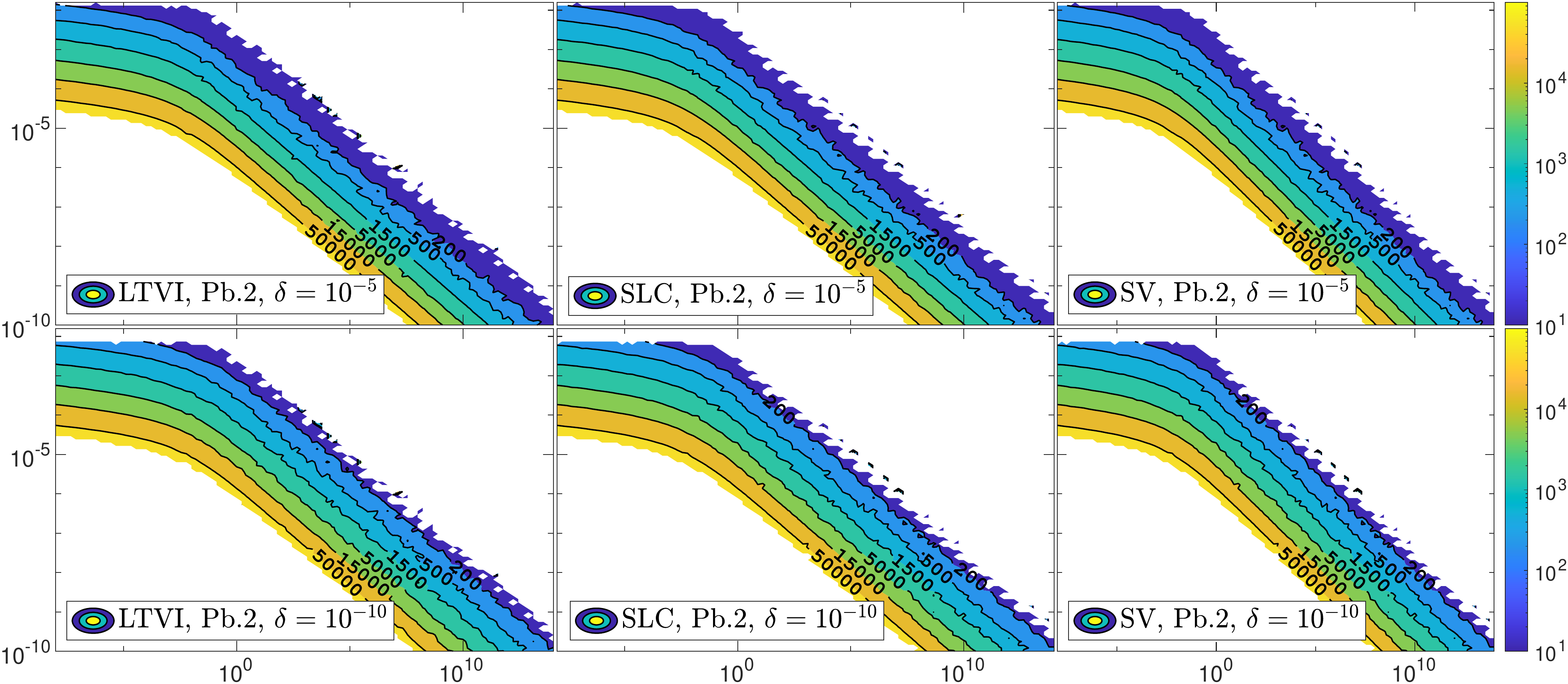} 	\vspace*{-8mm}
	\end{minipage} 
	\caption{Convergence regions in the $(C,h)$-plane for ExpoLTVI, ExpoSLC and ExpoSV applied to Problems~\ref{Problem: Quartic},~\ref{Problem: Polynomial and Log} with $\eta = 6$.  \label{fig: Ch_Robustness_Expo} }
\end{figure}

\begin{figure}[!h]
	\centering
		\includegraphics[width=0.93\textwidth]{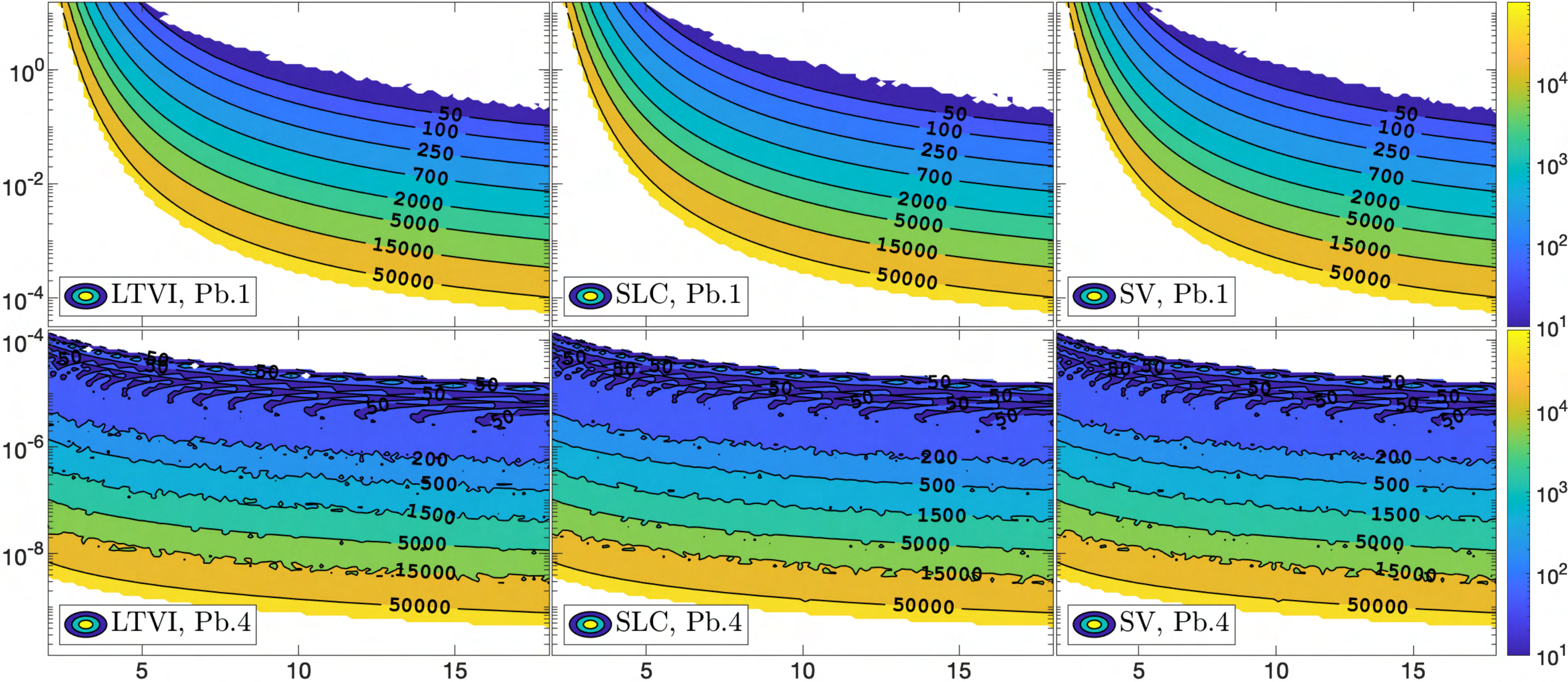}  \vspace*{-3.5mm}
	\caption{Convergence regions in the $(p,h)$-plane for PolyLTVI, PolySLC and PolySV applied to Problems~\ref{Problem: Quartic},~\ref{Problem: Ill-Conditioned Quadratic}.  \label{fig: ph_Robustness} }
\end{figure}

Regions of convergence in the $(C,h)$ and $(p,h)$ planes for the different algorithms were computed based on 100$\times$100 grids of points and are presented in Figures~\ref{fig: Ch_Robustness_Poly},~\ref{fig: Ch_Robustness_Expo},~\ref{fig: ph_Robustness}. We can see that the regions of fast convergence both in the $(C,h)$-plane and $(p,h)$-plane for the different algorithms are almost identical, regardless of the termination criteria. It seems that the three algorithms perform in a very similar way in terms of computational efficiency, robustness, and stability. 

\begin{figure}[!hp]
	\centering
	\begin{minipage}[b]{0.99\textwidth}
		\includegraphics[width=\textwidth]{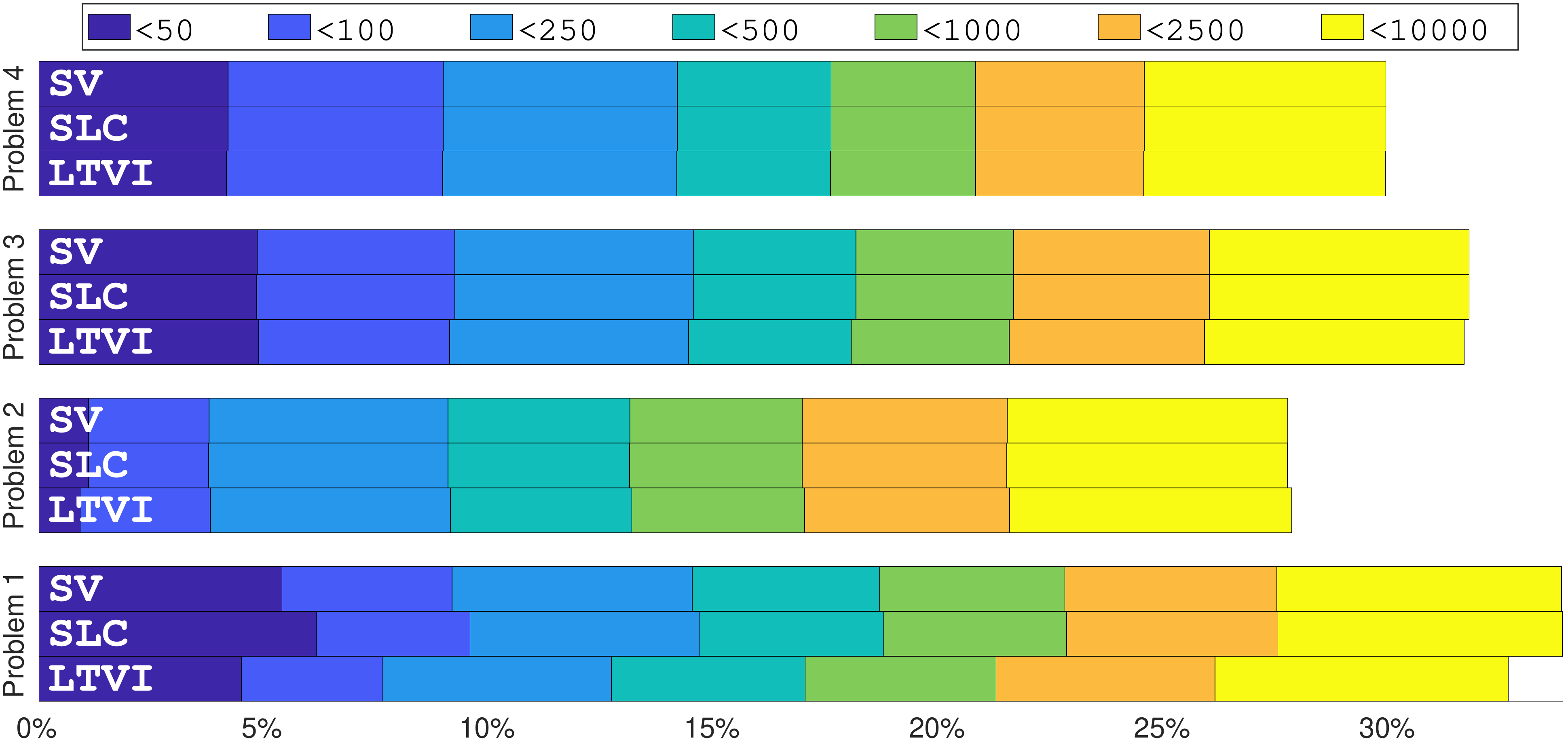} \vspace*{0.5mm}
	\end{minipage}
	\begin{minipage}[b]{0.99\textwidth}
		\includegraphics[width=\textwidth]{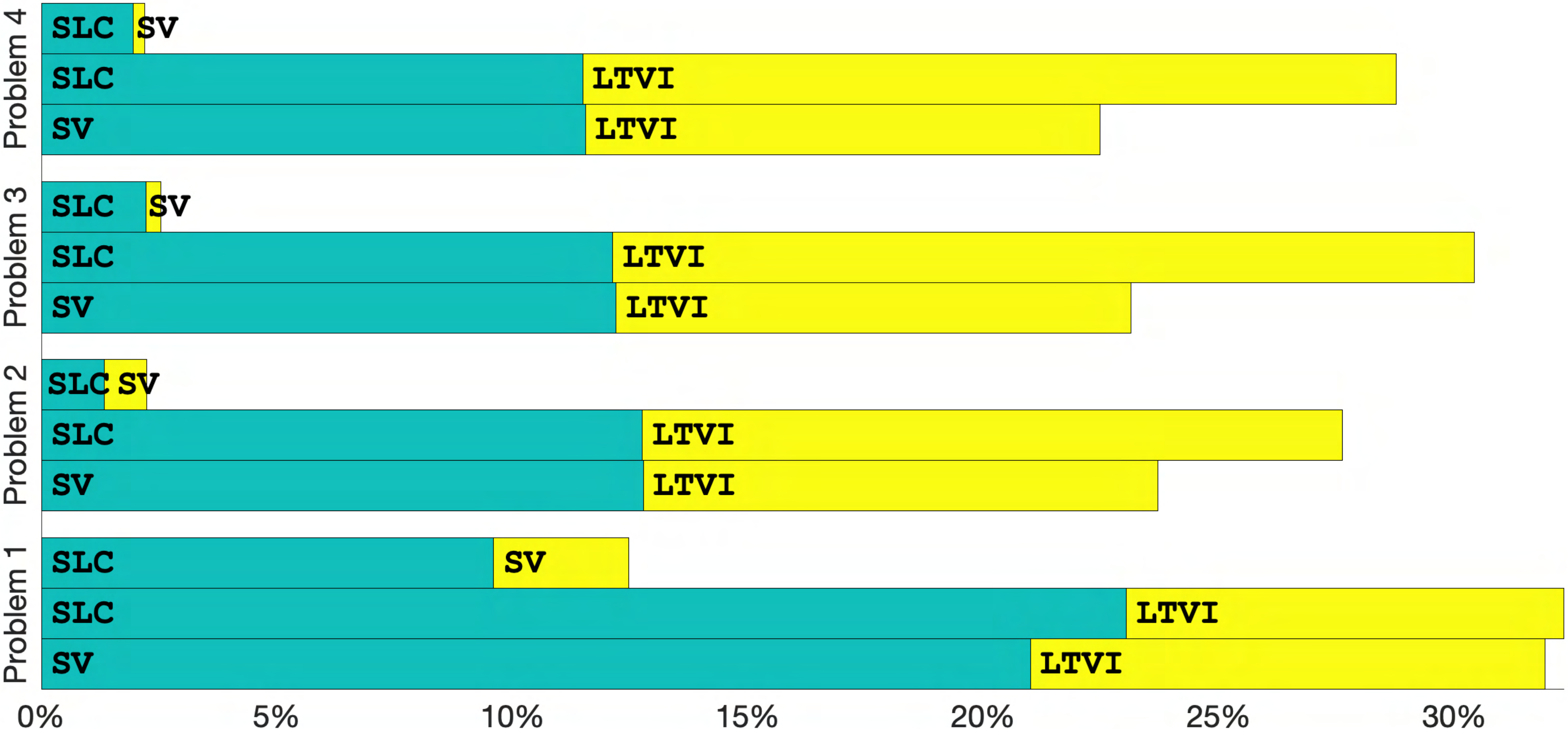} 	\vspace*{0.5mm}
	\end{minipage}
	\begin{minipage}[b]{0.99\textwidth}
		\includegraphics[width=\textwidth]{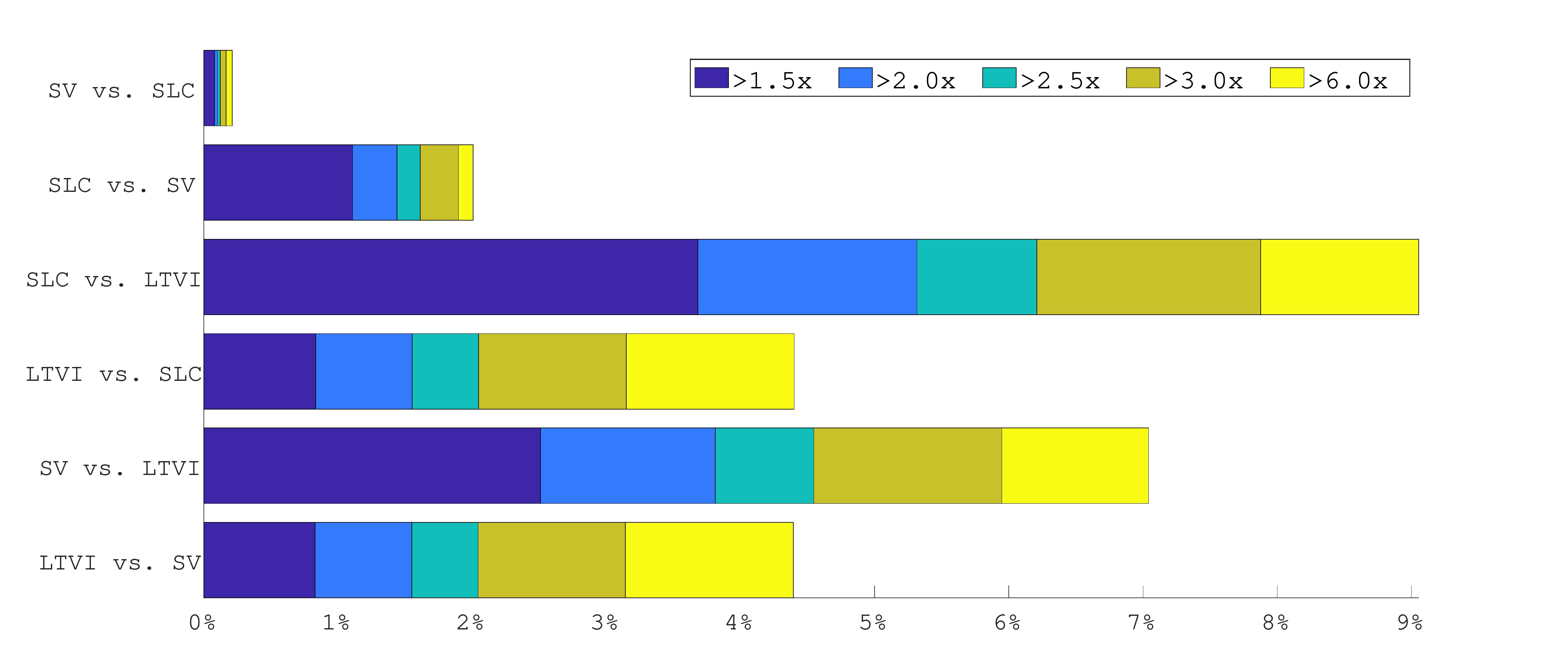}
	\end{minipage} 
	\caption{Results for the polynomial Bregman family with $\delta = 10^{-5}$.  \label{fig: Robustness3D_Poly_delta5} }
\end{figure}

\begin{figure}[!hp]
	\centering
	\begin{minipage}[b]{0.99\textwidth}
		\includegraphics[width=\textwidth]{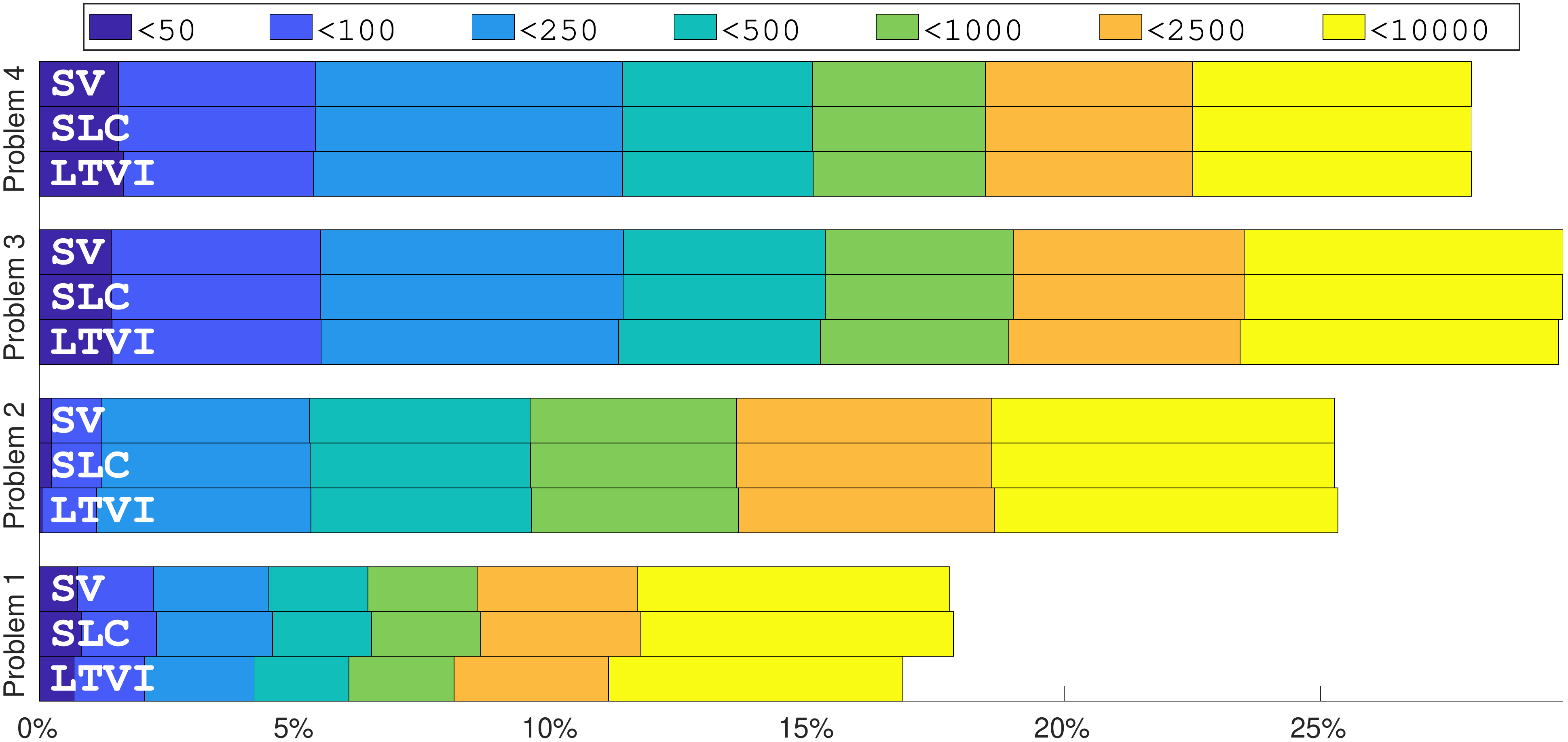} \vspace*{0.5mm}
	\end{minipage}
	\begin{minipage}[b]{0.99\textwidth}
		\includegraphics[width=\textwidth]{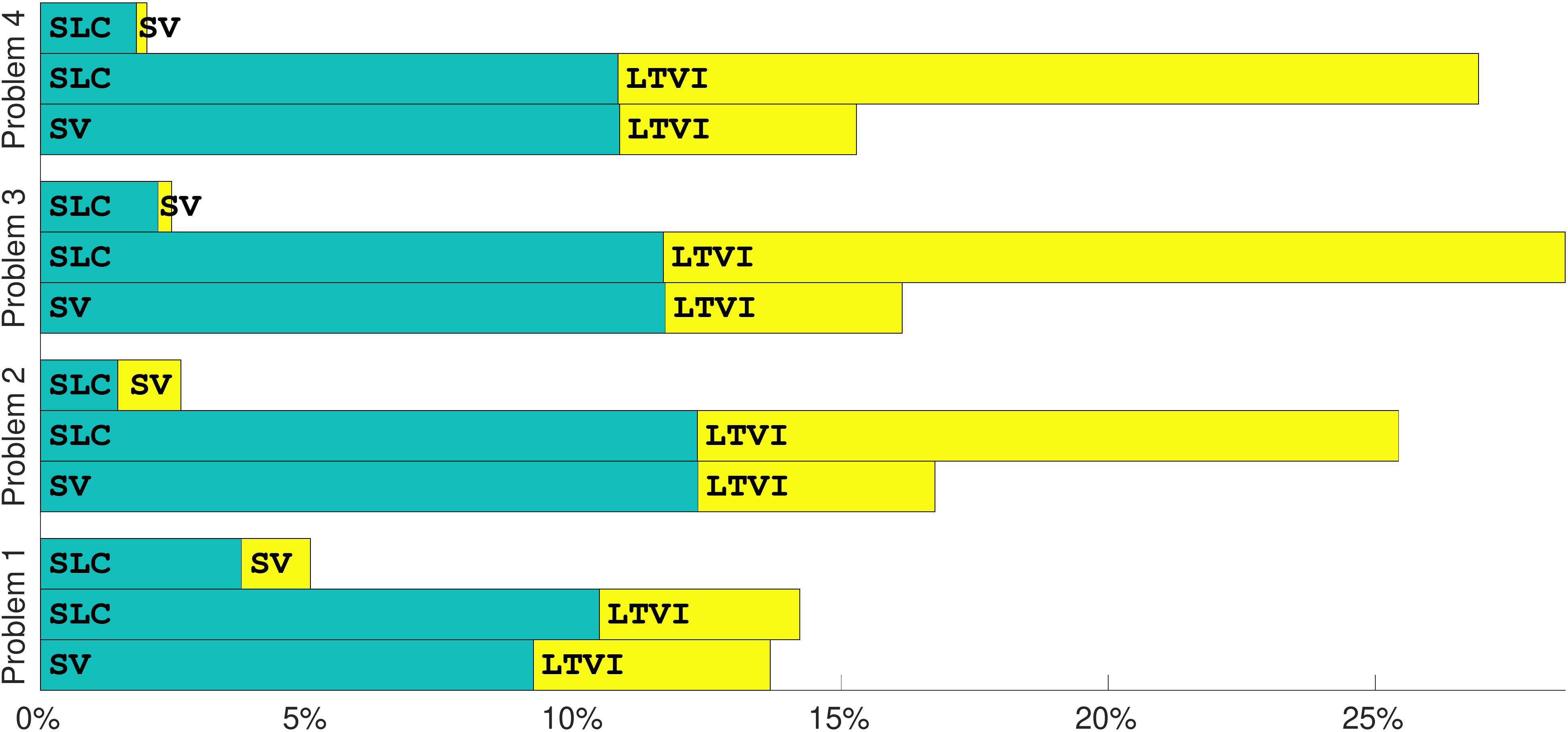} 	\vspace*{0.5mm}
	\end{minipage}
	\begin{minipage}[b]{0.99\textwidth}
		\includegraphics[width=\textwidth]{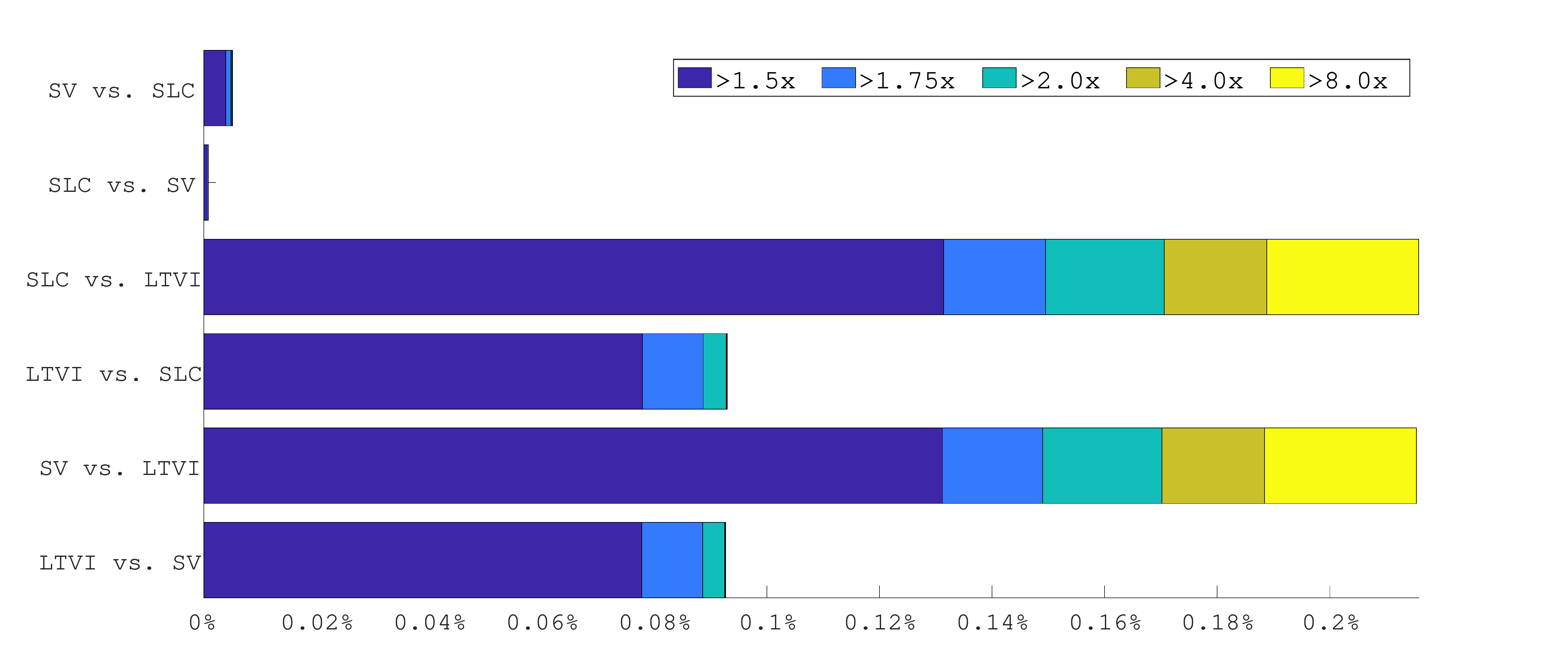}
	\end{minipage} 
	\caption{Results for the polynomial Bregman family with $\delta = 10^{-10}$.  \label{fig: Robustness3D_Poly_delta10} }
\end{figure}

\begin{figure}[!hp]
	\centering
	\begin{minipage}[b]{0.99\textwidth}
		\includegraphics[width=\textwidth]{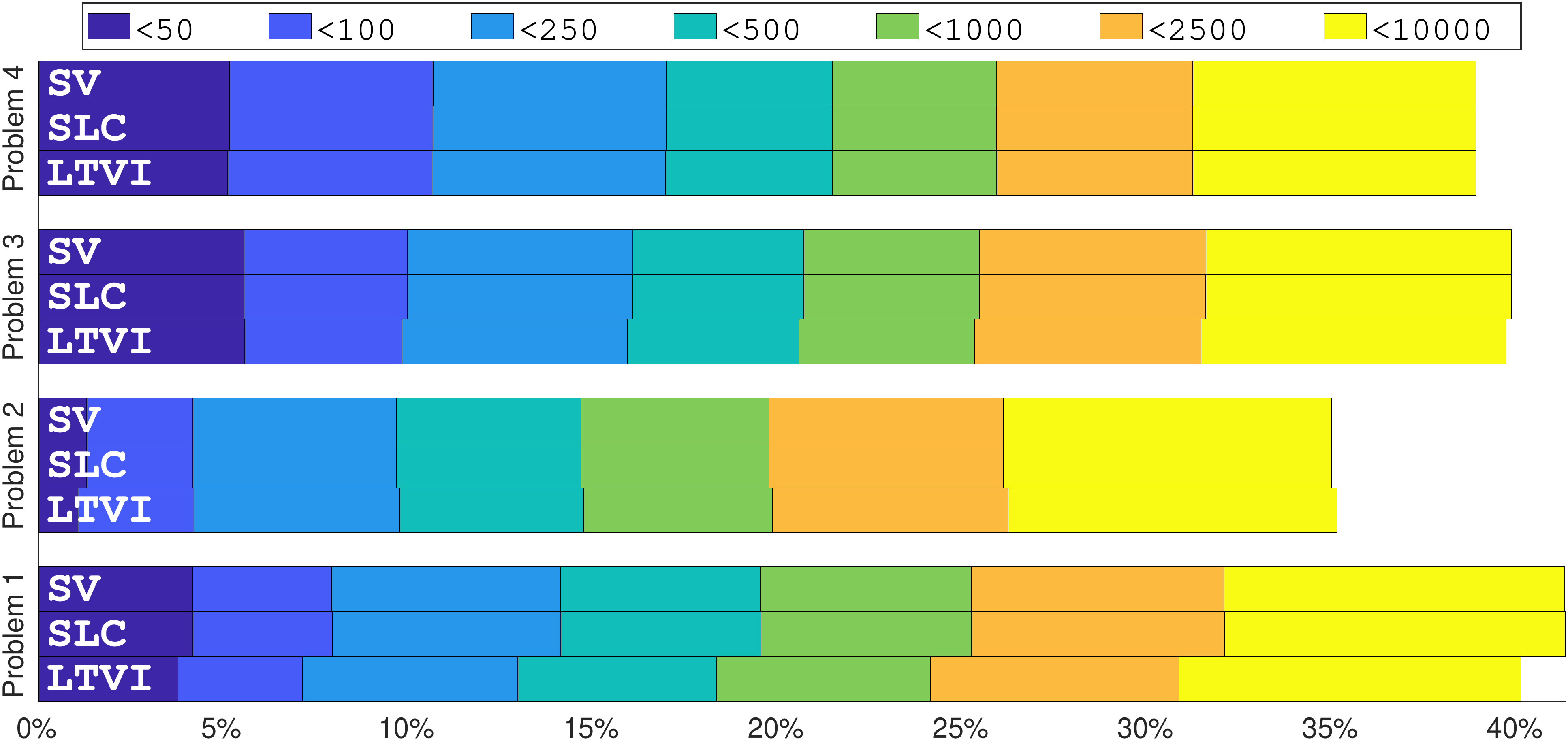} \vspace*{0.5mm}
	\end{minipage}
	\begin{minipage}[b]{0.99\textwidth}
		\includegraphics[width=\textwidth]{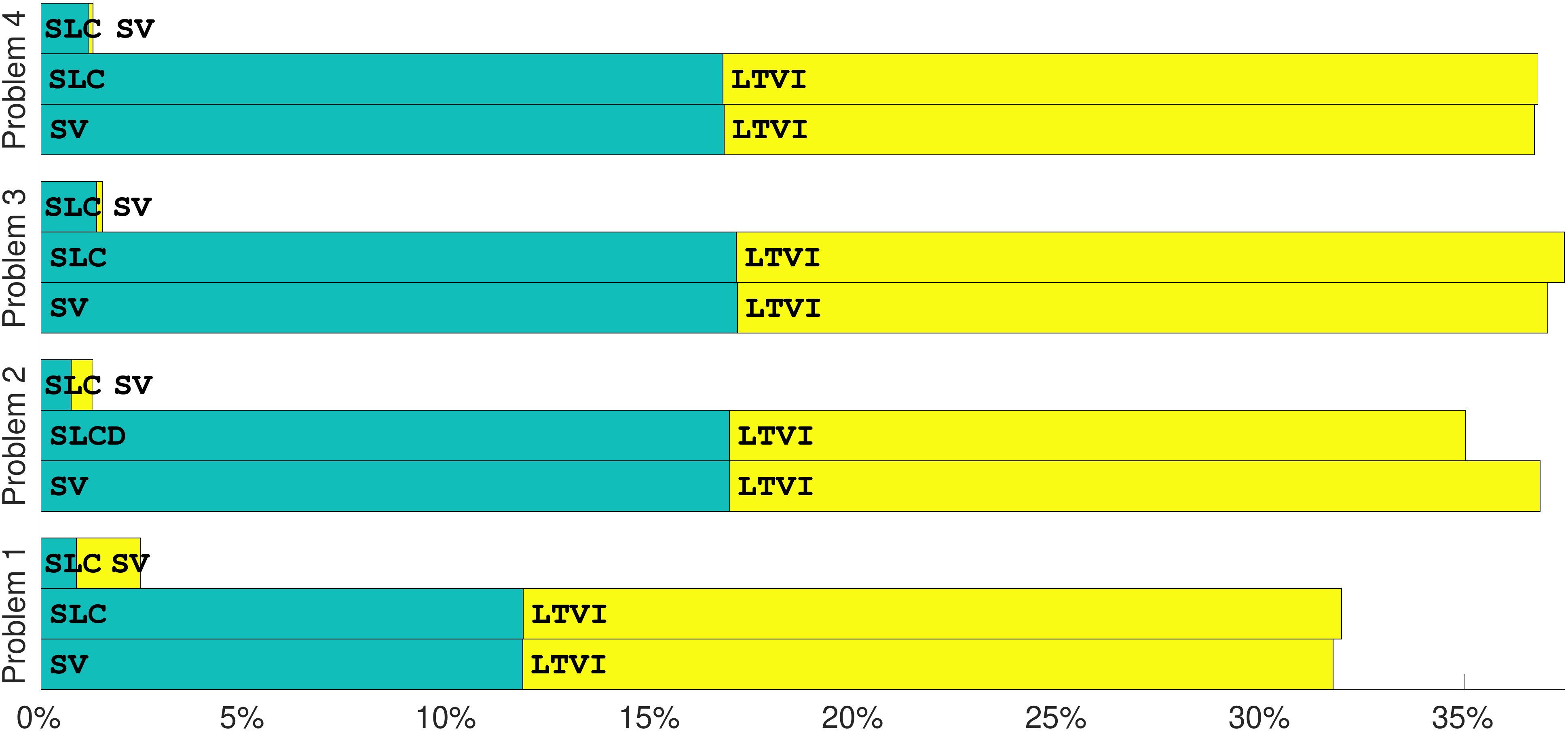} 	\vspace*{0.5mm}
	\end{minipage}
	\begin{minipage}[b]{0.99\textwidth}
		\includegraphics[width=\textwidth]{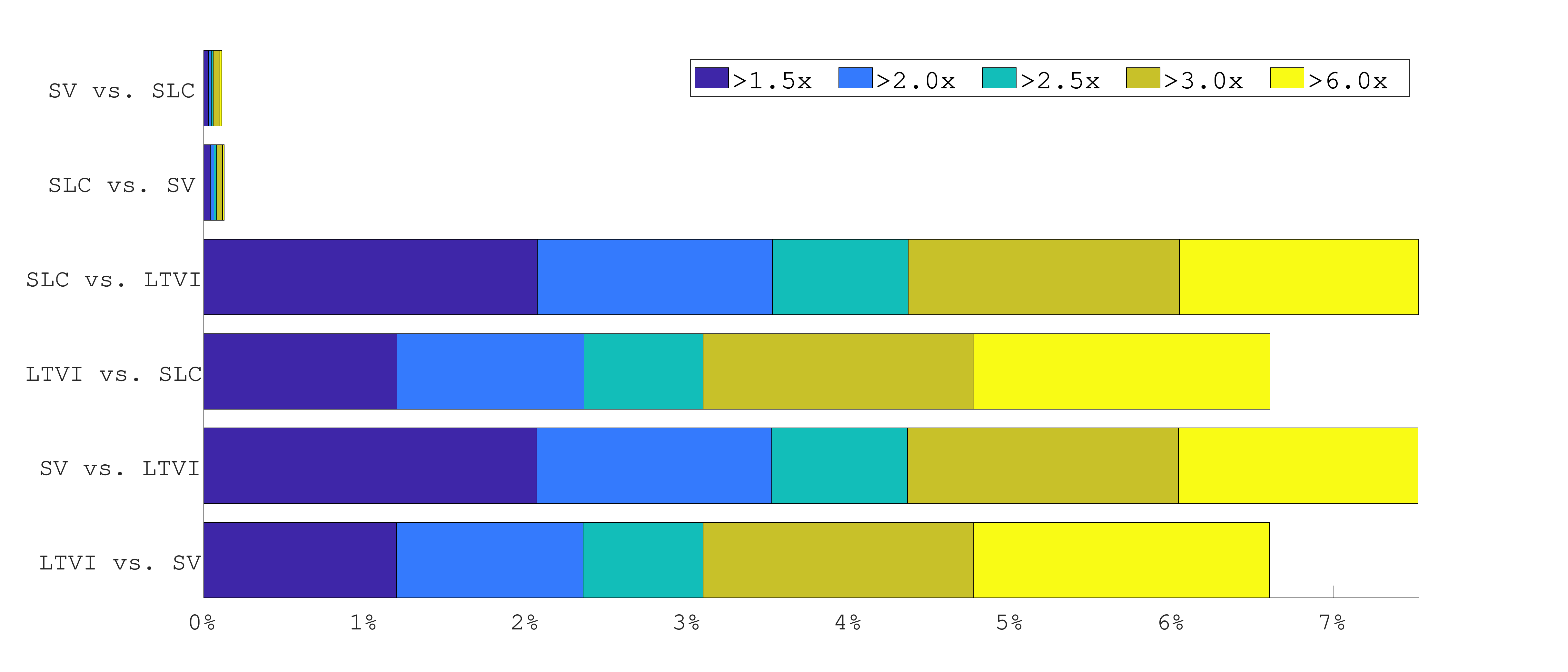}
	\end{minipage} 
	\caption{Results for the exponential Bregman family with $\delta = 10^{-5}$.  \label{fig: Robustness3D_Expo_delta5} }
\end{figure}

\begin{figure}[!hp]
	\centering
	\begin{minipage}[b]{0.99\textwidth}
		\includegraphics[width=\textwidth]{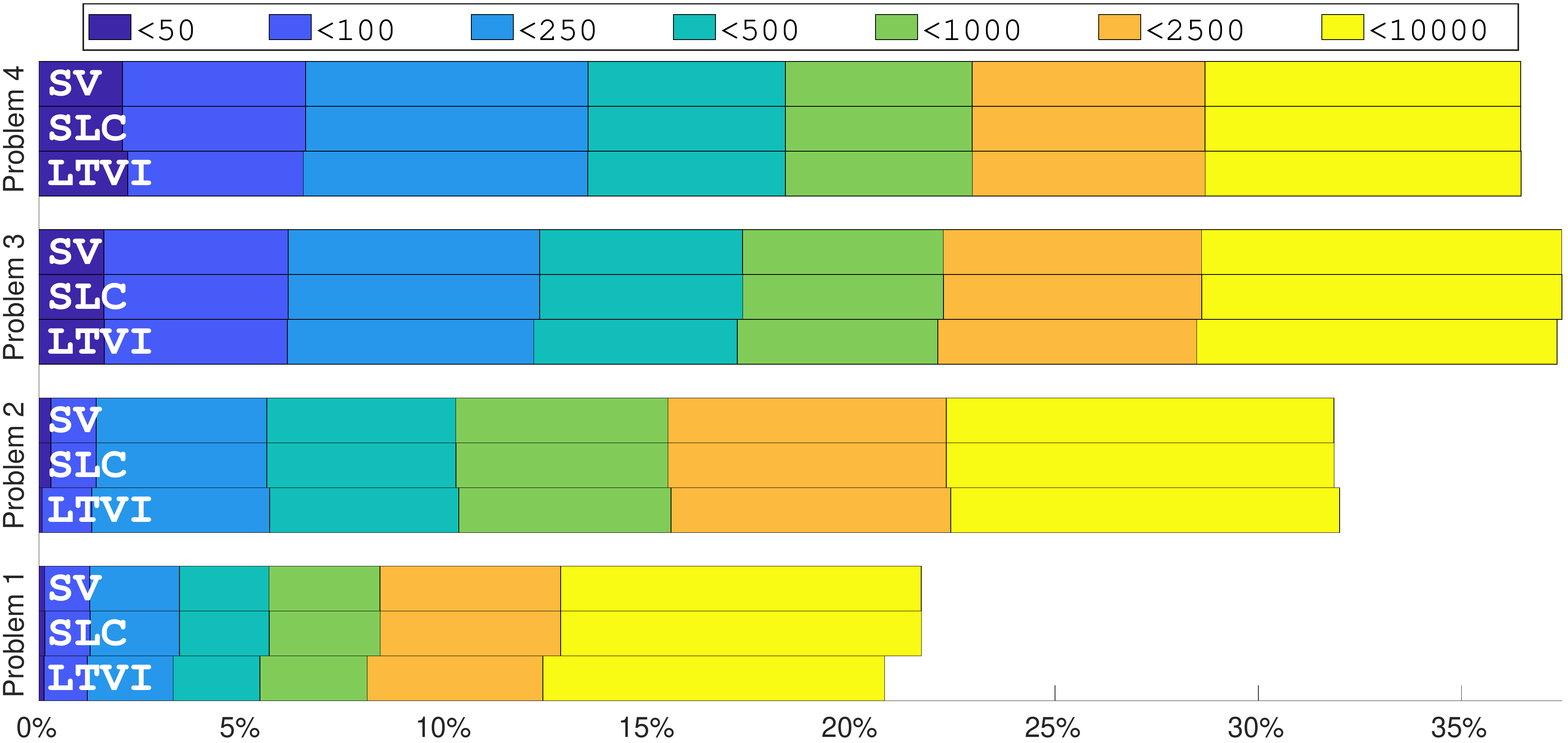} \vspace*{0.5mm}
	\end{minipage}
	\begin{minipage}[b]{0.99\textwidth}
		\includegraphics[width=\textwidth]{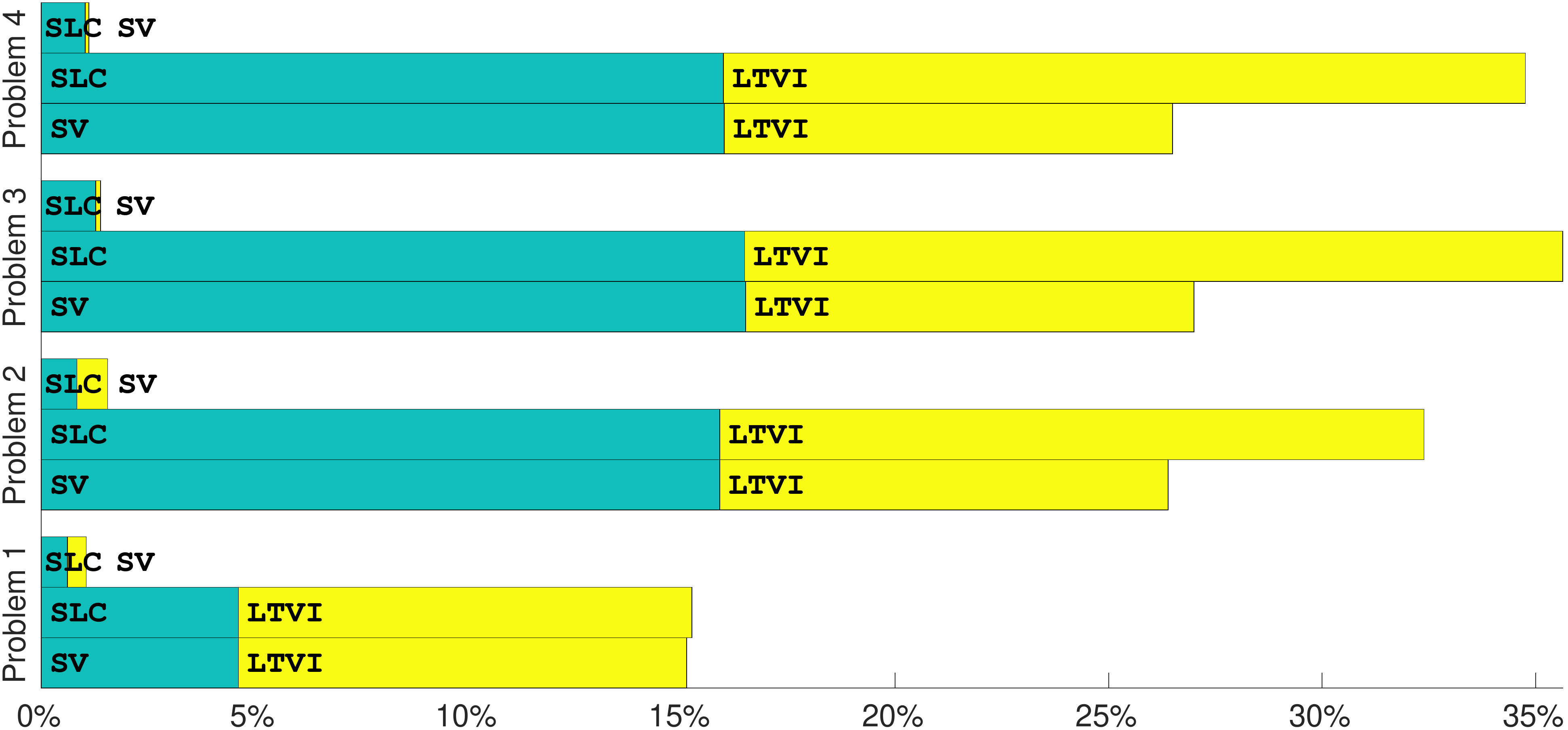} 	\vspace*{0.5mm}
	\end{minipage}
	\begin{minipage}[b]{0.99\textwidth}
		\includegraphics[width=\textwidth]{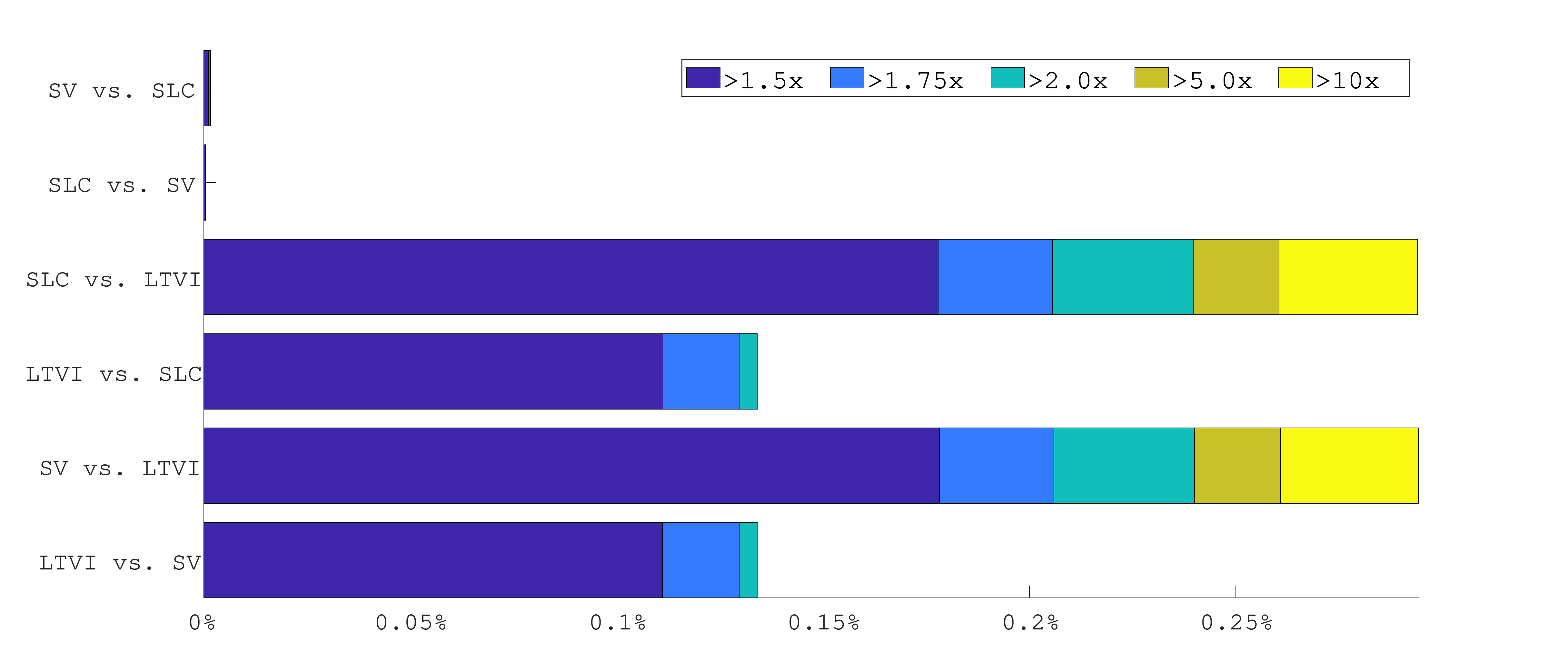}
	\end{minipage} 
	\caption{Results for the exponential Bregman family with $\delta = 10^{-10}$.  \label{fig: Robustness3D_Expo_delta10} }
\end{figure}

We pushed the numerical experimentation further and solved Problems~\ref{Problem: Quartic},~\ref{Problem: Polynomial and Log},~\ref{Problem: x log(x)},~\ref{Problem: Ill-Conditioned Quadratic} on a 3-dimensional grid of $100^3$ points in $(C,p,h)$-space. The results are displayed in Figures~\ref{fig: Robustness3D_Poly_delta5},~\ref{fig: Robustness3D_Poly_delta10},~\ref{fig: Robustness3D_Expo_delta5},~\ref{fig: Robustness3D_Expo_delta10}. 

The top barplots investigate the overall robustness and efficiency of the algorithms. They display the percentages of time each algorithm met the convergence criteria under certain numbers of iterations. For instance, the first bar in the top plot of Figure~\ref{fig: Robustness3D_Poly_delta5} shows that SV converged in $<50$ iterations roughly  5\% of the time, in $<100$ iterations close to 10\% of the time, and so on.  

Each bar in the middle barplots compares the performance of two algorithms for a specific problem, by displaying the percentages of times they outperformed each other. For instance, the last bar in the middle plot of Figure~\ref{fig: Robustness3D_Poly_delta5} shows that SV outperformed LTVI about 21\% of the time, while LTVI outperformed SV roughly 11\% of the time. 

The bottom barplots quantify the gain in efficiency of each algorithm versus the others, by displaying the speedups observed in terms of number of iterations required. For instance, the last bar in the bottom plot of Figure~\ref{fig: Robustness3D_Poly_delta5} shows that LTVI, when compared to SV, achieves a speedup $>1.5\times$ roughly 4.5\% of the time, $>2\times$ roughly 3.5\% of the time, etc... Note that for each bar, we only considered triplets $(C,h,p)$ for which both algorithms converged within $10000$ iterations. \\ 

Let us first focus on the polynomial Bregman family, based on the numerical results displayed in Figures~\ref{fig: Robustness3D_Poly_delta5} and~\ref{fig: Robustness3D_Poly_delta10}. The top plots confirm the earlier observation that the 3 algorithms have very similar regions of fast convergence, and are thus comparable in terms of robustness. The middle and bottom plots indicate that SLC outperforms SV more often than vice-versa, although both scenarios occur rather rarely. Although LTVI seems to outperform SLC/SV roughly as regularly as vice-versa, the speedups LTVI allows when compared to SLC/SV are not as significant and frequent as the slowdowns it entails. It should also be noted that as the termination criteria are made stricter, the differences and significant speedups between the methods become much rarer (although this is partially due to the fact that smaller tolerances require more iterations and we stopped iterating after 10000 iterations). Overall, the 3 algorithms perform very similarly, but the numerical results presented here suggest that SLC might be slightly better within the polynomial Bregman family.

Let us now focus on the exponential Bregman family, based on the numerical results displayed in Figures~\ref{fig: Robustness3D_Expo_delta5} and~\ref{fig: Robustness3D_Expo_delta10}. As was the case for the polynomial family, the 3 algorithms perform very similarly in terms of robustness, and the differences in performance between the algorithms become less significant as the convergence criteria are made stricter. SLC and SV perform almost identically on all problems, regardless of the convergence criteria. Now, it seems that LTVI outperforms SLC/SV slightly more often than vice-versa with the more relaxed tolerances, but with less significant speedups, and as the convergence criteria is made stricter, SLC/SV algorithms seem to outperform LTVI. Overall, the 3 algorithms perform very similarly, but the numerical results suggest that SLC/SV might be the slightly better choices for the exponential Bregman family.


To conclude this section, all the algorithms seem to perform very well and with very small discrepancies, but if we had to choose an algorithm to use to integrate the Bregman
dynamics, it seems that the SLC algorithms with momentum restarting are the slightly better choice. Algorithms~\ref{Alg: PolySLC},~\ref{Alg: ExpoSLC} show more detailed pseudocodes for the SLC algorithms: \\

\begin{algorithm}[H] \label{Alg: PolySLC}
	\DontPrintSemicolon
	
	\KwInput{An objective function $f : \mathbb{R}^d \rightarrow \mathbb{R}$. An initial guess $q\in \mathbb{R}^d$. Parameters $C,h,p>0$. }
	
	%\nonl	\hfill 
	
	$\mathfrak{q} \leftarrow 1, \quad  G \leftarrow \nabla f(q), \quad  r \leftarrow  - \frac{1}{2} Chp \mathfrak{q}^{2p-1 } G$
	
	\While{convergence criteria are not met}
	{		
		$\Delta q \leftarrow   hp \left( \mathfrak{q} + \frac{h}{2} \right)^{-p-1}  r$ \\
		
		$q  \leftarrow q + \Delta q$ \\ 
		
		$G \leftarrow \nabla f(q)$ \\
		
		\textbf{if}  $G^\top \Delta q >0$ \textbf{then} restart momentum: $r  \leftarrow 0$  \\
		
		$\mathfrak{q}  \leftarrow  \mathfrak{q} + h$ \\ 
		
		$r  \leftarrow  r - Chp \mathfrak{q}^{2p-1 } G$\\
	}

	\caption{\textbf{S}ymmetric \textbf{L}eapfrog \textbf{C}omposition of Component Dynamics for the \textbf{Poly}nomial Bregman dynamics, with Momentum \textbf{R}estarting (\textbf{PolySLC-R}) }
\end{algorithm} 

\hfill

\begin{algorithm}[H] \label{Alg: ExpoSLC}
	\DontPrintSemicolon
	
	\KwInput{An objective function $f : \mathbb{R}^d \rightarrow \mathbb{R}$. An initial guess $q\in \mathbb{R}^d$. Parameters $C,h,\eta >0$. }
	
	%\nonl	\hfill 
	
	$\mathfrak{q} \leftarrow 1, \quad  G \leftarrow \nabla f(q), \quad  r \leftarrow -  \frac{1}{2} C\eta h   e^{2\eta \mathfrak{q}} G $
	
	\While{convergence criteria are not met}
	{		
		$\Delta q \leftarrow   	\eta h e^{-\eta \left( \mathfrak{q} + \frac{h}{2} \right)}  r  $ \\
		
		$q  \leftarrow q + \Delta q$ \\ 
		
		$G \leftarrow \nabla f(q)$ \\
		
		\textbf{if}  $G^\top \Delta q >0$ \textbf{then} restart momentum: $r  \leftarrow 0$  \\
		
		$\mathfrak{q}  \leftarrow  \mathfrak{q} + h$ \\ 
		
		$r  \leftarrow  r -  C\eta h   e^{2\eta \mathfrak{q}} G $\\
	} 
	
	\caption{\textbf{S}ymmetric \textbf{L}eapfrog \textbf{C}omposition of Component Dynamics for the \textbf{Expo}nential Bregman dynamics, with Momentum \textbf{R}estarting (\textbf{ExpoSLC-R}) }
\end{algorithm}

\hfill \\

\section{Tuning the Algorithms} \label{section: Tuning the Algorithms}

We will now investigate how the PolySLC-R and ExpoSLC-R algorithms perform as the parameters $C,h,p,\eta $ are varied, and try to reduce the number of parameters needing tuning in practice. 

\subsection{Tuning PolySLC-R}

\hfill \\

\begin{figure}[!h]
	\vspace*{-2mm}
	\centering
	\includegraphics[width=0.91\textwidth]{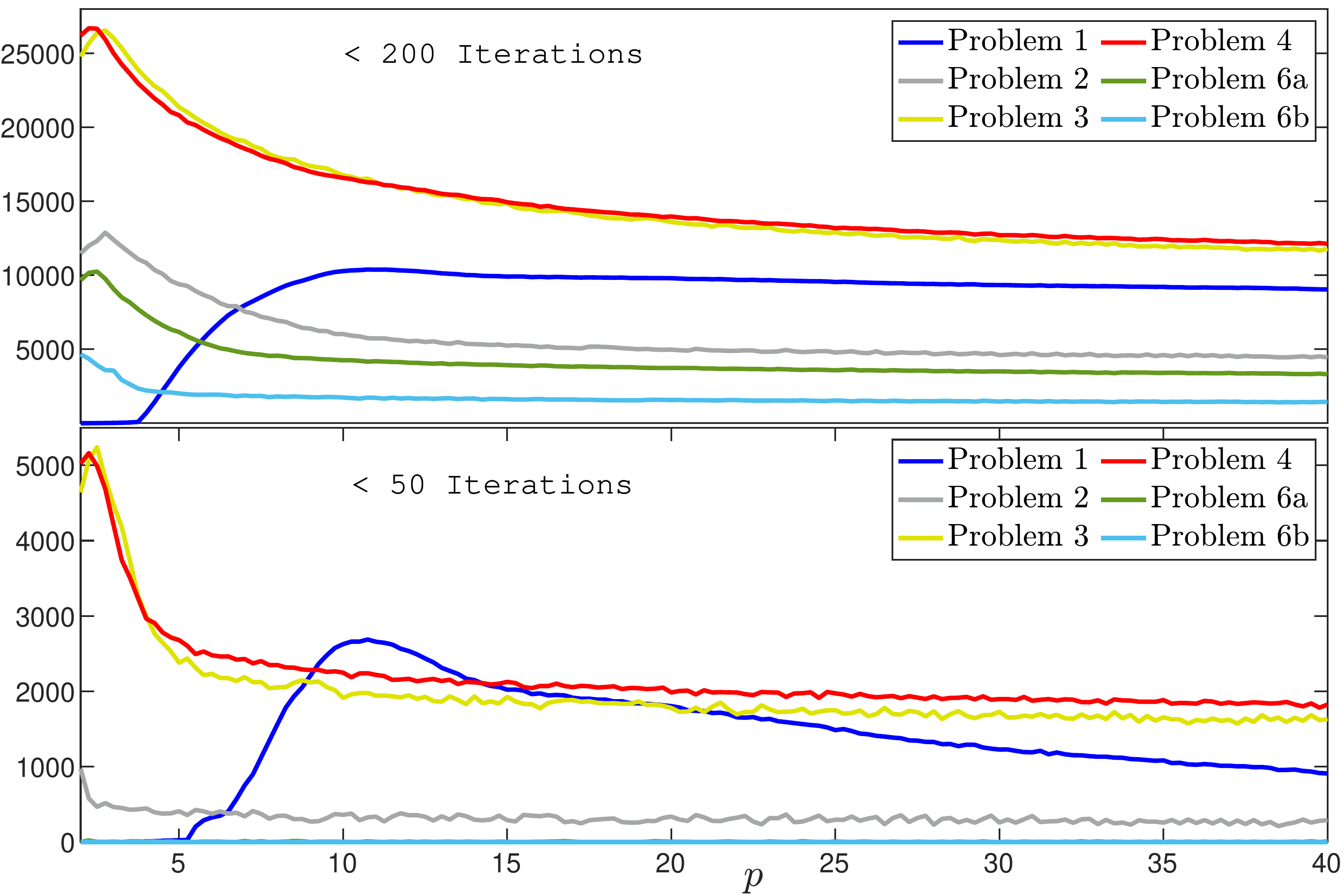}  \vspace*{-4mm}
	\caption{Number of $(C,h)$ pairs (out of $500^2$) for which convergence was achieved under 200 and 50 iterations using PolySLC-R, as the value of $p$ is varied.  \label{fig: PolySLC_p} }
\end{figure}

We solved Problems~\ref{Problem: Quartic},~\ref{Problem: Polynomial and Log},~\ref{Problem: x log(x)},~\ref{Problem: Ill-Conditioned Quadratic}  and two distinct randomly generated instances of Problem~\ref{Problem: Linear Regression} using PolySLC-R on a 3-dimensional grid of $500\times153\times500$ points in $(C,p,h)$-space (logarithmically-spaced in $C$ between $10^{-12}$ and $10^{12}$, logarithmically-spaced in $h$ between $10^{-6}$ and $10^3$, and linearly-spaced in $p$ between $2$ and $40$), and recorded the number of iterations needed to achieve convergence with~$\delta = 10^{-10}$. Figure~\ref{fig: PolySLC_p} displays the number of $(C,h)$ pairs for which convergence was achieved under 200 and 50 iterations for each value of $p$.  We can see that the value of $p$ does not seem to significantly affect the number of $(C,h)$ pairs that exhibit fast convergence, once it is taken to be sufficiently large, so tuning the parameter $p$ carefully might not be very helpful and necessary. For numerical stability reasons, which will be discussed in Section~\ref{section: Temporal Looping to Improve Numerical Stability}, it is preferable to use lower values of $p$, so we will set $p= 6$ since this is a small value of $p$ which performed very well in Figure~\ref{fig: PolySLC_p}. % \\

\begin{figure}[!h]
	\centering
	\includegraphics[width=0.91\textwidth]{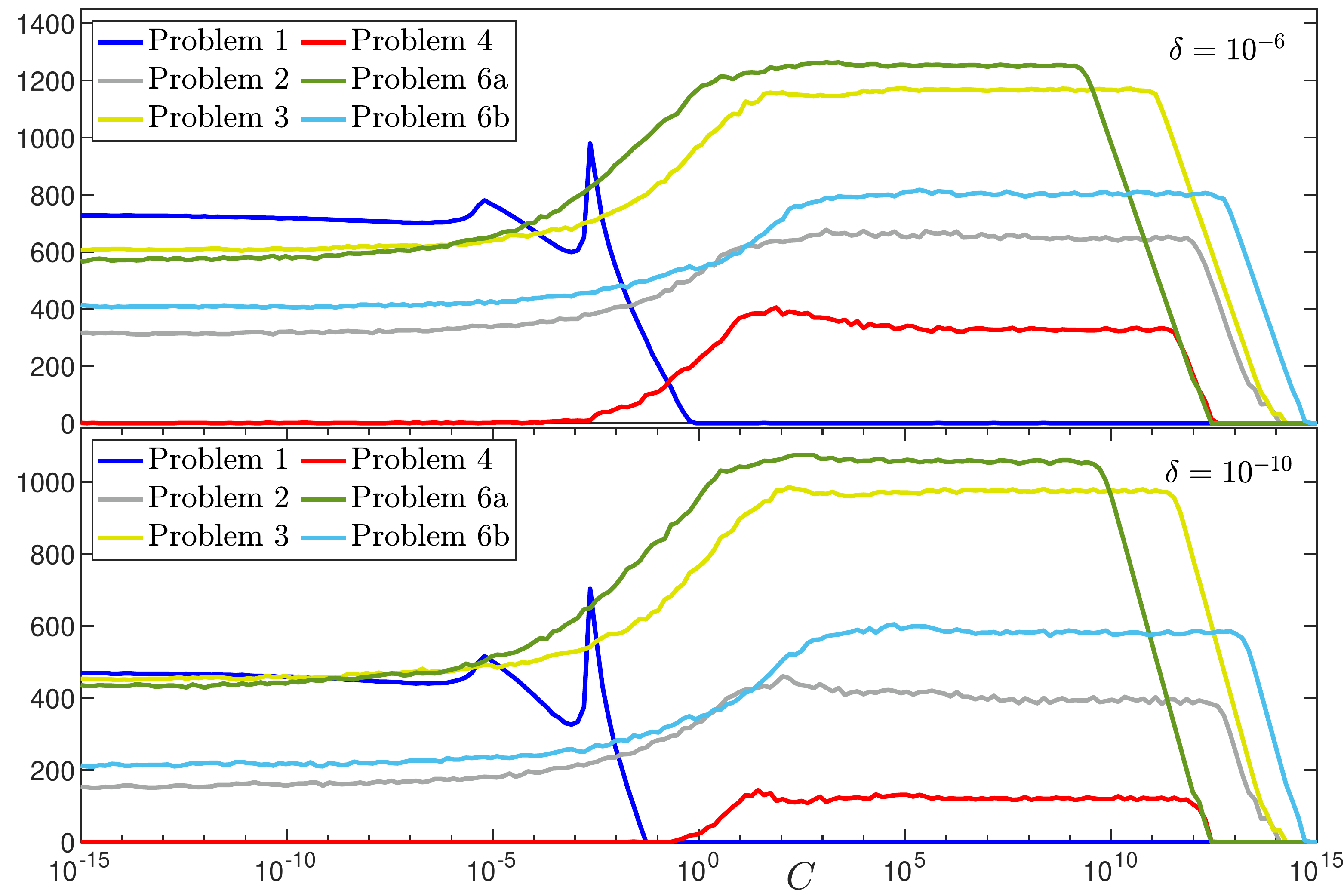}  \vspace*{-4mm}
	\caption{Number of values of $h$ (out of $10^4$) for which convergence was achieved under 200 iterations using PolySLC-R with $p=6$, as the value of $C$ is varied.  \label{fig: PolySLC_2D} }
\end{figure}

\begin{figure}[!h]
	\centering
	\includegraphics[width=0.9\textwidth]{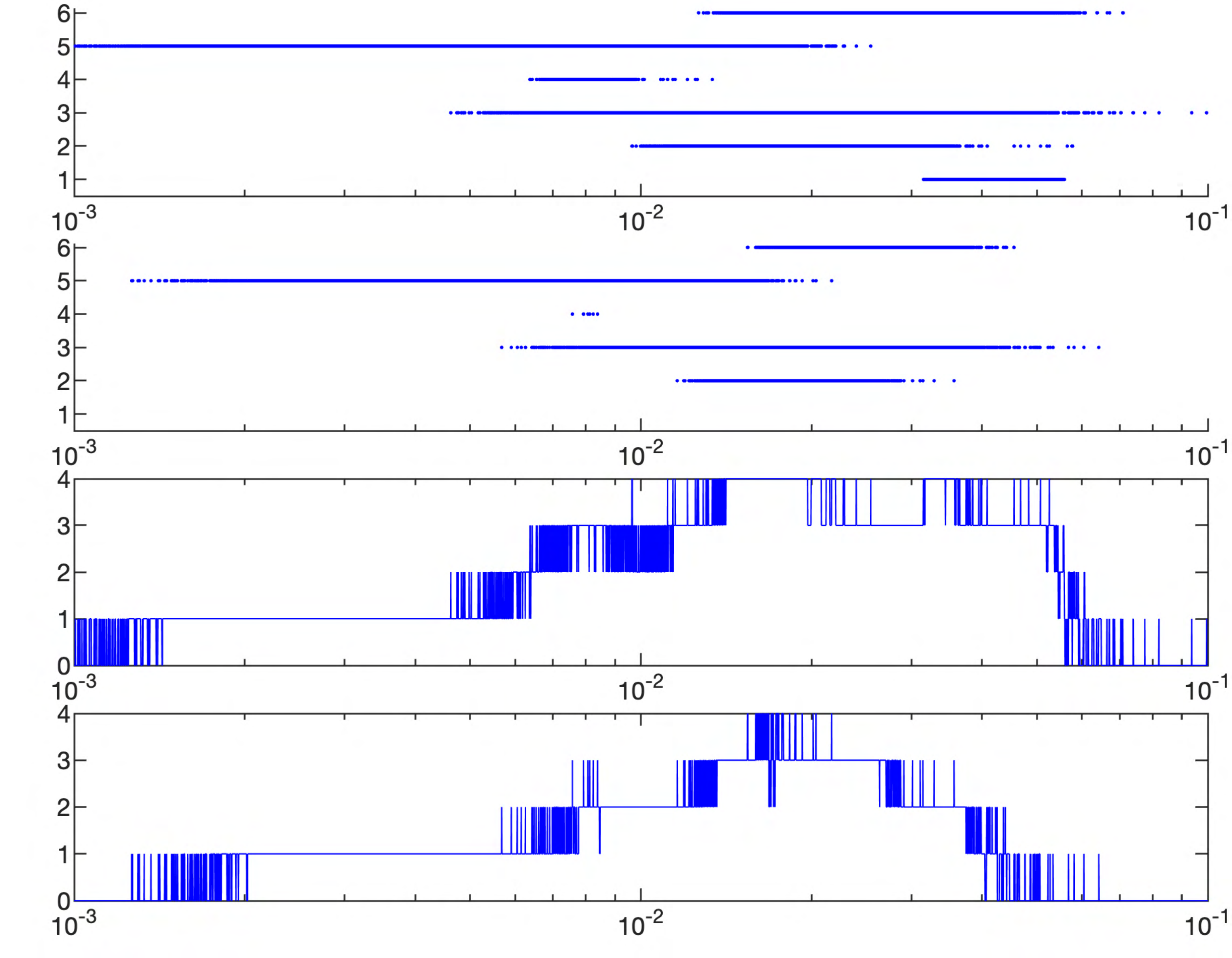}  \vspace*{-4mm}
	\caption{The top two plots display the values of $h$ for which PolySLC-R with $p=6$ and $C=0.1$ converged in less than 200 iterations for each of the six problems considered. The bottom two plots display the number of problems (out of 6) that PolySLC-R with $p=6$ and $C=0.1$ was able to solve in less than 200 iterations.  \label{fig: PolySLC_1D} }
\end{figure}

We solved the same problems using PolySLC-R with $p=6$ on a 2-dimensional grid of $200\times10000$ logarithmically-spaced points in $(C,h)$-space ($C$ between $10^{-15}$ and $10^{15}$, $h$ between $10^{-8}$ and $10^4$). The results, presented in Figure~\ref{fig: PolySLC_2D}, confirm the observation made in Section~\ref{section: Parameter C} that there is no universally optimal value of $C$. However, $C=0.1$ is an intermediate value which seems to work well for most problems, so we will set it as the default value, but it might need some tuning in practice. A similar experiment was conducted for $10^5$ logarithmically-spaced values of $h$ using PolySLC-R with $p=6$ and $C=0.1$, and Figure~\ref{fig: PolySLC_1D} shows that $h = 0.01$ could be a good default value.

\subsection{Tuning ExpoSLC-R}

\hfill \\

\begin{figure}[!h]
	\vspace*{-1mm}
	\centering
	\includegraphics[width=0.91\textwidth]{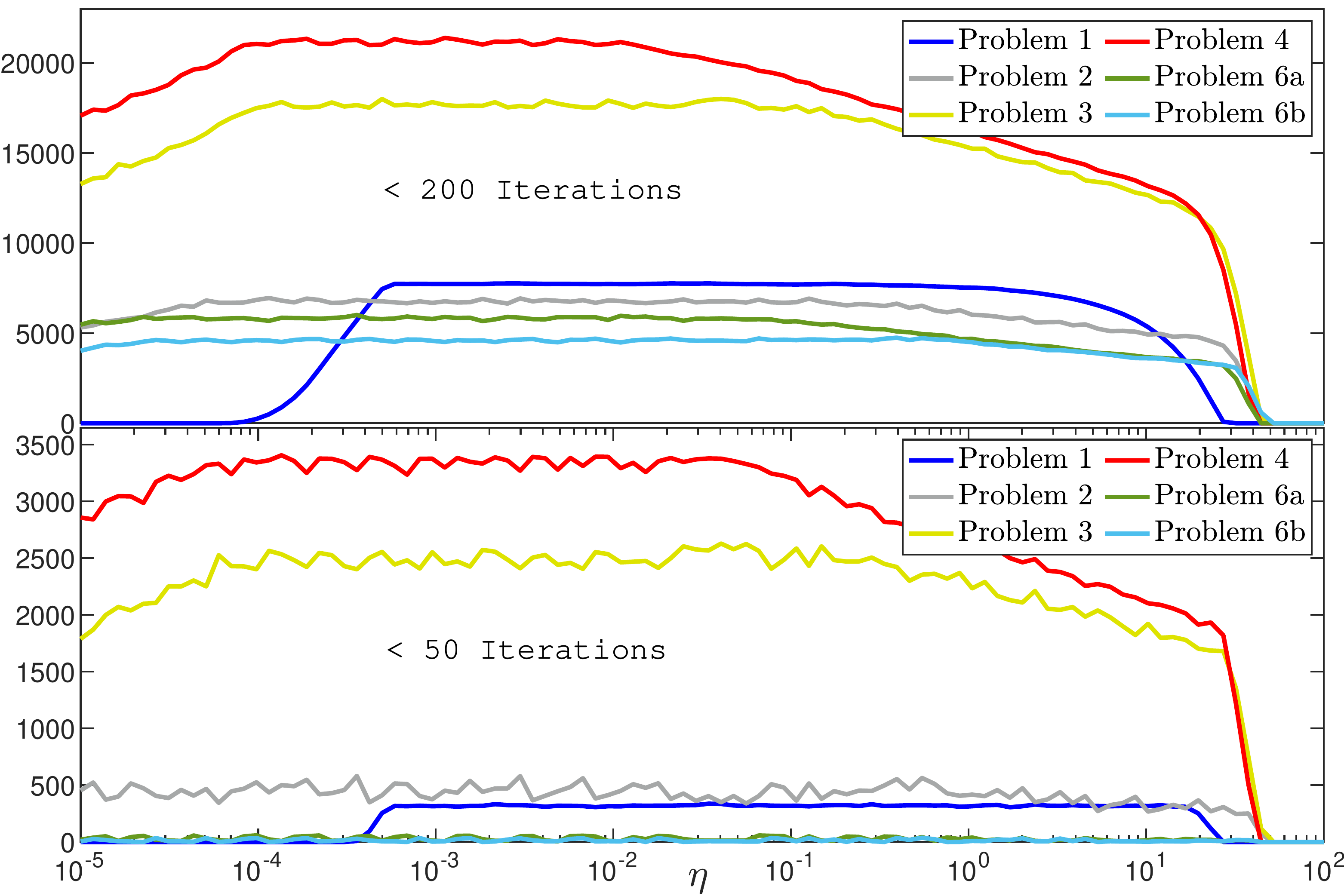}  \vspace*{-4mm}
	\caption{Number of $(C,h)$ pairs (out of $500^2$) for which convergence was achieved under 200 and 50 iterations using ExpoSLC-R, as the value of $\eta$ is varied.  \label{fig: ExpoSLC_eta} }
\end{figure}

We solved Problems~\ref{Problem: Quartic},~\ref{Problem: Polynomial and Log},~\ref{Problem: x log(x)},~\ref{Problem: Ill-Conditioned Quadratic}  and two distinct randomly generated instances of Problem~\ref{Problem: Linear Regression} using ExpoSLC-R on a 3-dimensional grid of logarithmically-spaced $500\times100\times500$ points in $(C,\eta,h)$-space ($C$ between $10^{-12}$ and $10^{12}$, $h$ between $10^{-6}$ and $10^3$, $\eta$ between $10^{-5}$ and $10^2$), and recorded the number of iterations needed to achieve convergence with~$\delta = 10^{-10}$. Figure~\ref{fig: ExpoSLC_eta} displays the number of $(C,h)$ pairs for which convergence was achieved under 200 and 50 iterations for each value of $\eta$. As for the polynomial Bregman algorithm, the value of $\eta$ does not seem to significantly affect the number of  $(C,h)$ pairs of fast convergence, as long as it falls between $0.001$ and $10$. Therefore, tuning the parameter $\eta$ carefully might not be very helpful and necessary, so we will fix it to $\eta= 0.01$.% \\

\begin{figure}[!h]
	\centering
	\includegraphics[width=0.91\textwidth]{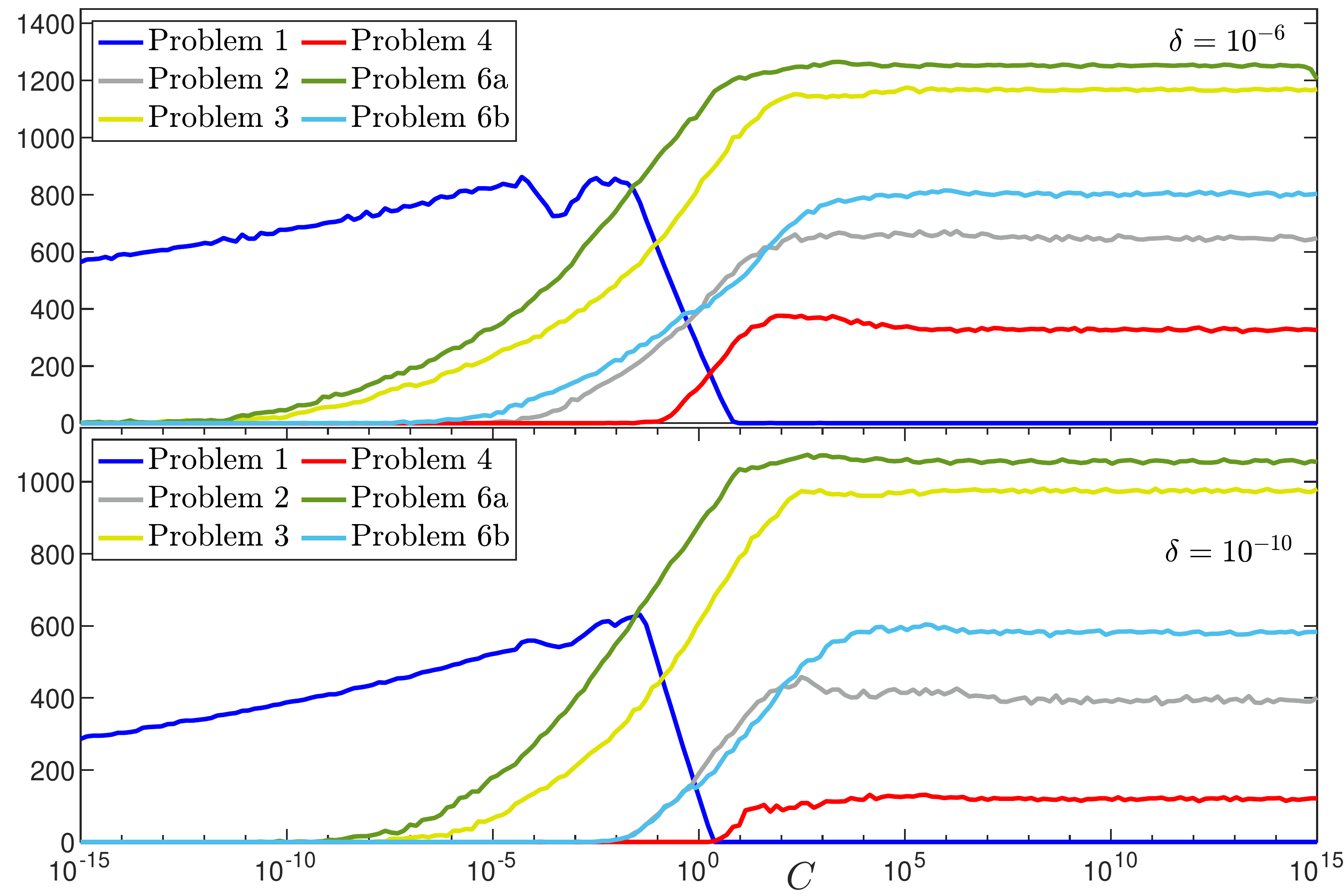}  \vspace*{-4mm}
	\caption{Number of values of $h$ (out of $10^4$) for which convergence was achieved under 200 iterations using ExpoSLC-R with $\eta =0.01$, as the value of $C$ is varied.  \label{fig: ExpoSLC_2D} }
\end{figure}

We then solved Problems~\ref{Problem: Quartic},~\ref{Problem: Polynomial and Log},~\ref{Problem: x log(x)},~\ref{Problem: Ill-Conditioned Quadratic}  and two distinct randomly generated instances of Problem~\ref{Problem: Linear Regression} using ExpoSLC-R with $\eta=0.01$ on a 2-dimensional grid of $200\times10000$ logarithmically-spaced points in $(C,h)$-space ($C$ between $10^{-15}$ and $10^{15}$, $h$ between $10^{-8}$ and $10^4$). The results, presented in Figure~\ref{fig: ExpoSLC_2D} confirm the observation made in Section~\ref{section: Parameter C} that there is no value of $C$ which is  universally optimal. However, $C=1$ seems to work well for most problems, so we will set it as the default value, but it might need some tuning in practice. \\

\begin{figure}[!h]
	\centering
	\includegraphics[width=1\textwidth]{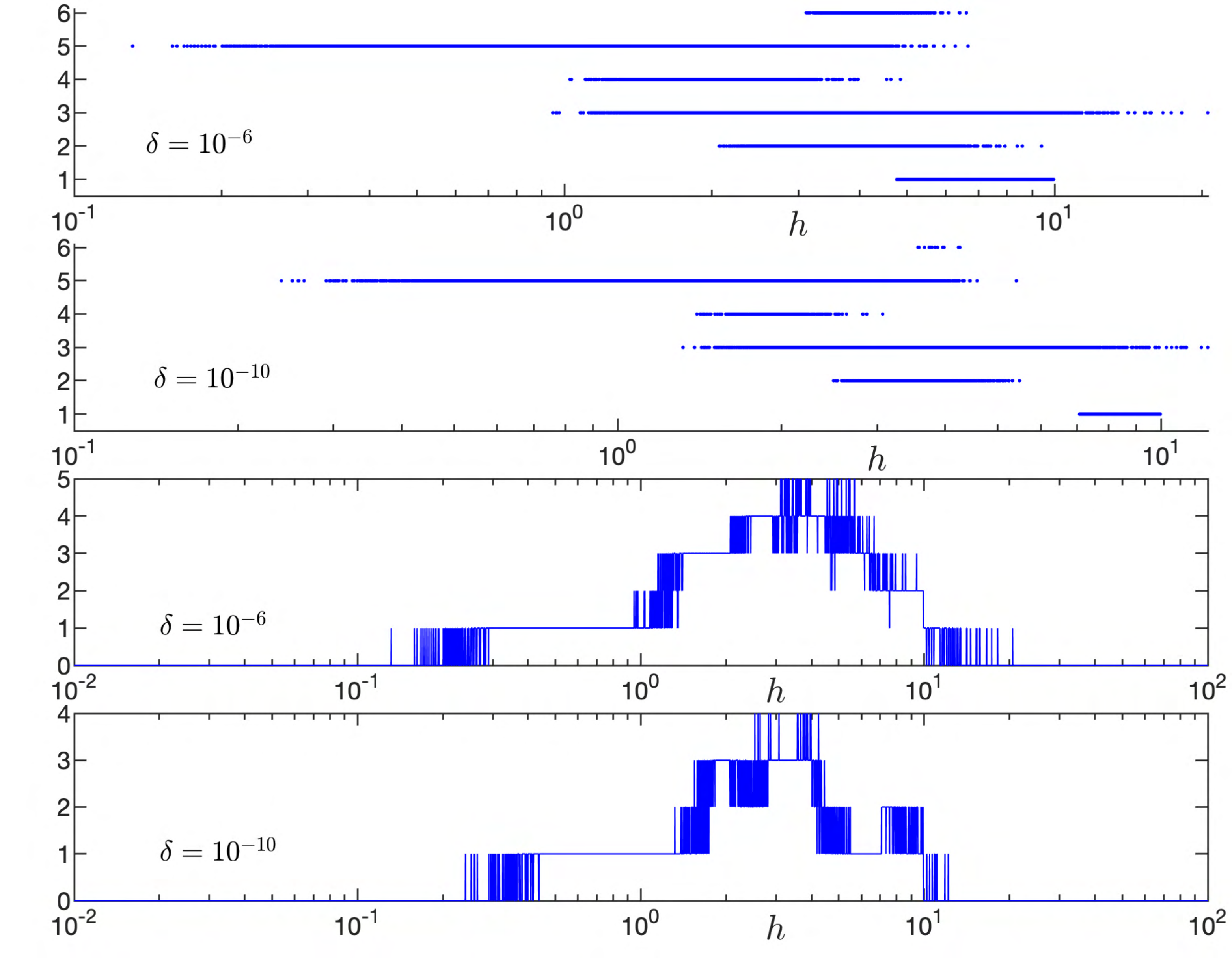} 
	\caption{The top two plots display the values of $h$ for which ExpoSLC-R with $\eta=0.01$ and $C=1$ converged in less than 200 iterations for each of the six problems considered. The bottom two plots display the number of problems (out of 6) that convergence was achieved in less than 200 iterations.  \label{fig: ExpoSLC_1D} }
\end{figure}

A similar experiment was conducted for $10^5$ logarithmically-spaced values of $h$ using PolySLC-R with $\eta=0.01$ and $C=1$, and Figure~\ref{fig: ExpoSLC_1D} shows that $h =4$ seems to perform well on most problems considered here.  

\hfill \\

\newpage 

\section{Temporal Looping to Improve Numerical Stability} \label{section: Temporal Looping to Improve Numerical Stability}

\hfill

There is an important caveat to the promising performance observed for the optimization algorithms constructed in this paper. The evolution of the variables $q,\mathfrak{q}$ and $r$ associated with the exponential Poincar\'e Hamiltonian,
\begin{equation}
	\bar{H}^{\eta}(\bar{q},\bar{r})  =    \frac{\eta}{2 e^{\eta \mathfrak{q}}}  \langle r,r\rangle + C\eta  e^{2\eta \mathfrak{q}} f(q) +  \mathfrak{r},
\end{equation}  
is guided by Hamilton's equations,
\begin{equation}
	\dot{q} =    \eta e^{-\eta \mathfrak{q}} r, \qquad   \dot{r} = -  C\eta  e^{2\eta \mathfrak{q}} \nabla f(q),  \qquad  \dot{\mathfrak{q}} = 1.
\end{equation}

\noindent From these equations of motion, we can see that the time variable $\mathfrak{q}$ grows linearly without bound, and as a result quantities like $e^{\eta \mathfrak{q}}$ grow exponentially without bound. More precisely, looking at the updates of the ExpoSLC algorithm,
\begin{align}
	r   \leftarrow  r -  C\eta h   e^{2\eta \mathfrak{q}} \nabla f(q) , \qquad  \mathfrak{q} \leftarrow \mathfrak{q}+h,  \qquad 
	q   \leftarrow q +  \Delta q = q +	\eta h e^{-\eta \left( \mathfrak{q} + \frac{h}{2} \right)}  r ,
\end{align}
we have at every iteration that the new update is given by
\begin{align} \label{eq: q-update for TL}
	(\Delta q )_{new} \leftarrow A (\Delta q)_{previous} + B e^{\eta \mathfrak{q}}\nabla f(q) , \end{align}
for some constants $A$ and $B$. 

\begin{figure}[!h]
	\centering
	\includegraphics[width=1\textwidth]{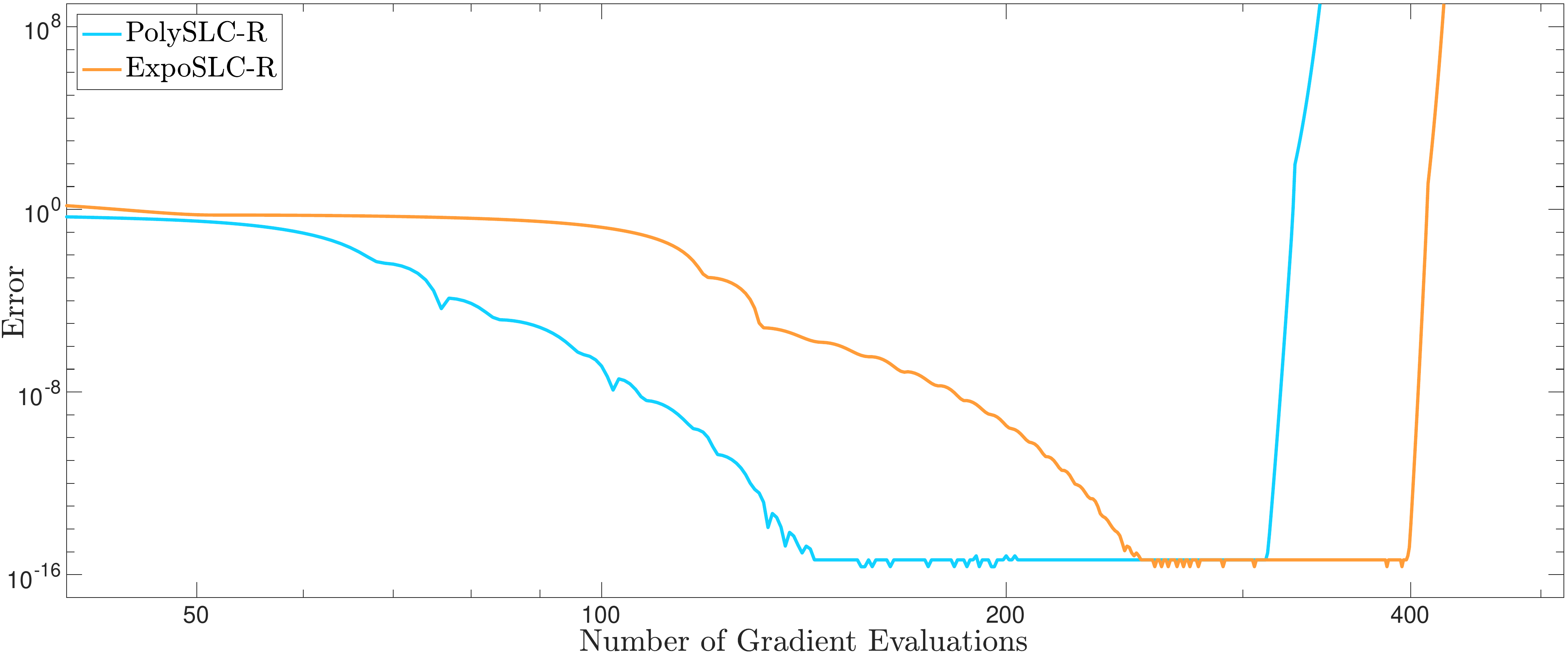} 
	\caption{The PolySLC-R and ExpoSLC-R algorithms applied to Problem~\ref{Problem: Polynomial and Log}. \label{fig: Instability} }
\end{figure}

If we could perform all the operations exactly, the gradient $\nabla f$ would converge to 0 with arbitrary precision and neutralize the unbounded growth of $e^{\eta \mathfrak{q}}$, and the quantity $B e^{\eta \mathfrak{q}} \nabla f(q)$ would remain very small. However, in practice, we can only perform operations with finite precision in floating-point arithmetic. As a result, $\nabla f$ only decays to 0 up to machine precision while $e^{\eta \mathfrak{q}}$ grows without bound. Eventually, $B e^{\eta \mathfrak{q}} \nabla f(q)$ becomes large again and the position variable $q$ moves away from the equilibrium it found near its optimal value. Something analogous happens in the polynomial family of Bregman dynamics, except that the unbounded growing exponential is replaced by an unbounded growing polynomial. This numerical instability phenomenon is illustrated in Figure~\ref{fig: Instability} which displays the evolution of the error $|f(x_k) - f(x^*)|$ when the SLC algorithms are applied to Problem~\ref{Problem: Polynomial and Log}. We see that both algorithms first achieve convergence to machine precision, stay at the minimizer for a few hundred iterations, and finally are expelled away from the minimizer due to numerical instability.

In all our numerical experiments so far, the algorithms stopped when they reached a desired convergence criterion, so we did not observe this numerical instability issue as it happens only after convergence is achieved. However, in practice, optimization algorithms are often terminated after a specified number of iterations instead of a specified convergence criterion. Thus, we need a strategy to avoid this numerical instability phenomenon. Since the numerical instability results from the limitation imposed by machine precision on accurately representing the decay to 0 of $\nabla f(q)$ while the term $e^{\eta \mathfrak{q}}$ grows without bound, it is natural to try to restrict the growth of the term $e^{\eta \mathfrak{q}}$, by restricting the growth of $\mathfrak{q}$ (and similarly in the polynomial case).  \\

One possibility is to reset time whenever a certain numerical instability criterion is met, via $\mathfrak{q} \leftarrow \beta \mathfrak{q}$ for some $\beta \in (0,1)$. A larger $\beta$ is preferable to keep enough momentum in case convergence to the minimizer was only suboptimal when the numerical instability criterion was met, or if the algorithm is used in an online fashion or with a stochastic or mini-batch approach. It is also preferable to avoid values of $\beta$ very close to 1, since the algorithm would then always remain close to numerical instability, and could possibly become unstable if the criterion is not chosen very carefully. In practice, taking $\beta$ between 0.6 and 0.95 works well, by ensuring that a reasonable amount of momentum is kept while avoiding the numerical instability region.  

Alternatively, one could reset time via $\mathfrak{q} \leftarrow  \mathfrak{q} - \nu h$ for some $\nu >1$. A smaller $\nu$ is preferable to retain momentum, while $\nu$ should not be too close to 1 to avoid numerical instability. \\

In practice, we will reset time via $\mathfrak{q} \leftarrow \max(\epsilon ,  \beta \mathfrak{q})$ or $\mathfrak{q} \leftarrow  \max( \epsilon , \mathfrak{q} - \nu h)$, where $\epsilon$ is a small positive number, to avoid very small or negative values of time~$\mathfrak{q}$. This phenomenon where the time variable $\mathfrak{q}$ is stuck in a loop by resetting $\mathfrak{q} \leftarrow \beta \mathfrak{q}$ or $\mathfrak{q} \leftarrow  \mathfrak{q} - \nu h$ whenever numerical instability is near will be referred to as \textbf{Temporal Looping}.  \\

Improving the ExpoSLC-R algorithm via temporal looping yields Algorithm~\ref{Alg: ExpoSLC-RTL}:  \vspace{3mm}

\begin{algorithm}[H] \label{Alg: ExpoSLC-RTL}
	\DontPrintSemicolon
	
	\KwInput{An objective function $f : \mathbb{R}^d \rightarrow \mathbb{R}$. An initial guess $q\in \mathbb{R}^d$. \\ \hspace{12.5mm} Parameters $C,h,p>0$, $\beta \in (0,1)$ or $\nu > 1$.  }
	
	\nonl	\hfill \vspace{-2.5mm}
	
	$\epsilon \leftarrow 0.001, \quad \mathfrak{q} \leftarrow 1, \quad  G \leftarrow \nabla f(q), \quad  r \leftarrow -  \frac{1}{2} C\eta h   e^{2\eta \mathfrak{q}} G $
	
	\While{convergence criteria are not met}
	{		
		$\Delta q \leftarrow   	\eta h e^{-\eta \left( \mathfrak{q} + \frac{h}{2} \right)}  r  $ \\
		
		$q  \leftarrow q + \Delta q$ \\ 
		
		$G \leftarrow \nabla f(q)$ \\
		
		\textbf{if}  $G^\top \Delta q >0$ \textbf{then} restart momentum: $r  \leftarrow 0$  \\
		
		\textbf{if}  \textit{numerical instability criterion is met} \textbf{then}  $\mathfrak{q} \leftarrow \max(\epsilon ,  \beta \mathfrak{q})$ or $\mathfrak{q} \leftarrow  \max( \epsilon , \mathfrak{q} - \nu h)$ \\
		
		$\mathfrak{q}  \leftarrow  \mathfrak{q} + h$ \\ 
		
		$r  \leftarrow  r -  C\eta h   e^{2\eta \mathfrak{q}} G $\\
	} 
	
	\caption{\textbf{ExpoSLC-RTL}: \textbf{S}ymmetric \textbf{L}eapfrog \textbf{C}omposition for the \textbf{Expo}nential Bregman dynamics, with \textbf{R}estarting and \textbf{T}emporal \textbf{L}ooping  }
\end{algorithm} 

\hfill  \\

In our numerical experiments with ExpoSLC-RTL, we use the instability criterion 
\begin{equation}
	C h^2  \eta^2 e^{\eta \mathfrak{q}} \| G \|  > e^{-\eta h} \| \Delta q \|,
\end{equation} 
to reset time. This criterion roughly ensures that \begin{equation}|B| e^{\eta \mathfrak{q}} \|\nabla f(q)\| < |A| \|  (\Delta q)_{previous} \| \end{equation} in equation~\eqref{eq: q-update for TL}, so that $(\Delta q)_{new}$ is not significantly larger in norm than $(\Delta q)_{previous}$.

\newpage

Improving the PolySLC-R algorithm via temporal looping yields Algorithm~\ref{Alg: PolySLC-RTL}:  \vspace{1.5mm}

\begin{algorithm}[H] \label{Alg: PolySLC-RTL}
	\DontPrintSemicolon
	
	\KwInput{An objective function $f : \mathbb{R}^d \rightarrow \mathbb{R}$. An initial guess $q\in \mathbb{R}^d$. \\ \hspace{12.5mm} Parameters $C,h,p>0$, $\beta \in (0,1)$ or $\nu > 1$. } 
	
	\nonl	\hfill \vspace{-2.5mm}
	
	$\epsilon \leftarrow 0.001, \quad \mathfrak{q} \leftarrow 1, \quad  G \leftarrow \nabla f(q), \quad  r \leftarrow  - \frac{1}{2} Chp \mathfrak{q}^{2p-1 } G$
	
	\While{convergence criteria are not met}
	{		
		$\Delta q \leftarrow   hp \left( \mathfrak{q} + \frac{h}{2} \right)^{-p-1}  r$ \\
		
		$q  \leftarrow q + \Delta q$ \\ 
		
		$G \leftarrow \nabla f(q)$ \\
		
		\textbf{if}  $G^\top \Delta q >0$ \textbf{then} restart momentum: $r  \leftarrow 0$  \\
		
		\textbf{if}  \textit{numerical instability criterion is met} \textbf{then}  $\mathfrak{q} \leftarrow \max(\epsilon ,  \beta \mathfrak{q})$ or $\mathfrak{q} \leftarrow  \max( \epsilon , \mathfrak{q} - \nu h)$ \\
		
		$\mathfrak{q}  \leftarrow  \mathfrak{q} + h$ \\ 
		
		$r  \leftarrow  r - Chp \mathfrak{q}^{2p-1 } G$\\
	}

	\caption{\textbf{PolySLC-RTL}: \textbf{S}ymmetric \textbf{L}eapfrog \textbf{C}omposition for the \textbf{Poly}nomial Bregman dynamics, with \textbf{R}estarting and \textbf{T}emporal \textbf{L}ooping }
\end{algorithm} 
\hfill 

In our numerical experiments, we have chosen the numerical instability criterion
\begin{equation}
	C  h^2 p^2  (\mathfrak{q}+h)^{p+1} \| G \|  >  \mathfrak{q} \| \Delta q \|, 
\end{equation}
which roughly ensures that the new position update is not significantly larger than the previous one. % \\

Figure~\ref{fig: Instability Fix} shows that temporal looping takes care of the numerical instability issue experienced earlier in Figure~\ref{fig: Instability} for Problem~\ref{Problem: Polynomial and Log}:
\begin{figure}[!h]
	\centering
	\begin{minipage}[b]{0.98\textwidth}
		\includegraphics[width=\textwidth]{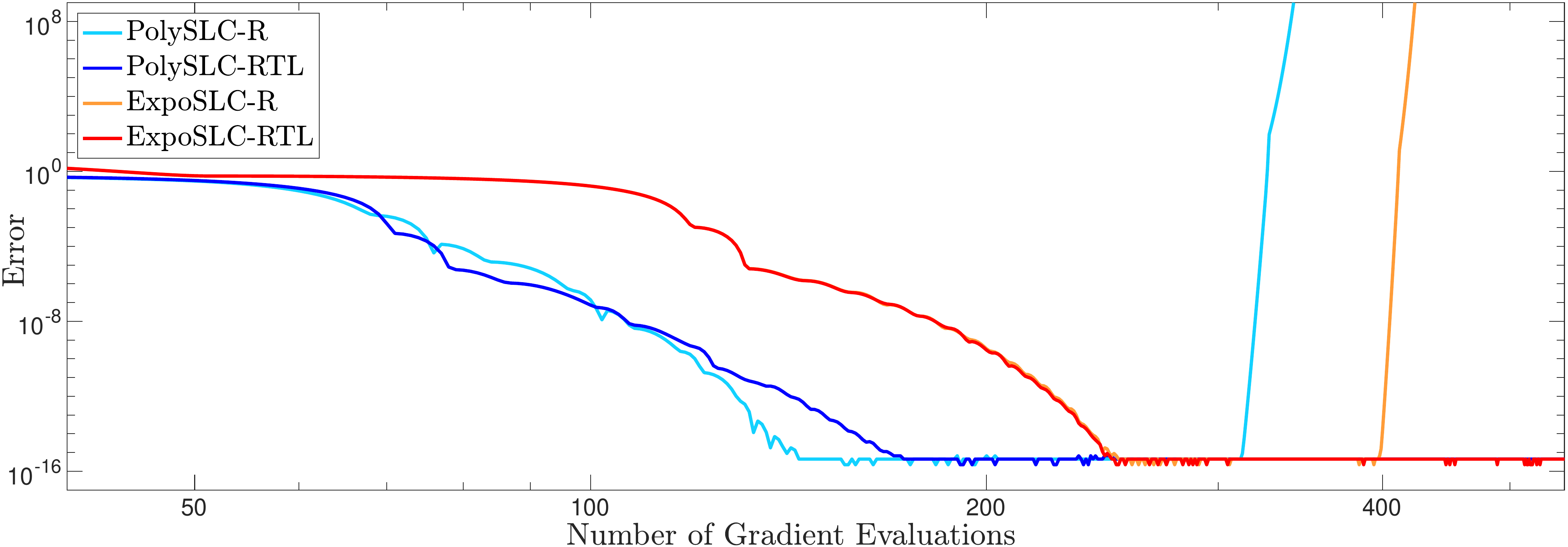} 
	\end{minipage}
	\begin{minipage}[b]{0.98\textwidth}
		\includegraphics[width=\textwidth]{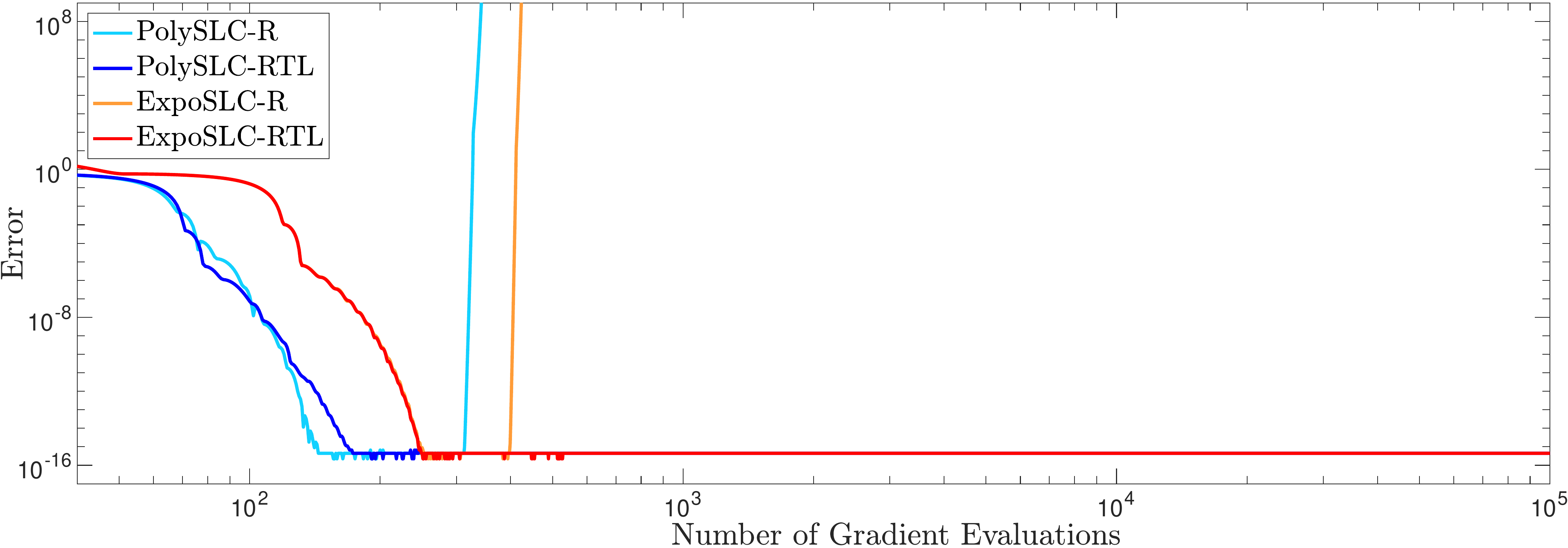} 	
	\end{minipage}  \vspace*{-2.5mm}
	\caption{The effect of temporal looping in PolySLC-R and ExpoSLC-R. \label{fig: Instability Fix} }
\end{figure}

%\newpage 

\begin{figure}[!h] 
	\centering
	\includegraphics[width=1\textwidth]{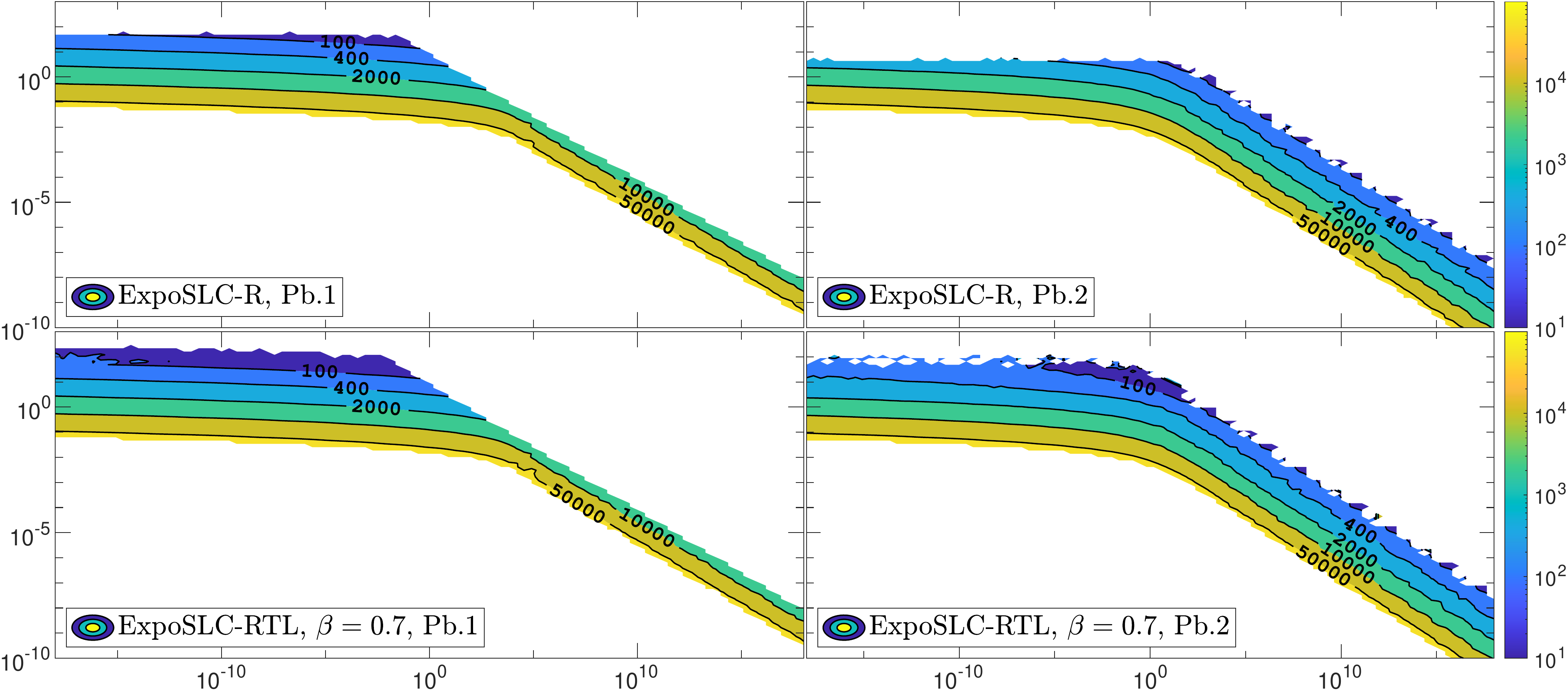} 	\vspace*{-7.5mm} \caption{Contour plot of the number of iterations required to achieve convergence ($\delta = 10^{-8}$) in the $(C,h)$-plane, for the ExpoSLC-R and ExpoSLC-RTL algorithms, when applied with $\eta = 0.01$ to Problems~\ref{Problem: Quartic} and \ref{Problem: Polynomial and Log}.} \label{fig: ExpoTLRobustness} 
	\vspace*{-5mm}
\end{figure}

\begin{figure}[!h] 
	\centering
	\includegraphics[width=1\textwidth]{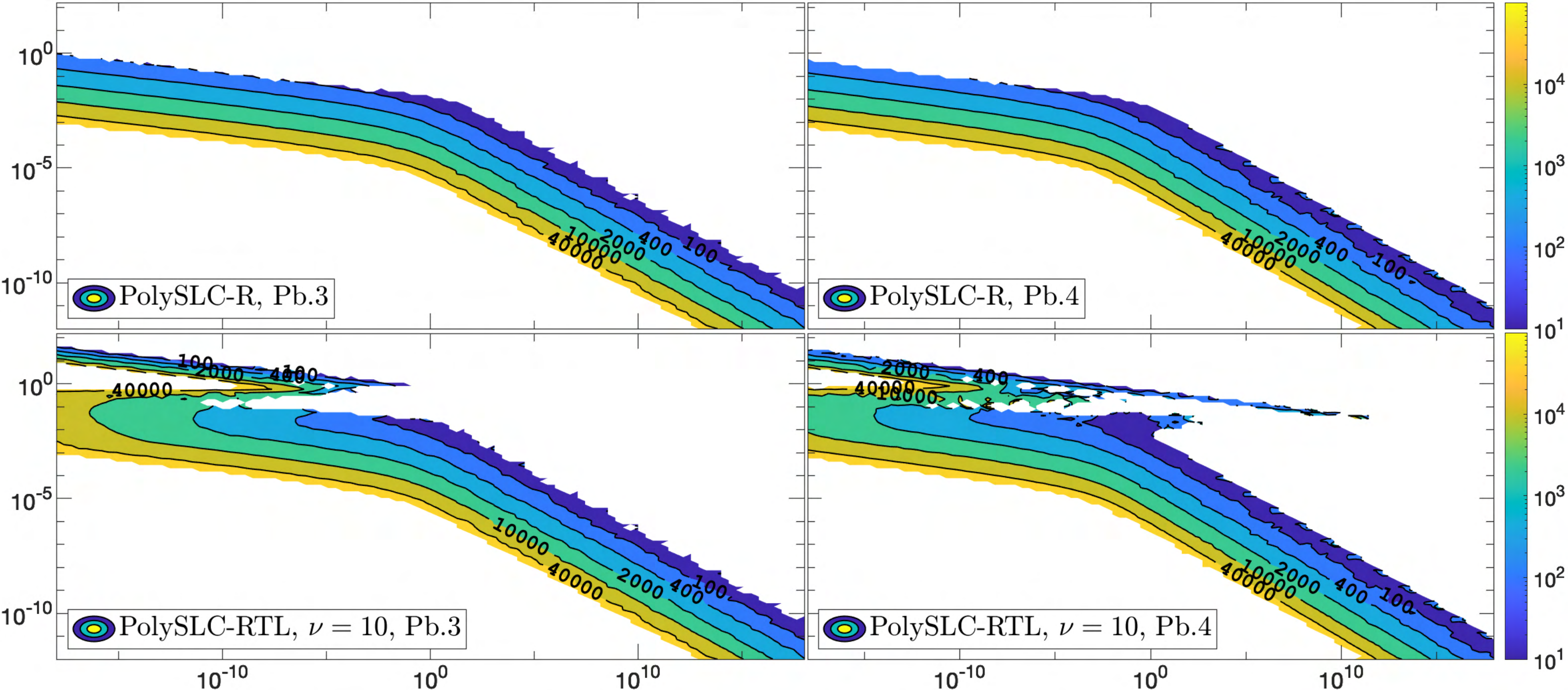} 	\vspace*{-7.5mm} \caption{Contour plot of the number of iterations required to achieve convergence ($\delta = 10^{-8}$) in the $(C,h)$-plane, for the PolySLC-R and PolySLC-RTL algorithms, when applied with $p=10$ to Problems~\ref{Problem: x log(x)} and \ref{Problem: Ill-Conditioned Quadratic}.}\label{fig: PolyTLRobustness} 
\end{figure}

It can be seen from Figures~\ref{fig: ExpoTLRobustness} and \ref{fig: PolyTLRobustness} that temporal looping, with the $\mathfrak{q} \leftarrow \max(\epsilon ,  \beta \mathfrak{q})$ scheme or $\mathfrak{q} \leftarrow  \max( \epsilon , \mathfrak{q} - \nu h)$ scheme, does not negatively affect the performance of the algorithms, although the algorithms with temporal looping might sometimes require a larger number of iterations to achieve convergence for a fixed $(C,h)$-pair. Indeed the regions of fast convergence might be shifted slightly, but remained at least as large if not larger when using temporal looping.


Overall, we have seen that temporal looping can be very helpful to deal with post-convergence numerical instability, and that it does not affect negatively the initial performance of the algorithm. Note that temporal looping could be improved by tuning the parameters $\beta$ or $\nu$, or by designing a better suited numerical instability criterion. 
%\newpage 

%\newpage 

\section{Testing for Machine Learning Applications} \label{section: ML Testing}

%\hfill 

We now test our algorithms on more challenging optimization problems for machine learning with a variety of model architectures, loss functions, and applications. For most of the problems considered in this section, the gradients are evaluated in a mini-batch fashion. For reference, we will also solve these optimization problems using gradient descent and the most commonly used optimizer in machine learning, Adam~\cite{ADAM}:  \vspace{2mm}

\begin{center}
	\textbf{\underline{ADAM}} 
\end{center}
\begin{align*}
	m_{k+1} & = \beta_1 m_{k} + (1-\beta_1) \nabla f(x_{k})  \\ 
	v_{k+1}  & = \beta_2 v_{k} + (1-\beta_2) \nabla f(x_{k}) \odot \nabla f(x_{k})
	\\ 
	x_{k+1}  & = x_{k} - h \left[  (1-\beta_1^{k+1}) \left( \sqrt{ (1-\beta_2^{k+1})^{-1} v_{k+1} } + \epsilon \right) \right]^{-1}m_{k+1} 
\end{align*}

\vspace{2mm}

\noindent Here,  $u \odot v$ denotes elementwise multiplication, and the update for $x_{k+1}$ is performed elementwise. The variable $\epsilon$ present in the updates of the Adam algorithm is there to avoid numerical instability associated with division by 0 (with default value $\epsilon = 10^{-8}$ in \cite{ADAM} and PyTorch). The three parameters of Adam are $\beta_1, \beta_2$ used for computing running averages of gradients, and the learning rate $h$ (with default values $\beta_1 = 0.9$, $\beta_2 = 0.999$, $h = 0.001$ in \cite{ADAM}, PyTorch and TensorFlow).  \\

The ExpoSLC-RTL and PolySLC-RTL algorithms have been implemented under the more evocative names \textbf{eBrAVO} and \textbf{pBrAVO} (\textbf{Br}egman \textbf{A}ccelerated \textbf{V}ariational \textbf{O}ptimizer) in Python so that they can be used within the TensorFlow and PyTorch frameworks. These algorithms are available at \href{https://github.com/vduruiss/AccOpt_via_GNI}{github.com/vduruiss/AccOpt\_via\_GNI}, and can be called in a similar way as the Adam algorithm in TensorFlow and PyTorch: 
\begin{lstlisting}[language=Python]
	optimizer = tf.keras.optimizers.Adam(learning_rate = 0.001)
	optimizer = BrAVO_tf.eBravo(learning_rate = 1)
\end{lstlisting}
\vspace{-0.5mm}
\begin{lstlisting}[language=Python]
	optimizer = torch.optim.Adam(model.parameters(), lr = 0.01)
	optimizer = BrAVO_torch.eBravo(model.parameters(), lr = 1)
\end{lstlisting}

\hfill 

The purpose of this section is not to do a very careful computational comparison of the BrAVO algorithms with commonly used optimization algorithms in machine learning but rather to show that the BrAVO algorithms can be used conveniently within the PyTorch and TensorFlow frameworks for numerous concrete machine learning applications, and that they might be worth considering and improving in the future. A very careful computational comparison of optimization algorithms for machine learning is a much more ambitious goal which is beyond the scope of this paper, and would be more meaningful once the implementation of the BrAVO algorithms within the PyTorch and TensorFlow frameworks has been highly-optimized.  \\

\begin{figure}[!h]
	\centering
	\includegraphics[width=1\textwidth]{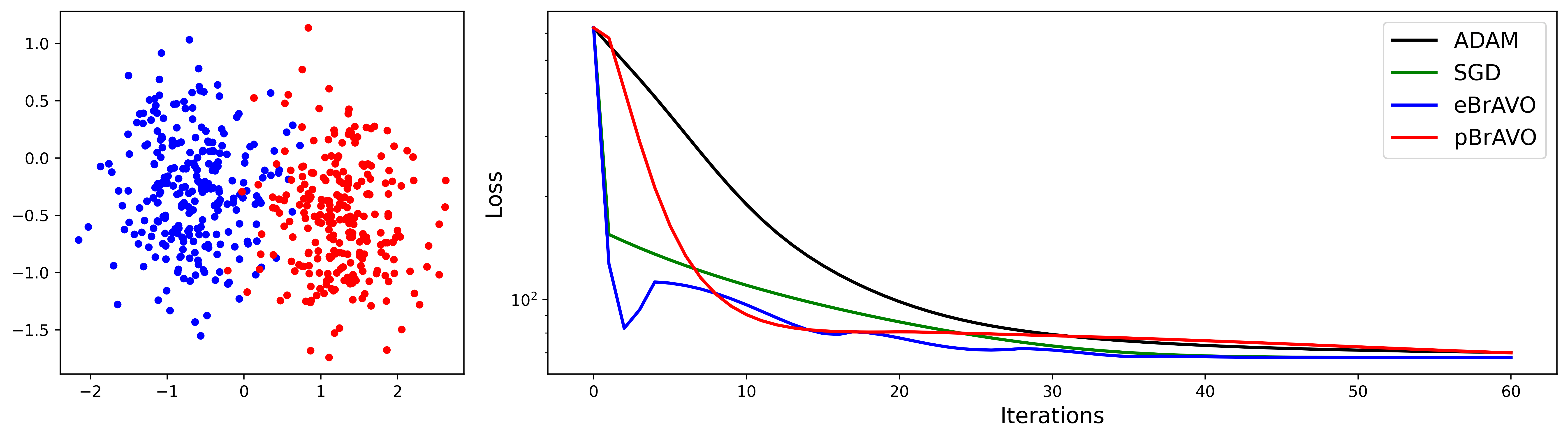} 	\vspace*{-8.5mm}
	\caption{Comparison of algorithms applied to a Binary Classification Problem~\ref{Problem: Binary Classification}.  \label{fig: Binary Classification}  }
\end{figure}

\begin{figure}[!h]
	\centering
	\includegraphics[width=1\textwidth]{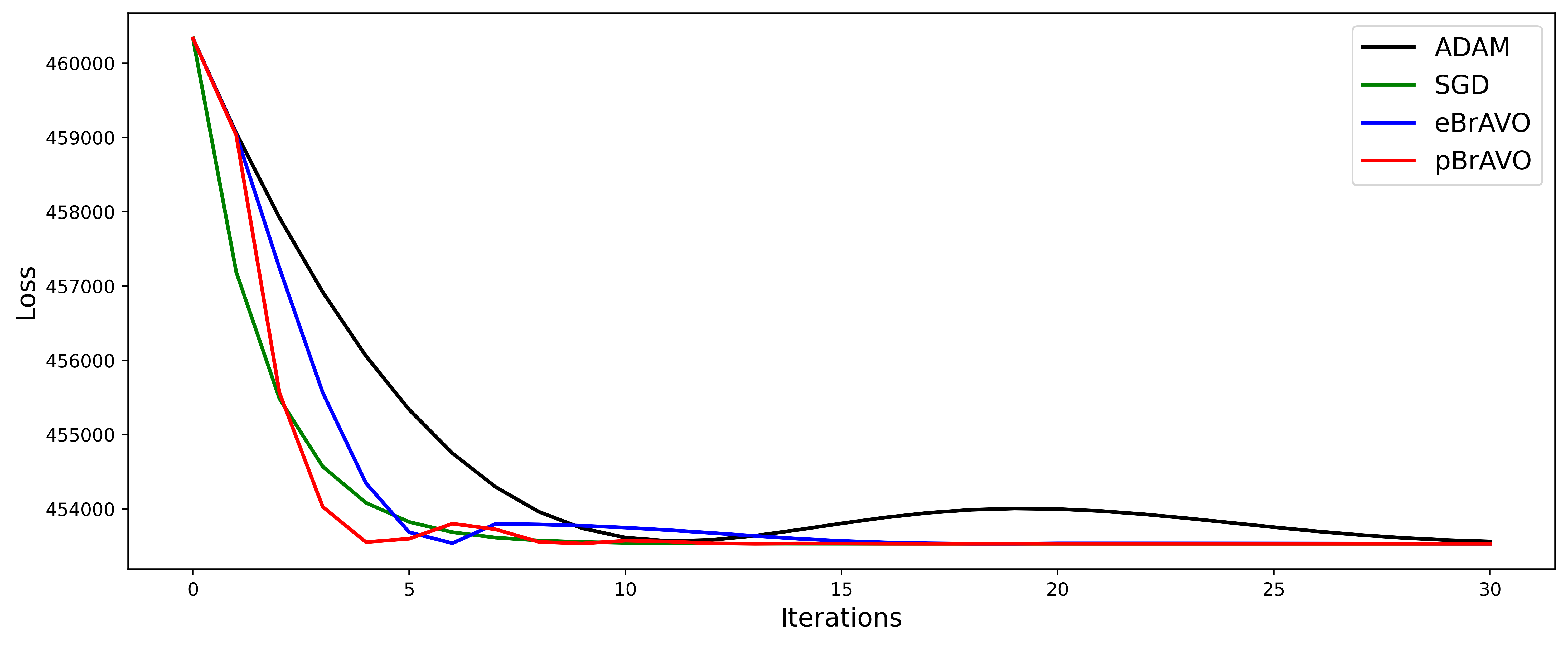} 	\vspace*{-9mm}
	\caption{Comparison of algorithms applied to a Location Problem~\ref{Problem: Location}.  \label{fig: Location} \vspace{1.8mm} }
\end{figure}

We have first tested the performance of our algorithms with automatic differentiation on instances of the Binary Classification Problem~\ref{Problem: Binary Classification} and the Fermat--Weber Location Problem~\ref{Problem: Location}. Figure~\ref{fig: Binary Classification} shows the evolution of the loss function~\eqref{eq: Binary Classification Loss} when formulating a model separating blue and red regions of 2-dimensional space using a line based on the displayed 500 randomly generated points. Figure~\ref{fig: Location} shows the evolution of the loss function~\eqref{eq: Location Loss} when solving the Fermat--Weber Location Problem~\ref{Problem: Location} with 5000 randomly generated vectors in $\mathbb{R}^{1000}$ and 5000 randomly generated corresponding scalar weights. We can see from Figures~\ref{fig: Binary Classification} and \ref{fig: Location} that our algorithms solve the binary classification and location problems with an accuracy and efficiency comparable to those of the Adam and standard gradient descent~(SGD) algorithms implemented in TensorFlow. \\

\newpage

	\vspace*{-8mm}

Next, we tested our algorithm on the popular multi-label image classification problem based on the Fashion--MNIST dataset \cite{FashionMNIST}: \textit{`Fashion--MNIST is a dataset of Zalando's article images consisting of a training set of 60,000 examples and a test set of 10,000 examples. Each example is a 28$\times$28 grayscale image, associated with a label from 10 classes (t-shirt/top, trouser, pullover, dress, coat, sandal, shirt, sneaker, bag, ankle boot)'}. We use {\ttfamily nn.CrossEntropyLoss()} as the loss function, and the following neural network architecture in PyTorch as our classification model: \\

\vspace{-1.4mm}
\lstset{basicstyle=\ttfamily \scriptsize }
\begin{mdframed}[backgroundcolor=backcolour]\vspace{-2.5mm}\begin{lstlisting}[language=Python] 		
		Layer (type)		Output Shape			Parameters
		========================================================
		dense (Dense)		[-1, 784] 			0
		dense_1 (Dense)	[-1, 64]				50,240
		dense_2 (Dense)	[-1, 64]				0
		dense_3 (Dense)	[-1, 64]				4,160
		========================================================
		Total Number of Parameters: 55,050 \end{lstlisting} \vspace{-2.5mm} \end{mdframed}

\vspace{3mm}

\begin{figure}[!h]
	\centering
	\includegraphics[width=1\textwidth]{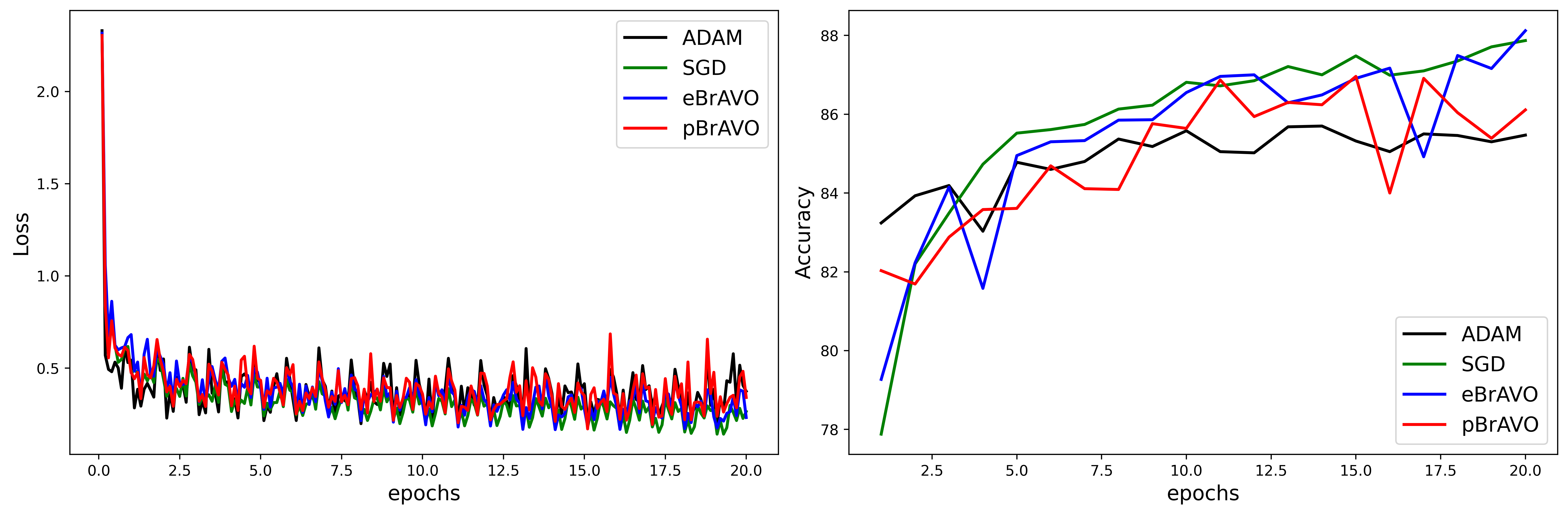}\vspace*{-2mm}
	\caption{Evolution of the loss function and accuracy (in \%) of Adam, standard gradient descent (SGD) and the BrAVO algorithms,   when applied to the Fashion--MNIST multi-label classification problem.  \label{fig: MNIST} }
\end{figure}

\begin{figure}[!h]
	\centering
	\includegraphics[width=1\textwidth]{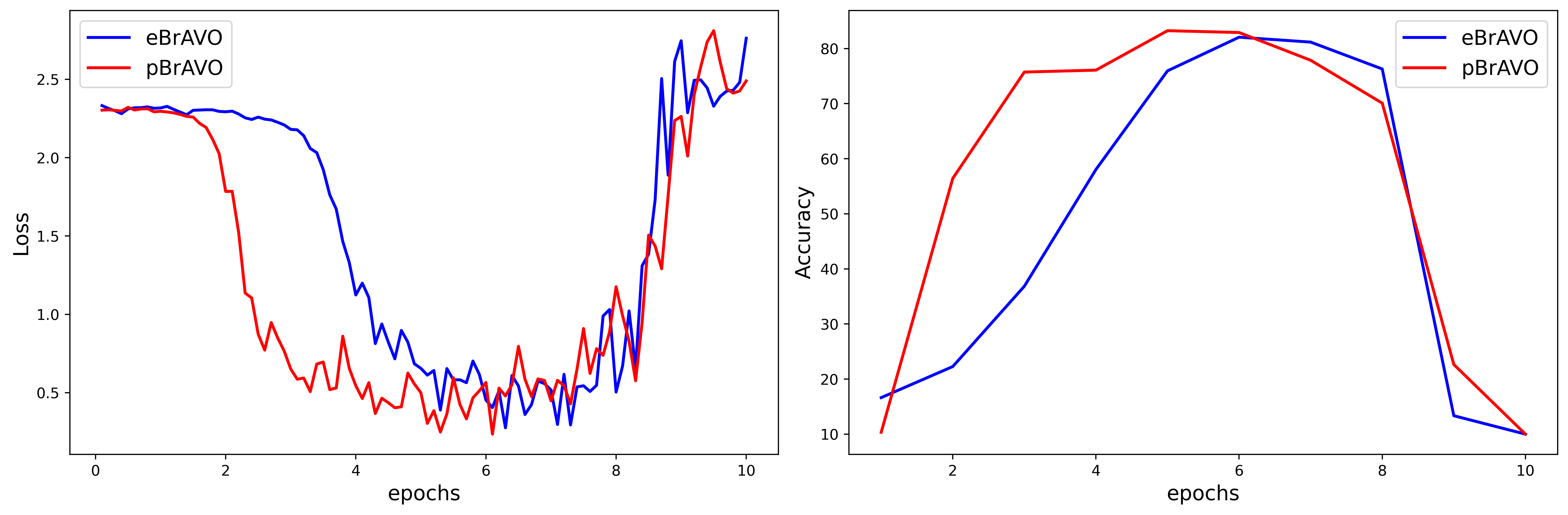}\vspace*{-2mm} 
	\caption{Convergence and loss of convergence for the BrAVO algorithms without momentum restarting and temporal looping, when applied to the Fashion--MNIST multi-label classification problem.  \label{fig: unstable MNIST}  \vspace{3mm}}
\end{figure}

Figure~\ref{fig: MNIST} shows that the BrAVO algorithms achieve comparable accuracy and efficiency on the Fashion--MNIST classification problem as the Adam and gradient descent (SGD) algorithms. Note that the momentum restarting scheme and the temporal looping strategy are essential to the good behavior of the algorithms. Indeed, we can see from Figure~\ref{fig: unstable MNIST} that without them, the algorithms eventually lose convergence due to numerical instability. Note as well that these strategies can also allow for larger time-steps which usually translates into faster convergence. \\

%\newpage 

%\hfill 


\begin{figure}[!h]
	\centering
	\includegraphics[width=1\textwidth]{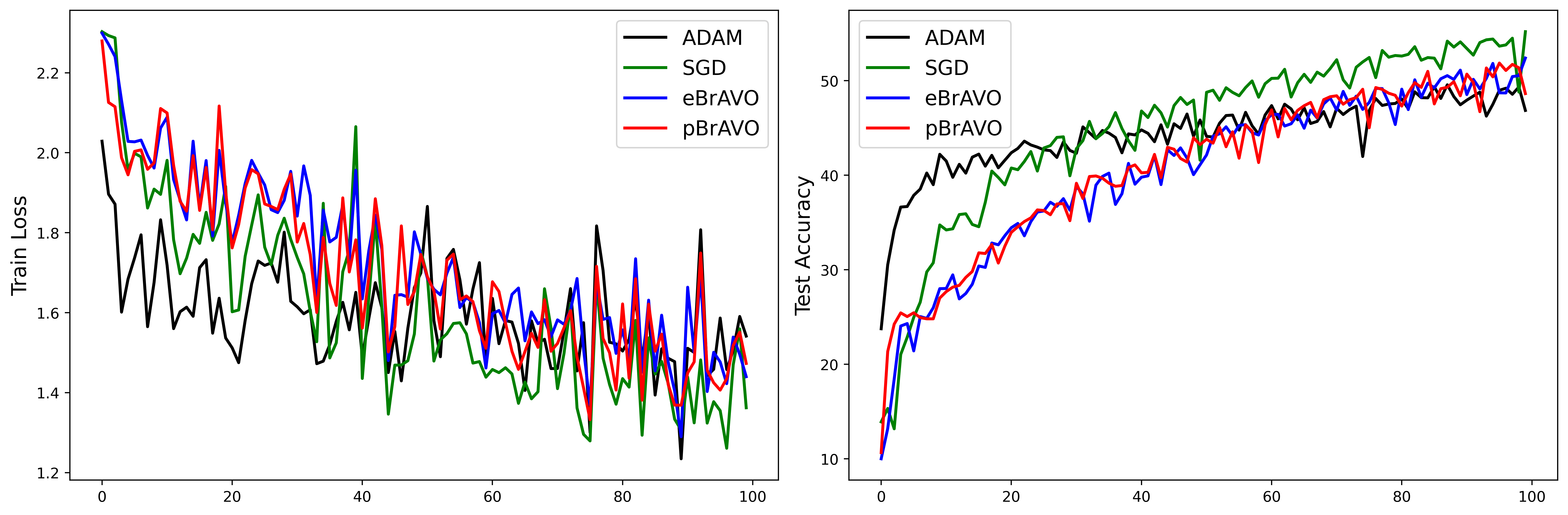}\vspace*{-3.2mm} 
	\caption{Evolution over 20 epochs of the loss function and accuracy of various algorithms when applied to the CIFAR--10 multi-label image classification problem.  \label{fig: CIFAR} } 
\end{figure}

We then tested our algorithm on another popular multi-label image classification problem based on the CIFAR-10 dataset \cite{CIFAR}: \textit{`the CIFAR-10 dataset consists of 60000 32$\times$32 color images in 10 mutually exclusive classes (airplane, automobile, bird, cat, deer, dog, frog, horse, ship, truck), with 6000 images per class'}. The results are displayed in Figure~\ref{fig: CIFAR}. 

We used {\ttfamily nn.CrossEntropyLoss()} as the loss function to minimize and a small Convolutional Neural Network in PyTorch very similar to the LeNet-5 architecture first described in \cite{LeCun1998}:  \vspace{3mm}

\lstset{basicstyle=\ttfamily \scriptsize }
\begin{mdframed}[backgroundcolor=backcolour]\vspace{-1.5mm}\begin{lstlisting}[language=Python] 		
		Layer (type)		Output Shape			Parameters
		====================================================
		Conv2d-1			[-1, 6, 28, 28] 		456
		Conv2d-2			[-1, 16, 10, 10]		2,416
		Linear-3			[-1, 120]						48,120
		Linear-4			[-1, 84]						10,164
		Linear-5			[-1, 10]						850
		====================================================
		Total Number of Parameters: 62,006 \end{lstlisting} \vspace{-2.5mm} \end{mdframed}

\newpage 

Let us now consider the Natural Language Processing problem of constructing a multi-label text classifier which can provide suggestions for the most appropriate subject areas for arXiv papers based on their abstracts. The code and architecture used are based on the Keras tutorial~\cite{KerasNLP}. An arXiv paper can belong to multiple categories, so the prediction task can be divided into a series of multiple binary classification problems, and we can use the Binary Cross Entropy loss. We will use the following neural network architecture:  \vspace{2mm}

\lstset{basicstyle=\ttfamily\footnotesize }
\begin{mdframed}[backgroundcolor=backcolour]\vspace{-2.5mm}\begin{lstlisting}[language=Python]
model = keras.Sequential()
model.add(layers.Dense(units=256, activation='relu'))
model.add(layers.Dense(units=256, activation='relu'))
model.add(layers.Dense(units=lookup.vocabulary_size(), activation='sigmoid')) \end{lstlisting} \vspace{-2.5mm} \end{mdframed} \vspace{3mm}

\begin{figure}[!h]
\centering
\includegraphics[width=1\textwidth]{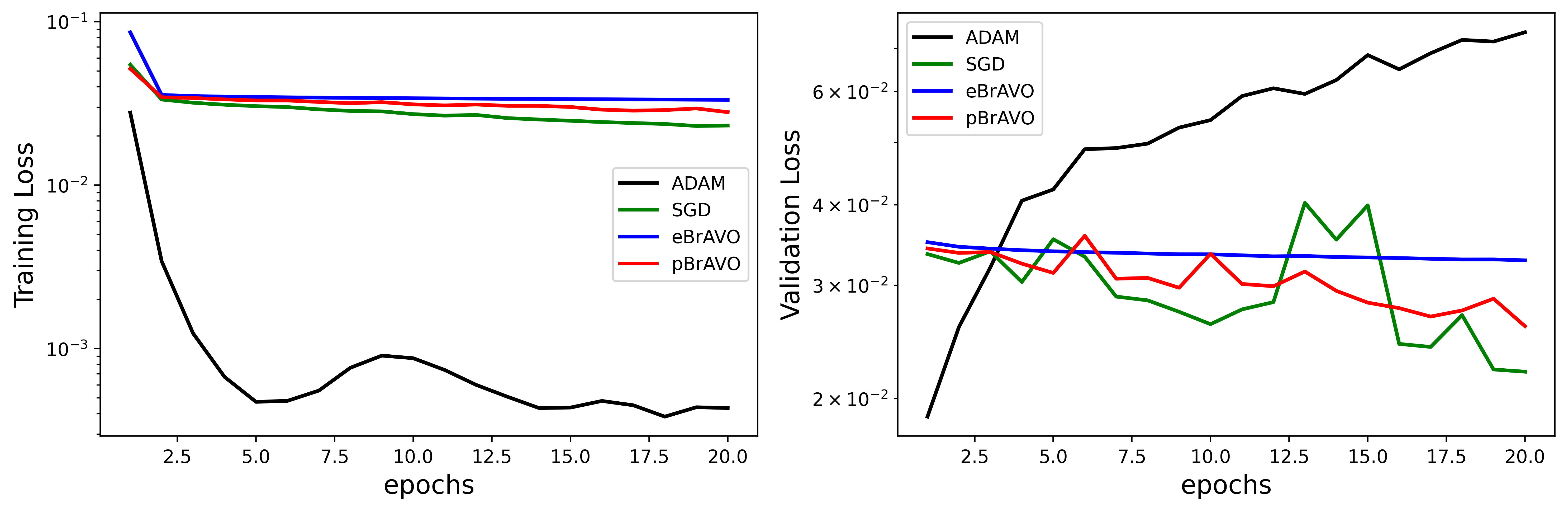} \vspace*{-7mm}
\caption{Evolution of the Binary Cross Entropy loss function on training and validation datasets for several algorithms, when applied to the Natural Language Processing problem of multi-label text classification of arXiv papers. \label{fig: NLP} }
\end{figure}

The evolution of the loss on the training and validation datasets is displayed in Figure~\ref{fig: NLP}. Although the Adam optimizer achieves the smallest loss on the training dataset, the resulting optimized model does not outperform the models generated by the other optimizers on the validation dataset. Its validation loss actually worsens as the epoch number increases (unlike for the other algorithms) which indicates that the optimized model might be suffering from overfitting.

Next, we consider timeseries forecasting for weather prediction, based on the Keras tutorial~\cite{KerasWeather}. We use the Mean Squared Error (MSE) as the loss function and the following Long Short-Term Memory (LSTM) model (with 5,153 parameters):  \vspace{2mm}

\lstset{basicstyle=\ttfamily\scriptsize }
\begin{mdframed}[backgroundcolor=backcolour]\vspace{-2.5mm}\begin{lstlisting}[language=Python]
inputs = layers.Input(shape=(inputs.shape[1], inputs.shape[2]))
lstm_out = layers.LSTM(32)(inputs)
outputs = layers.Dense(1)(lstm_out)
model = keras.Model(inputs=inputs, outputs=outputs)\end{lstlisting} \vspace{-2.5mm} \end{mdframed} \vspace{4mm}

\begin{figure}[!h]
\centering
\includegraphics[width=0.95\textwidth]{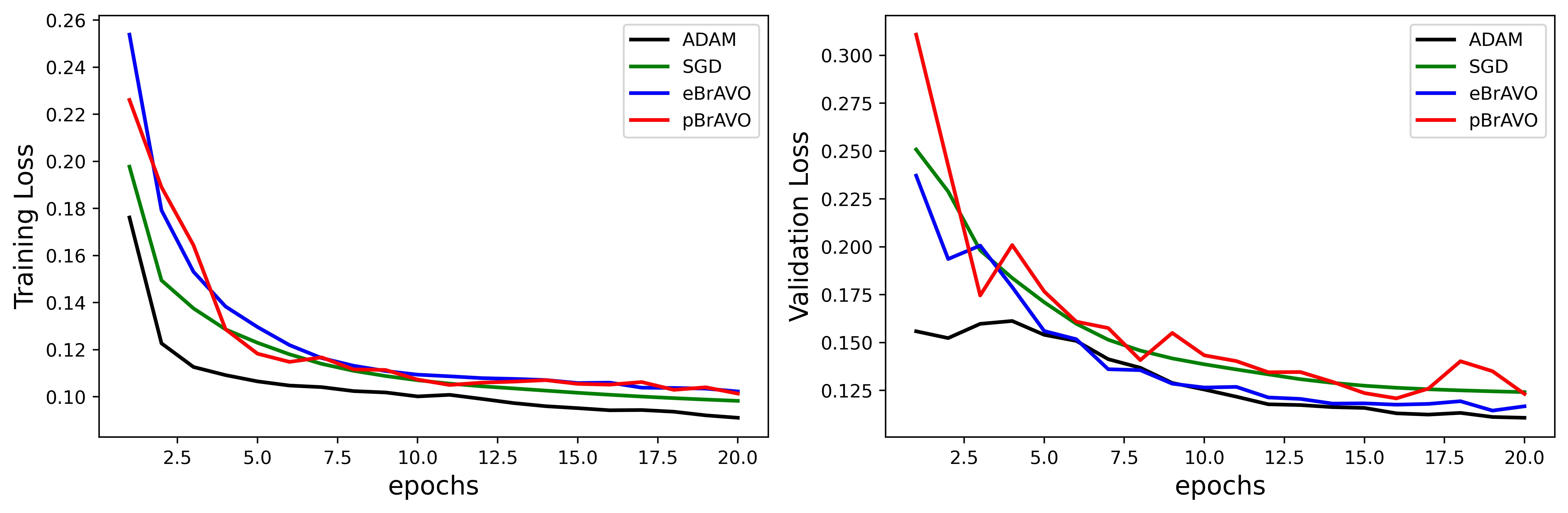}  \vspace*{-4mm}
\caption{MSE evolution on training and validation datasets for several algorithms, when used to optimize a LSTM model for timeseries forecasting for weather prediction. \label{fig: Weather} }
\vspace{-1mm}
\end{figure}

The evolution of the mean squared error on the training and validation sets is displayed in Figure~\ref{fig: Weather}. We can see that the four different algorithms generate similar losses on the training and validation datasets.

\begin{figure}[!h]
\vspace*{-1.8mm}
\centering
\includegraphics[width=0.85\textwidth]{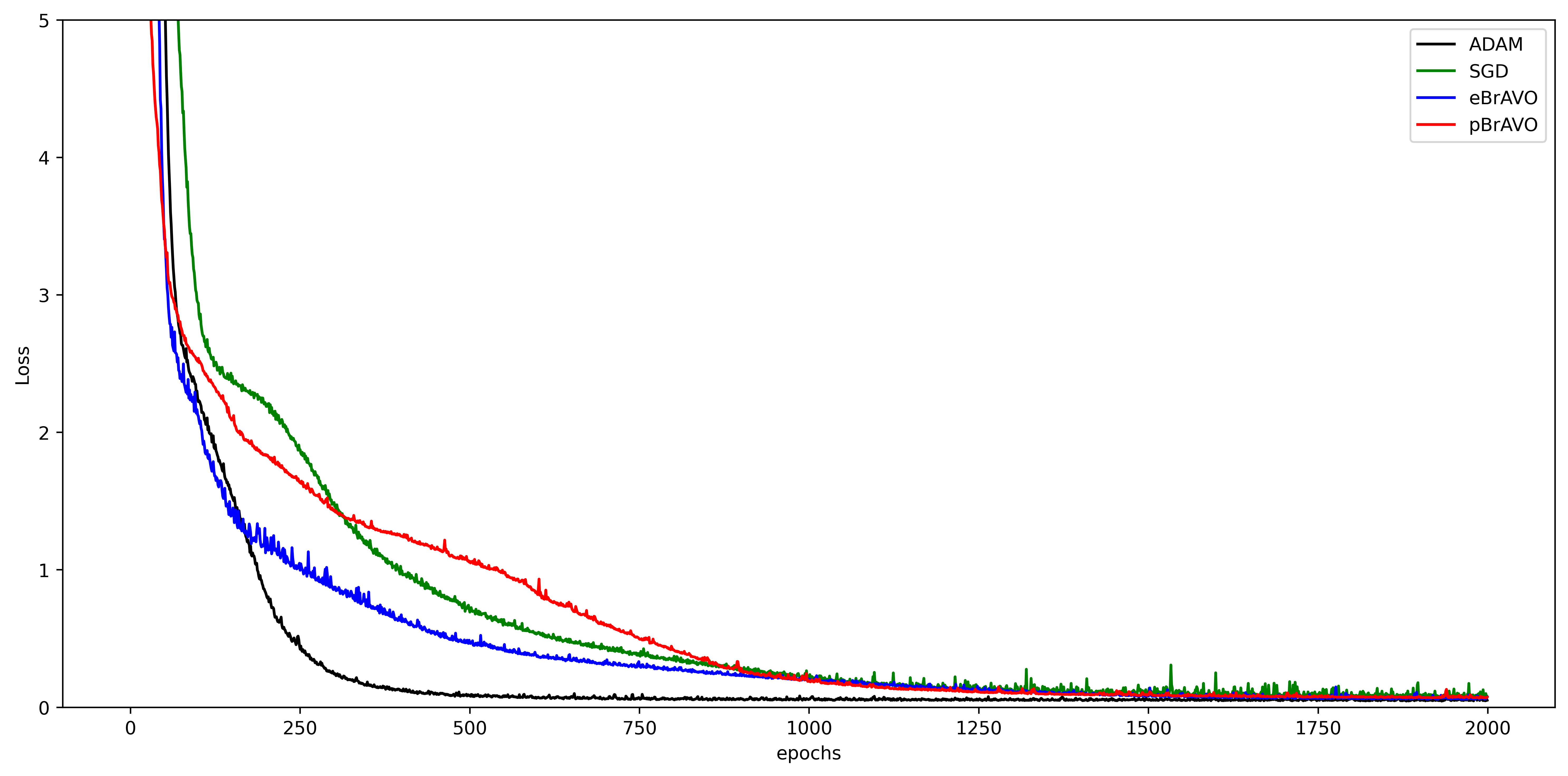} 
\vspace*{-5mm}
\caption{Evolution of the mean squared error for various algorithms, when applied to the problem of fitting a model to a set of 500 data points. \label{fig: DataFittingLoss} }
\end{figure}

\begin{figure}[!h]
	\vspace{-3mm}
\centering
\includegraphics[width=0.89\textwidth]{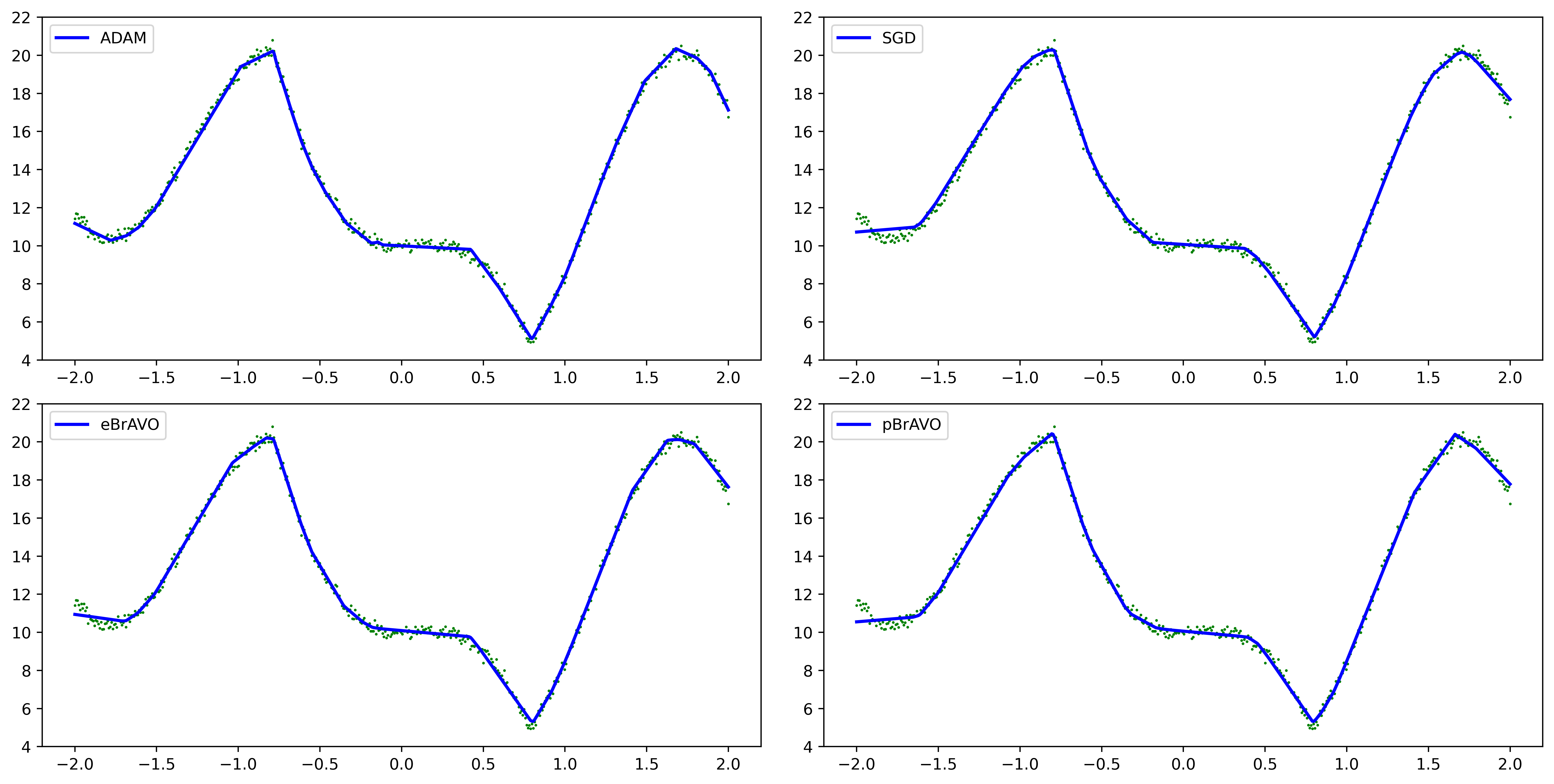} 	\vspace*{-3.7mm} 
\caption{Models obtained after 2000 epochs using various algorithms to fit the 500 data points displayed in green. \label{fig: DataFitting} }   \vspace{-1mm}
\end{figure}

Then, we solved a data fitting problem: given 500 data points from a noisy version of the function $10 x | \cos{2x} |+ 10 \exp{(-\sin{x})}$ on the interval $[-2,2]$, we wish to obtain a model which fits these data points as well as possible. We used the following neural network architecture (with 4,355 parameters) and loss function in TensorFlow:  \vspace{2.5mm}

\lstset{basicstyle=\ttfamily\footnotesize }
\begin{mdframed}[backgroundcolor=backcolour]\vspace{-2.5mm}\begin{lstlisting}[language=Python]
model = keras.Sequential()
model.add(layers.Dense(units = 1, activation = 'linear', input_shape=[1]))
model.add(layers.Dense(units = 64, activation = 'relu'))
model.add(layers.Dense(units = 64, activation = 'relu'))
model.add(layers.Dense(units = 1, activation = 'linear')) \end{lstlisting}  \end{mdframed} \vspace{3.5mm} 
The results are displayed in Figures~\ref{fig: DataFittingLoss} and~\ref{fig: DataFitting}. We see that all the algorithms achieve very small mean squared error, and all generate models, plotted as blue curves in Figure~\ref{fig: DataFitting}, which fit the green data points very well.  \\

\begin{figure}[!h]
\centering
\includegraphics[width=1\textwidth]{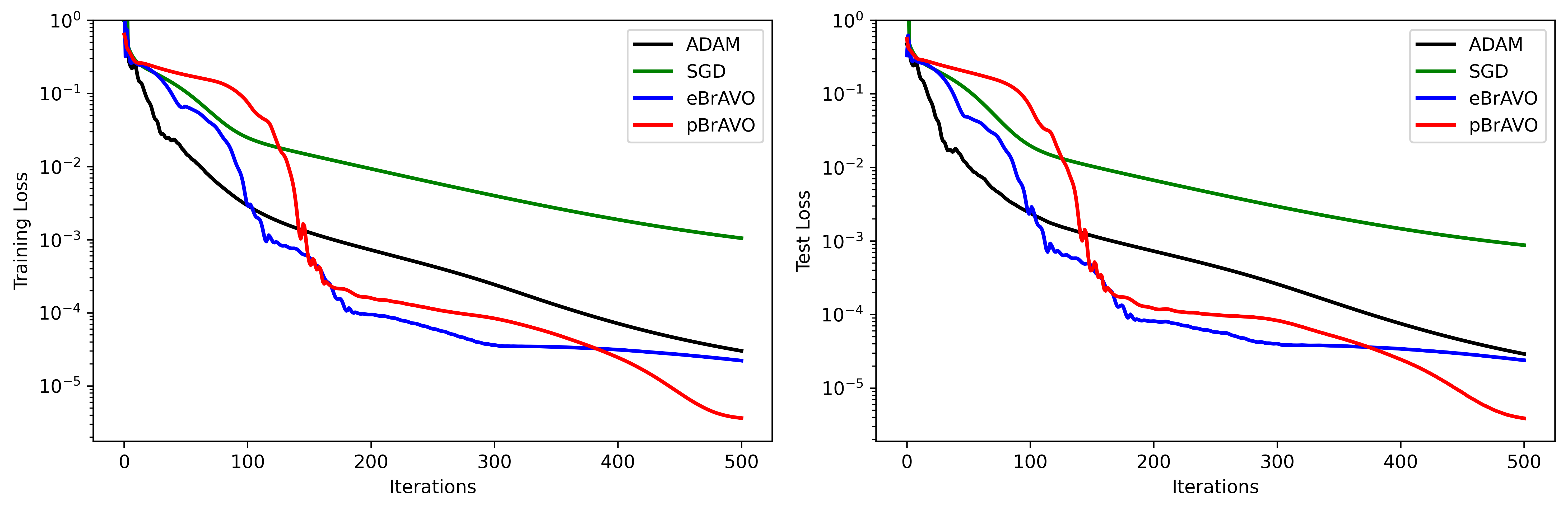} 
\vspace*{-7.5mm}
\caption{Evolution of the training and test losses for various algorithms, when learning the 231,310 parameters of a neural ODE network for dynamics learning. \label{fig: SO3Learning} } \vspace{-1mm}
\end{figure}

Finally, we test our algorithms for dynamics learning and control on the rotation group $SO(3)$. We consider the same problem as in~\cite{Duong21,DuruisseauxLieFVINs}, where we wish to learn the dynamics of a fully-actuated pendulum with dynamics given by $\ddot{\varphi} = -15\sin{\varphi} + 3u$, where $\varphi$ is the angle of the pendulum with respect to its vertically downward position and $u$ is a scalar control input. The data is collected from an OpenAI Gym environment, provided by \cite{Zhong2020}. We can see from Figure~\ref{fig: SO3Learning} that Adam and the BrAVO algorithms can achieve good training and test losses on this system identification problem using the Hamiltonian-based neural ODE network from~\cite{Duong21} (with 231,310 parameters), inspired by \cite{HNN2019,Zhong2020}. Note that we were unable to tune SGD to obtain a similar performance.

Overall, we have demonstrated that the BrAVO algorithms can be used conveniently within the PyTorch and TensorFlow frameworks, and that they can perform very well on more challenging optimization problems arising in machine learning applications, with a variety of model architectures, loss functions, and applications. We reiterate that this was the main purpose of this section, and that it is not our intention to make a very careful computational comparison of the BrAVO algorithms with other optimization algorithms that are commonly used by the machine learning community.   

A very careful computational comparison of optimization algorithms for machine learning is a much more ambitious goal which is beyond the scope of this paper. Such a comparison would be more meaningful once the current rudimentary implementation of the BrAVO algorithms within the PyTorch and TensorFlow frameworks has been highly-optimized, to take advantage of hardware architectures and highly-optimized PyTorch/TensorFlow operations. Aside from the quality of the implementation, other practical aspects of the algorithm could be investigated and improved further before carrying a careful comparison, for instance by looking into ways to boost the performance of the temporal looping technique or of the momentum restarting scheme.

An important advantage of our methods is that they are derived by discretizing continuous-time dynamical systems. We might be able to derive theoretical results about the algorithm by considering the associated continuous-time dynamical system and the discretization used. Furthermore, by considering the associated continuous-time dynamical system, we may be able to leverage numerous results from the theory of differential equations, dynamical systems, and geometric numerical integration. As a first example, in Section~\ref{section: Parameter C}, we exploited perturbation theory for continuous-time dynamical systems to gain insight into the effect of the parameter $C$ on the performance of the algorithms, which enabled us to improve tuning. As a second example, numerous ideas from the continuous-time theory of dynamical systems have been exploited in \cite{Alvarez2002,Attouch2020,Attouch2021,Attouch2022} and in particular the notions of dissipation, of viscous and Hessian-driven damping, and of inertia, in second-order differential equations. As a last example, the notion of momentum itself is better understood as a physical property of a continuous-time dynamical system, and we can also gain a lot of insight into the mechanism allowing the accelerated convergence towards the minimizer by considering these dynamical systems. There might be many other ways in which the performance of our algorithms can be improved by leveraging the associated continuous-time dynamical system.  \\

\section*{Conclusion and Future Directions}

	In this paper, we have discussed practical considerations which can significantly boost the computational performance and ease the tuning of symplectic accelerated optimization algorithms that are constructed by integrating Lagrangian and Hamiltonian systems coming from the variational framework for optimization introduced in~\citep{WiWiJo16}. 
	
	We showed that momentum restarting can lead to a significant gain in computational efficiency and robustness by reducing the undesirable effect of oscillations, and that a temporal looping strategy helps to avoid instability issues caused by numerical precision without impairing the computational performance of the algorithms. We also observed that time-adaptivity and the choice of symplectic integrator hardly make a difference once a momentum restarting scheme is incorporated in the optimization algorithms. This observation, along with other numerical experiments designed to study the effects of the different parameters, has provided insights that allowed to inform and ease the tuning process by simplifying the algorithms and by reducing the number of parameters to tune. 
	
	Overall, we have designed symplectic accelerated optimization algorithms whose computational efficiency and stability have been improved using temporal looping and momentum restarting, and which are now more user-friendly. We tested these algorithms on machine learning optimization problems with numerous different model architectures, loss functions, and applications, and saw that they can achieve very good results when tuned properly. %\\

	Preliminary experiments suggest that the benefits of momentum restarting and temporal looping uncovered in this paper extend to the Riemannian manifold framework for accelerated optimization introduced in~\cite{Duruisseaux2022Riemannian}. We intend to explore that direction to improve the computational efficiency and stability of symplectic accelerated optimization algorithms on Riemannian manifolds. 
	
	It would be nice to have further theoretical guarantees about the convergence of the discrete algorithm. However, this could be very difficult to obtain because momentum methods lack contraction, are nondescending, and are highly oscillatory~\cite{Orvieto2019}. While it is hoped that the continuous analysis will eventually guide the convergence analysis of the discrete-time algorithms, this does not appear to be a straightforward exercise, as one would first need to reconcile the arbitrarily fast rates of convergence of the continuous-time trajectories with Nesterov's barrier theorem of $\mathcal{O}(1/k^2)$ for discrete-time algorithms. We note however that some theoretical guarantees for certain integrators applied to the polynomial subfamily were obtained in the case where $p > 2$ in~\cite{Jadbabaie2018}, although this was already a very complicated task achieved under additional assumptions on the objective function and its derivatives. In the future, we intend to try to build upon the results of~\cite{Jadbabaie2018} to derive more general theoretical guarantees for our discrete algorithms, and see how momentum restarting and temporal looping affect those guarantees.

	The temporal looping technique could also be improved by designing different numerical instability criteria. Instead of temporal looping strategies, one could also try to implement very popular techniques in machine learning such as decaying learning rates via a learning rate scheduler, or to progressively increase the batch size~\cite{Smith2018}, or a combination of these different approaches.
	
	The current implementation of the algorithms within the PyTorch and TensorFlow frameworks is rather rudimentary, and can certainly be improved to reduce computational time by taking advantage of hardware architectures and highly-optimized PyTorch/TensorFlow operations. With the same objective in mind, one could also replace the gradient scheme for momentum restarting by the function scheme if the latter can be implemented more efficiently. 
	
	Once the algorithms have been improved further, possibly leveraging the theory of continuous-time dynamical systems, and once the implementation of the algorithms has been highly-optimized, it would be very interesting to perform a very careful computational comparison with other popular algorithms on many different types of problems to see whether the BrAVO algorithms can outperform the state-of-the-art algorithms on certain classes of machine learning problems.  \\

	\section*{Data Availability Statement}
 
Simple implementations of the optimization algorithms in MATLAB and Python, and more sophisticated Python code implementations which allow the optimizers to be called conveniently within the TensorFlow and PyTorch frameworks can be found at
\begin{framed}
\begin{center}
	\href{https://github.com/vduruiss/AccOpt_via_GNI}{github.com/vduruiss/AccOpt\_via\_GNI}
\end{center} 
\end{framed}

\hfill

\section*{Acknowledgments} 

%\hl{
The authors were supported in part by NSF under grants DMS-1411792, DMS-1345013, DMS-1813635, CCF-2112665, by AFOSR under grant FA9550-18-1-0288, and by the DoD under grant HQ00342010023 (Newton Award for Transformative Ideas during the COVID-19 Pandemic).

\newpage

\appendix

\section{List of Time-Adaptive Algorithms} \label{appendix: List of Time-Adaptive Algorithms}
\hfill 

\begin{minipage}{0.4\textwidth}
\centering	\textbf{\underline{PolyHTVI}}
	
	\begin{equation*}
		\begin{aligned}
			\mathfrak{q}_{k+1}  & = \mathfrak{q}_k + h \frac{p}{\mathring{p}} \mathfrak{q}_k^{1-\mathring{p}/p} \\
			r_{k+1} & = r_k -  \frac{p^2}{\mathring{p}} Ch \mathfrak{q}_k^{2p-\mathring{p}/p} \nabla f(q_k) \\
			q_{k+1} &= q_k +  \frac{p^2}{\mathring{p} } h\mathfrak{q}_k^{-p-\mathring{p}/p} r_{k+1}
		\end{aligned}
	\end{equation*} 
	
	\hfill

	\textbf{\underline{PolySLC}}
	
	\begin{equation*}
		\begin{aligned}
									r & \leftarrow  r - \frac{Cp^2}{2\mathring{p}} h \mathfrak{q}^{2p-\mathring{p}/p } \nabla f(q) \\
			\mathfrak{q} & \leftarrow  \left( \mathfrak{q}^{\mathring{p}/p} + \frac{h}{2} \right)^{p/\mathring{p}} \\
			q &\leftarrow  q + \frac{hp^2}{\mathring{p} \mathfrak{q}^{p+\mathring{p}/p}}  r \\
			\mathfrak{q} & \leftarrow  \left( \mathfrak{q}^{\mathring{p}/p} + \frac{h}{2} \right)^{p/\mathring{p}}\\
				r & \leftarrow  r - \frac{Cp^2}{2\mathring{p}} h \mathfrak{q}^{2p-\mathring{p}/p } \nabla f(q) 
		\end{aligned} 
	\end{equation*}

\end{minipage}
\hfill
\begin{minipage}{0.52\textwidth}
\centering	\textbf{\underline{PolyLTVI}}
	
	\begin{equation*}
		\begin{aligned}   \mathfrak{q}_{k+1}  & = \mathfrak{q}_k + h \frac{p}{\mathring{p}} \mathfrak{q}_k^{1-\mathring{p}/p}  \\ q_{k+1} & = q_k +  \frac{ hp^3  }{ \mathring{p} ^2\mathfrak{q}_k^{p -1 + 2 \mathring{p} / p}}  r_k - \frac{ Ch^2p^4  }{ \mathring{p} ^2} \mathfrak{q}_k^{p - 2 \mathring{p} / p } \nabla f(q_k)   \\  r_{k+1}     & =   \frac{\mathring{p} ^2\mathfrak{q}_k^{p  +  \mathring{p} / p}}{hp^3 \mathfrak{q}_{k+1}^{1-\mathring{p}/p} } (q_{k+1}-q_k)      
		\end{aligned} 
	\end{equation*}
	
	\hfill \\

	\textbf{\underline{PolySV}} 
	
	\begin{equation*}
		\begin{aligned}
			r_{k+\frac{1}{2}}&=r_k- \frac{p^2}{2\mathring{p}} Ch \mathfrak{q}_k^{2p-\mathring{p}/p} \nabla f(q_k) \\
			\text{Solve} &	\text{ }\text{ } \mathfrak{q}_{k+1}   = \mathfrak{q}_k  +   \frac{hp}{2\mathring{p}} \left( \mathfrak{q}_k^{1-\mathring{p}/p}  + \mathfrak{q}_{k+1}^{1-\mathring{p}/p}\right)   \\ q_{k+1}&=q_k+\frac{hp^2}{2 \mathring{p}} \left( \mathfrak{q}_k^{-p-\mathring{p}/p} +  \mathfrak{q}_{k+1}^{-p-\mathring{p}/p}    \right) r_{k+\frac{1}{2}}   \\
			r_{k+1}&= r_{k+\frac{1}{2}} - \frac{p^2}{2\mathring{p}} Ch \mathfrak{q}_{k+1}^{2p-\mathring{p}/p} \nabla f(q_{k+1}) 
		\end{aligned} 
	\end{equation*}

\end{minipage}

\hfill  

\hrulefill 

\hfill

\begin{minipage}{0.4\textwidth}
\centering	\textbf{\underline{ExpoHTVI}}
	
	\begin{equation*}
		\begin{aligned}
			\mathfrak{q}_{k+1}  & = \mathfrak{q}_k +  \frac{\eta}{\mathring{\eta}} h  \\
			r_{k+1} & = r_k -  \frac{\eta^2}{\mathring{\eta}} Ch e^{2 \eta \mathfrak{q}_k} \nabla f(q_k) \\
			q_{k+1} &= q_k +  \frac{\eta^2}{\mathring{\eta} } he^{-\eta \mathfrak{q}_k} r_{k+1}
		\end{aligned}
	\end{equation*} 
	
	\hfill \\

	\textbf{\underline{ExpoSLC}}
	
	\begin{equation*}
		\begin{aligned}
					r & \leftarrow  r -  \frac{Ch\eta^2}{2\mathring{\eta}} e^{2\eta \mathfrak{q}} \nabla f(q)  \\
			\mathfrak{q} &\leftarrow  \mathfrak{q} + \frac{h\eta}{2\mathring{\eta}} \\
			q & \leftarrow  q + 	\frac{h\eta^2}{ \mathring{\eta}e^{\eta \mathfrak{q}}}  r   \\
			\mathfrak{q} &\leftarrow  \mathfrak{q} + \frac{h\eta}{2\mathring{\eta}}  \\
						r & \leftarrow  r -  \frac{Ch\eta^2}{2\mathring{\eta}} e^{2\eta \mathfrak{q}} \nabla f(q)  \\
		\end{aligned} 
	\end{equation*}

\end{minipage}
\hfill
\begin{minipage}{0.52\textwidth}
\centering	\textbf{\underline{ExpoLTVI}}
	
	\begin{equation*}
		\begin{aligned}   \mathfrak{q}_{k+1}  & = \mathfrak{q}_k + \frac{\eta}{\mathring{\eta}}h   \\ q_{k+1} & = q_k +  \frac{h \eta^3  }{ \mathring{\eta} ^2}  e^{-\eta \mathfrak{q}_k}  r_k - \frac{ Ch^2\eta^4  }{ \mathring{\eta} ^2} e^{\eta \mathfrak{q}_k} \nabla f(q_k)  \\  r_{k+1}     & =   \frac{\mathring{\eta}^2}{h\eta^3} e^{\eta \mathfrak{q}_k} (q_{k+1} - q_k)  
		\end{aligned} 
	\end{equation*}
	
	\hfill \\
	
	\hfill 
	
	\textbf{\underline{ExpoSV} } 
	
	\begin{equation*}
		\begin{aligned}
			r_{k+\frac{1}{2}}&=r_k- \frac{\eta^2}{\mathring{2\eta}} Ch e^{2 \eta \mathfrak{q}_k} \nabla f(q_k) \\
			\mathfrak{q}_{k+1}  & = \mathfrak{q}_k +   \frac{\eta}{\mathring{\eta}}h   \\ q_{k+1}&=q_k+\frac{h\eta^2}{2\mathring{\eta} } \left(    e^{-\eta \mathfrak{q}_{k+1}} + e^{-\eta \mathfrak{q}_k} \right) r_{k+\frac{1}{2}} \\
			r_{k+1}&=r_{k+\frac{1}{2}}- \frac{\eta^2}{\mathring{2\eta}} Ch e^{2 \eta \mathfrak{q}_{k+1}} \nabla f(q_{k+1})
		\end{aligned} 
	\end{equation*}
\end{minipage}

\newpage

\hfill 

\begin{minipage}{0.4\textwidth}
\centering	\textbf{\underline{ExpoToPolyHTVI}}
	
	\begin{equation*}
		\begin{aligned}
			\mathfrak{q}_{k+1}  & = \left( 1+ \frac{\eta h}{p} \right)\mathfrak{q}_k  \\
			r_{k+1} & = r_k -   \frac{\eta^2}{p} Ch \mathfrak{q}_k e^{2\eta \mathfrak{q}_k}  \nabla f(q_k)\\
			q_{k+1} &= q_k +  \frac{ h \eta^2}{ p e^{ \eta \mathfrak{q}_k }  } \mathfrak{q}_k r_{k+1}
		\end{aligned}
	\end{equation*} 
	
	\hfill \\ 
	
	\hfill

	\textbf{\underline{ExpoToPolySLC}}
	
	\begin{equation*}
		\begin{aligned}
						r & \leftarrow  r -  \frac{Ch \mathfrak{q}\eta^2}{2p} e^{2\eta \mathfrak{q}} \nabla f(q)   \\
			\mathfrak{q} & \leftarrow  \mathfrak{q} e^{\frac{\eta h}{ 2p}} \\
			q & \leftarrow  q + 	\frac{h \mathfrak{q} \eta^2}{ p e^{\eta \mathfrak{q}}}  r  \\
			\mathfrak{q} & \leftarrow  \mathfrak{q} e^{\frac{\eta h}{ 2p}}  \\
						r & \leftarrow  r -  \frac{Ch \mathfrak{q}\eta^2}{2p} e^{2\eta \mathfrak{q}} \nabla f(q)   \\
		\end{aligned} 
	\end{equation*}

\end{minipage}
\hfill
\begin{minipage}{0.52\textwidth}
\centering	\textbf{\underline{ExpoToPolyLTVI}}
	
	\begin{equation*}
		\begin{aligned}   \mathfrak{q}_{k+1}  & = \left( 1+ \frac{\eta h}{p} \right)\mathfrak{q}_k , \\ q_{k+1} & = q_k +  \frac{h  \mathfrak{q}_k^2 \eta^3 }{ p^2 e^{\eta \mathfrak{q}_k} }  r_k - \frac{ Ch^2\eta^4  }{ p^2}\mathfrak{q}_k^2 e^{\eta \mathfrak{q}_k} \nabla f(q_k),  \\  r_{k+1}     & = \frac{ p(p+\eta h)   }{h\eta^3 \mathfrak{q}_k^2}  e^{\eta \mathfrak{q}_k} (q_{k+1} - q_k)  . 
		\end{aligned} 
	\end{equation*}
	
	\hfill \\
	
	\hfill 
	
	\textbf{\underline{ExpoToPolySV} } 
	
	\begin{equation*}
		\begin{aligned}
		r_{k+\frac{1}{2}}&=r_k-  \frac{ \eta^2}{2p} C h \mathfrak{q}_k e^{2\eta \mathfrak{q}_k}  \nabla f(q_k) ,\\
			\mathfrak{q}_{k+1}  & = \frac{2p+\eta h}{2p-\eta h}\mathfrak{q}_k  ,  \\ q_{k+1}&= q_k+\frac{h\eta^2}{2p } \left(   \mathfrak{q}_k  e^{-\eta \mathfrak{q}_k}  + \mathfrak{q}_{k+1}  e^{-\eta \mathfrak{q}_{k+1}}   \right) r_{k+\frac{1}{2}} ,  \\
			r_{k+1}&=r_{k+\frac{1}{2}} -   \frac{ \eta^2}{2p} C h \mathfrak{q}_{k+1} e^{2\eta \mathfrak{q}_{k+1}}  \nabla f(q_{k+1}),
		\end{aligned} 
	\end{equation*}

\end{minipage}

\hfill  

\hrulefill

\hfill

\begin{minipage}{0.4\textwidth}
\centering	\textbf{\underline{PolyToExpoHTVI}}
	
	\begin{equation*}
		\begin{aligned}
			\mathfrak{q}_{k+1}  & = \mathfrak{q}_k   + h\frac{p}{\eta} e^{-\frac{\eta}{p} \mathfrak{q}_k }  \\
			r_{k+1} & = r_k -   \frac{Chp^2}{\eta e^{ \frac{\eta}{p} \mathfrak{q}_k} }  \mathfrak{q}_k^{2p-1}   \nabla f(q_k) \\
			q_{k+1} &= q_k + h 	  \frac{p^2 }{  \eta \mathfrak{q}_k^{p+1}  } e^{-\frac{\eta }{p} \mathfrak{q}_k} r_{k+1}
		\end{aligned}
	\end{equation*} 
	
	\hfill \\

	\hfill 
	
	\textbf{\underline{PolyToExpoSLC} }
	
	\begin{equation*}
		\begin{aligned}
				r & \leftarrow  r - h  \frac{Cp^2}{2\eta} \mathfrak{q}^{2p-1} e^{-\frac{\eta }{p} \mathfrak{q}}  \nabla f(q)   \\
					\mathfrak{q} &\leftarrow  \frac{p}{\eta}  \log{\left( e^{\frac{\eta }{p} \mathfrak{q}} +\frac{h}{2} \right)}  \\
			q & \leftarrow  q + 	\frac{hp^2  }{ \eta \mathfrak{q}^{p+1}}  e^{-\frac{\eta }{p} \mathfrak{q}}   r  \\
			\mathfrak{q} &\leftarrow  \frac{p}{\eta}  \log{\left( e^{\frac{\eta }{p} \mathfrak{q}} +\frac{h}{2} \right)}   \\
						r & \leftarrow  r - h  \frac{Cp^2}{2\eta} \mathfrak{q}^{2p-1} e^{-\frac{\eta }{p} \mathfrak{q}}  \nabla f(q)   
		\end{aligned} 
	\end{equation*}

\end{minipage}
\hfill
\begin{minipage}{0.52\textwidth}
\centering	\textbf{\underline{PolyToExpoLTVI} }
	
	\begin{equation*}
		\begin{aligned}   \mathfrak{q}_{k+1}  & = \mathfrak{q}_k   + h\frac{p}{\eta} e^{-\frac{\eta}{p} \mathfrak{q}_k } \\ q_{k+1} & = q_k +  \frac{ hp^3  }{   \eta^2 \mathfrak{q}_k^{p+1} e^{\frac{2\eta}{p} \mathfrak{q}_k}  }   r_k - \frac{ Ch^2p^4  }{ \eta^2 e^{\frac{2\eta}{p} \mathfrak{q}_k}  }\mathfrak{q}_k^{p-2}  \nabla f(q_k) \\  r_{k+1}     & = \frac{ \eta^2 \mathfrak{q}_k^{p+1}  }{h p^3}  e^{\frac{\eta}{p} (\mathfrak{q}_{k+1} + \mathfrak{q}_k)} (q_{k+1} - q_k)  
		\end{aligned} 
	\end{equation*}
	
	\hfill \\
	
	\hfill \\
	
	\hfill 
	
	\textbf{\underline{PolyToExpoSV} } 
	
	\begin{equation*}
		\begin{aligned}
		r_{k+\frac{1}{2}} &=r_k-  \frac{p^2}{2\eta} Ch \mathfrak{q}_k^{2p-1}  e^{- \frac{\eta}{p} \mathfrak{q}_k}  \nabla f(q_k) \\
			\text{Solve} &	\text{ }\text{ }	\mathfrak{q}_{k+1}    = \mathfrak{q}_k + \frac{ hp}{2\eta} \left(  e^{-\frac{\eta}{p} \mathfrak{q}_k } + e^{-\frac{\eta}{p} \mathfrak{q}_{k+1} }   \right)     \\ q_{k+1}&= q_k+\frac{hp^2}{2\eta} \left(    \mathfrak{q}_k^{-p-1} e^{-\frac{\eta }{p} \mathfrak{q}_k} +    \mathfrak{q}_{k+1}^{-p-1} e^{-\frac{\eta }{p} \mathfrak{q}_{k+1}}  \right)  r_{k+\frac{1}{2}}   \\
			r_{k+1}&=r_{k+\frac{1}{2}} -  \frac{p^2}{2\eta} Ch \mathfrak{q}_{k+1}^{2p-1}  e^{- \frac{\eta}{p} \mathfrak{q}_{k+1}}  \nabla f(q_{k+1})
		\end{aligned} 
	\end{equation*}

\end{minipage}

\newpage

\hfill 

\section{List of Non-Adaptive Algorithms} \label{appendix: List of Non-Adaptive Algorithms}

\hfill 

\hfill

\begin{minipage}{0.4\textwidth}
	\centering	\textbf{\underline{PolyHTVI}}
	
	\begin{equation*}
		\begin{aligned}
			r_{k+1} & = r_k -   Chp \mathfrak{q}_k^{2p-1} \nabla f(q_k) \\
			q_{k+1} &= q_k +  hp \mathfrak{q}_k^{-p-1} r_{k+1} \\
						\mathfrak{q}_{k+1}  & = \mathfrak{q}_k + h
		\end{aligned}
	\end{equation*} 
	
	\hfill \\
	
	\hfill \\

	\textbf{\underline{PolySLC}}
	
	\begin{equation*}
		\begin{aligned}
						r & \leftarrow  r - \frac{1}{2} Chp \mathfrak{q}^{2p-1 } \nabla f(q) \\
			\mathfrak{q} & \leftarrow  \mathfrak{q} +  \frac{h}{2} \\
			q &\leftarrow  q + hp \mathfrak{q}^{-p-1}  r \\
			\mathfrak{q} & \leftarrow  \mathfrak{q} +  \frac{h}{2} \\
				r & \leftarrow  r - \frac{1}{2} Chp \mathfrak{q}^{2p-1 } \nabla f(q) \\
		\end{aligned} 
	\end{equation*}

\end{minipage}
\hfill
\begin{minipage}{0.52\textwidth}
	\centering	\textbf{\underline{PolyLTVI}}
	
	\begin{equation*}
	\begin{aligned}  
	  q_{k+1} & = q_k + hp \mathfrak{q}_k^{-p-1}  r_k - Ch^2p^2 \mathfrak{q}_k^{p - 2} \nabla f(q_k)   \\  r_{k+1}     & =   \frac{\mathfrak{q}_k^{p 1}}{hp } (q_{k+1}-q_k)       \\
	  \mathfrak{q}_{k+1}  & = \mathfrak{q}_k + h  
	\end{aligned} 
\end{equation*}

	\hfill \\
	
	\hfill \\

	\textbf{\underline{PolySV}} 
	
	\begin{equation*}
		\begin{aligned}
			r_{k+\frac{1}{2}}&=r_k - \frac{1}{2} Chp \mathfrak{q}_k^{2p-1} \nabla f(q_k) \\
 \mathfrak{q}_{k+1}   & = \mathfrak{q}_k  +  h   \\ q_{k+1}&=q_k+\frac{h}{2} p \left( \mathfrak{q}_k^{-p-1} +  \mathfrak{q}_{k+1}^{-p-1}    \right) r_{k+\frac{1}{2}}   \\
			r_{k+1}&= r_{k+\frac{1}{2}} - \frac{1}{2} Chp \mathfrak{q}_{k+1}^{2p-1} \nabla f(q_{k+1}) 
		\end{aligned} 
	\end{equation*}
	
	\hfill 
	
	\hfill 
	
	\hfill

\end{minipage}

\hfill   \\

\hrulefill 

\hfill   \\

\begin{minipage}{0.4\textwidth}
	\centering	\textbf{\underline{ExpoHTVI}}
	
	\begin{equation*}
		\begin{aligned}
			r_{k+1} & = r_k -   C\eta h e^{2 \eta \mathfrak{q}_k} \nabla f(q_k) \\
			q_{k+1} &= q_k +  \eta he^{-\eta \mathfrak{q}_k} r_{k+1} \\
						\mathfrak{q}_{k+1}  & = \mathfrak{q}_k +  h  
		\end{aligned}
	\end{equation*} 

	\hfill \\

\hfill \\

	\textbf{\underline{ExpoSLC}}
	
	\begin{equation*}
		\begin{aligned}
						r & \leftarrow  r -  \frac{1}{2} C\eta h   e^{2\eta \mathfrak{q}} \nabla f(q)  \\
			\mathfrak{q} &\leftarrow  \mathfrak{q} + \frac{h}{2} \\
			q & \leftarrow  q + 	\eta h e^{-\eta \mathfrak{q}}  r   \\
						\mathfrak{q} &\leftarrow  \mathfrak{q} + \frac{h}{2} \\
					r & \leftarrow  r -  \frac{1}{2} C\eta h   e^{2\eta \mathfrak{q}} \nabla f(q)  \\
		\end{aligned} 
	\end{equation*}

\end{minipage}
\hfill
\begin{minipage}{0.52\textwidth}
	\centering	\textbf{\underline{ExpoLTVI}}
	
	\begin{equation*}
		\begin{aligned}    q_{k+1} & = q_k + h\eta  e^{-\eta \mathfrak{q}_k}  r_k - C \eta^2 h^2 e^{\eta \mathfrak{q}_k} \nabla f(q_k)  \\  r_{k+1}     & =   \frac{e^{\eta \mathfrak{q}_k}}{\eta h}  (q_{k+1} - q_k)  \\
			\mathfrak{q}_{k+1}  & = \mathfrak{q}_k + h  
		\end{aligned} 
	\end{equation*}
	
	\hfill \\

\hfill \\
	
	\textbf{\underline{ExpoSV} } 
	
	\begin{equation*}
		\begin{aligned}
			r_{k+\frac{1}{2}}&=r_k- \frac{1}{2} C \eta h e^{2 \eta \mathfrak{q}_k} \nabla f(q_k) \\
			\mathfrak{q}_{k+1}  & = \mathfrak{q}_k +   h   \\ q_{k+1}&=q_k+\frac{1}{2} \eta h \left(    e^{-\eta \mathfrak{q}_{k+1}} + e^{-\eta \mathfrak{q}_k} \right) r_{k+\frac{1}{2}} \\
			r_{k+1}&=r_{k+\frac{1}{2}}- \frac{1}{2} C \eta h e^{2 \eta \mathfrak{q}_{k+1}} \nabla f(q_{k+1})
		\end{aligned} 
	\end{equation*}

	\hfill 

\hfill 

\hfill 

\end{minipage}

\newpage

\bibliography{PracticalOpt}

\def\cprime{$'$}
\begin{thebibliography}{91}
\providecommand{\natexlab}[1]{#1}
\providecommand{\url}[1]{\texttt{#1}}
\expandafter\ifx\csname urlstyle\endcsname\relax
  \providecommand{\doi}[1]{doi: #1}\else
  \providecommand{\doi}{doi: \begingroup \urlstyle{rm}\Url}\fi

\bibitem[Ahn and Sra(2020)]{Sra2020}
K.~Ahn and S.~Sra.
\newblock From {N}esterov's estimate sequence to {R}iemannian acceleration.
\newblock In \emph{Proceedings of Thirty Third Conference on Learning Theory},
  volume 125 of \emph{Proceedings of Machine Learning Research}, pages 84--118.
  PMLR, 09--12 Jul 2020.

\bibitem[Alecsa and László(2021)]{Alecsa2021}
C.~D. Alecsa and S.~C. László.
\newblock {T}ikhonov regularization of a perturbed heavy ball system with
  vanishing damping.
\newblock \emph{SIAM Journal on Optimization}, 31\penalty0 (4):\penalty0
  2921--2954, 2021.
\newblock \doi{10.1137/20M1382027}.

\bibitem[Alimisis et~al.(2020{\natexlab{a}})Alimisis, Orvieto, B\'ecigneul, and
  Lucchi]{Alimisis2020-1}
F.~Alimisis, A.~Orvieto, G.~B\'ecigneul, and A.~Lucchi.
\newblock Practical accelerated optimization on {R}iemannian manifolds.
\newblock 2020{\natexlab{a}}.
\newblock URL \url{https://arxiv.org/abs/2002.04144}.

\bibitem[Alimisis et~al.(2020{\natexlab{b}})Alimisis, Orvieto, B\'ecigneul, and
  Lucchi]{alimisis2020}
F.~Alimisis, A.~Orvieto, G.~B\'ecigneul, and A.~Lucchi.
\newblock A continuous-time perspective for modeling acceleration in
  {R}iemannian optimization.
\newblock In \emph{Proceedings of the 23rd International AISTATS Conference},
  volume 108 of \emph{PMLR}, pages 1297--1307, 2020{\natexlab{b}}.

\bibitem[Alimisis et~al.(2021)Alimisis, Orvieto, B{\'e}cigneul, and
  Lucchi]{Alimisis2021}
F.~Alimisis, A.~Orvieto, G.~B{\'e}cigneul, and A.~Lucchi.
\newblock Momentum improves optimization on {R}iemannian manifolds.
\newblock In \emph{AISTATS}, 2021.

\bibitem[Alvarez et~al.(2002)Alvarez, Attouch, Bolte, and Redont]{Alvarez2002}
F.~Alvarez, H.~Attouch, J.~Bolte, and P.~Redont.
\newblock A second-order gradient-like dissipative dynamical system with
  {H}essian-driven damping: Application to optimization and mechanics.
\newblock \emph{Journal de Mathématiques Pures et Appliquées}, 81\penalty0
  (8):\penalty0 747--779, 2002.
\newblock ISSN 0021-7824.
\newblock \doi{10.1016/S0021-7824(01)01253-3}.

\bibitem[Attouch and Chbani(2016)]{Attouch2016}
H.~Attouch and Z.~Chbani.
\newblock Combining fast inertial dynamics for convex optimization with
  {T}ikhonov regularization.
\newblock 2016.

\bibitem[Attouch and Czarnecki(2017)]{Attouch2017}
H.~Attouch and M.~Czarnecki.
\newblock Asymptotic behavior of gradient-like dynamical systems involving
  inertia and multiscale aspects.
\newblock \emph{Journal of Differential Equations}, 262\penalty0 (3):\penalty0
  2745--2770, 2017.
\newblock ISSN 0022-0396.
\newblock \doi{10.1016/j.jde.2016.11.009}.

\bibitem[Attouch et~al.(2020)Attouch, Chbani, Fadili, and Riahi]{Attouch2020}
H.~Attouch, Z.~Chbani, J.~Fadili, and H.~Riahi.
\newblock First-order optimization algorithms via inertial systems with
  {H}essian driven damping.
\newblock \emph{Mathematical Programming}, Nov 2020.
\newblock \doi{10.1007/s10107-020-01591-1}.

\bibitem[Attouch et~al.(2021)Attouch, Chbani, Fadili, and Riahi]{Attouch2021}
H.~Attouch, Z.~Chbani, J.~M. Fadili, and H.~Riahi.
\newblock {Convergence of iterates for first-order optimization algorithms with
  inertia and Hessian driven damping}.
\newblock \emph{Optimization}, 2021.
\newblock \doi{10.1080/02331934.2021.2009828}.

\bibitem[Attouch et~al.(2022)Attouch, Balhag, Chbani, and Riahi]{Attouch2022}
H.~Attouch, A.~Balhag, Z.~Chbani, and H.~Riahi.
\newblock Fast convex optimization via inertial dynamics combining viscous and
  {H}essian-driven damping with time rescaling.
\newblock \emph{Evolution Equations and Control Theory}, 11\penalty0
  (2):\penalty0 487--514, 2022.

\bibitem[Attri et~al.(2020)Attri, Sharma, Takach, and Shah]{KerasWeather}
P.~Attri, Y.~Sharma, K.~Takach, and F.~Shah.
\newblock Timeseries forecasting for weather prediction.
\newblock \emph{Keras Tutorial}, 2020.
\newblock URL
  \url{https://keras.io/examples/timeseries/timeseries_weather_forecasting/}.

\bibitem[Beck and Teboulle(2009)]{Beck2009}
A.~Beck and M.~Teboulle.
\newblock Gradient-based algorithms with applications to signal-recovery
  problems.
\newblock \emph{Convex Optimization in Signal Processing and Communications},
  pages 42--88, 2009.
\newblock \doi{10.1017/CBO9780511804458.003}.

\bibitem[Benettin(1994)]{Benettin1994}
A.~Benettin, G.and~Giorgilli.
\newblock On the {H}amiltonian interpolation of near-to-the identity symplectic
  mappings with application to symplectic integration algorithms.
\newblock \emph{Journal of Statistical Physics}, 74:\penalty0 1117--1143, 03
  1994.
\newblock \doi{10.1007/BF02188219}.

\bibitem[Bertsekas(2009)]{Bertsekas2009}
D.~Bertsekas.
\newblock \emph{Convex Optimization Algorithms}.
\newblock Athena Scientific, 2009.

\bibitem[Betancourt et~al.(2018)Betancourt, Jordan, and
  Wilson]{JordanSymplecticOptimization}
M.~Betancourt, {M. I.} Jordan, and A.~Wilson.
\newblock On symplectic optimization.
\newblock 2018.
\newblock URL \url{https://arxiv.org/abs/1802.03653}.

\bibitem[Blanes and Casas(2017)]{Blanes2017}
S.~Blanes and F.~Casas.
\newblock \emph{A Concise Introduction to Geometric Numerical Integration}.
\newblock 2017.
\newblock ISBN 9781482263442.
\newblock \doi{10.1201/b21563}.

\bibitem[Boltyanski et~al.(1999)Boltyanski, Martini, Soltan, and
  Soltan]{Boltyanski1999}
V.~Boltyanski, H.~Martini, V.~Soltan, and V.P. Soltan.
\newblock \emph{Geometric Methods and Optimization Problems}.
\newblock Combinatorial Optimization. Springer US, 1999.
\newblock \doi{10.1007/978-1-4615-5319-9}.

\bibitem[Boyd and Vandenberghe(2004)]{Boyd2004}
S.~Boyd and L.~Vandenberghe.
\newblock \emph{Convex Optimization}.
\newblock Cambridge University Press, 2004.
\newblock \doi{10.1017/CBO9780511804441}.

\bibitem[Calvo and Sanz-Serna(1993)]{CalSan93}
J.~P. Calvo and J.~M. Sanz-Serna.
\newblock The development of variable-step symplectic integrators, with
  application to the two-body problem.
\newblock \emph{SIAM J. Sci. Comp.}, 14\penalty0 (4):\penalty0 936--952, 1993.

\bibitem[Campos et~al.(2021)Campos, Mahillo, and de~Diego]{Campos2021}
C.~M. Campos, A.~Mahillo, and D.~Mart{\'i}n de~Diego.
\newblock A discrete variational derivation of accelerated methods in
  optimization.
\newblock 2021.
\newblock URL \url{https://arxiv.org/abs/2106.02700}.

\bibitem[Cauchy(1847)]{Cauchy1847}
A.~L. Cauchy.
\newblock M\'ethode g\'en\'erale pour la r\'esolution des syst\`emes
  d’\'equations simultan\'ees.
\newblock \emph{Acad. Sci. Paris}, 25:\penalty0 536--538, 1847.

\bibitem[Dai et~al.(2004)Dai, Liao, and Li]{Dai2004}
Y.-H. Dai, L.-Z. Liao, and D.~Li.
\newblock On restart procedures for the conjugate gradient method: Theory and
  practice in optimization.
\newblock \emph{Numerical Algorithms}, 35, 04 2004.
\newblock \doi{10.1023/B:NUMA.0000021761.10993.6e}.

\bibitem[Dea\~no et~al.(January 2018)Dea\~no, Huybrechs, and
  Iserles]{IserlesBook}
A.~Dea\~no, D.~Huybrechs, and A.~Iserles.
\newblock \emph{Computing Highly Oscillatory Integrals}.
\newblock SIAM, Philadelphia, January 2018.

\bibitem[Donghwan and Fessler(2018)]{Donghwan2018}
K.~Donghwan and J.~Fessler.
\newblock Adaptive restart of the optimized gradient method for convex
  optimization.
\newblock \emph{Journal of Optimization Theory and Applications}, 178, 07 2018.
\newblock \doi{10.1007/s10957-018-1287-4}.

\bibitem[Drezner and Hamacher(2002)]{Drezner2002}
Z.~Drezner and H.W. Hamacher.
\newblock \emph{Facility Location: Applications and Theory}.
\newblock Springer Berlin Heidelberg, 2002.
\newblock ISBN 9783540213451.

\bibitem[Duong and Atanasov(2021)]{Duong21}
T.~Duong and N.~Atanasov.
\newblock Hamiltonian-based neural {ODE} networks on the {SE(3)} manifold for
  dynamics learning and control.
\newblock In \emph{Proceedings of Robotics: Science and Systems}, July 2021.
\newblock \doi{10.15607/RSS.2021.XVII.086}.

\bibitem[Duruisseaux and Leok(2022{\natexlab{a}})]{Duruisseaux2022Constrained}
V.~Duruisseaux and M.~Leok.
\newblock Accelerated optimization on {R}iemannian manifolds via discrete
  constrained variational integrators.
\newblock \emph{Journal of Nonlinear Science}, 32\penalty0 (42),
  2022{\natexlab{a}}.
\newblock URL \url{https://doi.org/10.1007/s00332-022-09795-9}.

\bibitem[Duruisseaux and Leok(2022{\natexlab{b}})]{Duruisseaux2022Projection}
V.~Duruisseaux and M.~Leok.
\newblock Accelerated optimization on {R}iemannian manifolds via projected
  variational integrators.
\newblock 2022{\natexlab{b}}.
\newblock URL \url{https://arxiv.org/abs/2201.02904}.

\bibitem[Duruisseaux and Leok(2022{\natexlab{c}})]{Duruisseaux2022Riemannian}
V.~Duruisseaux and M.~Leok.
\newblock A variational formulation of accelerated optimization on {R}iemannian
  manifolds.
\newblock \emph{SIAM Journal on Mathematics of Data Science}, 4\penalty0
  (2):\penalty0 649--674, 2022{\natexlab{c}}.
\newblock URL \url{https://doi.org/10.1137/21M1395648}.

\bibitem[Duruisseaux and Leok(2023)]{Duruisseaux2022Lagrangian}
V.~Duruisseaux and M.~Leok.
\newblock Time-adaptive {L}agrangian variational integrators for accelerated
  optimization on manifolds.
\newblock \emph{Journal of Geometric Mechanics}, 15\penalty0 (1):\penalty0
  224--255, 2023.
\newblock ISSN 1941-4889.
\newblock URL \url{https://doi.org/10.3934/jgm.2023010}.

\bibitem[Duruisseaux et~al.(2021)Duruisseaux, Schmitt, and
  Leok]{duruisseaux2020adaptive}
V.~Duruisseaux, J.~Schmitt, and M.~Leok.
\newblock Adaptive {H}amiltonian variational integrators and applications to
  symplectic accelerated optimization.
\newblock \emph{SIAM Journal on Scientific Computing}, 43\penalty0
  (4):\penalty0 A2949--A2980, 2021.
\newblock URL \url{https://doi.org/10.1137/20M1383835}.

\bibitem[Duruisseaux et~al.(2023)Duruisseaux, Duong, Leok, and
  Atanasov]{DuruisseauxLieFVINs}
V.~Duruisseaux, T.~Duong, M.~Leok, and N.~Atanasov.
\newblock Lie group forced variational integrator networks for learning and
  control of robot systems.
\newblock \emph{5th Learning for Dynamics and Control Conference (L4DC)}, 2023.
\newblock URL \url{https://arxiv.org/pdf/2211.16006.pdf}.

\bibitem[Farr(2009)]{Farr2009}
W.~M. Farr.
\newblock Variational integrators for almost-integrable systems.
\newblock \emph{Celestial Mechanics and Dynamical Astronomy}, 102\penalty0
  (2):\penalty0 105--118, 2009.

\bibitem[Fercoq and Qu(2016)]{Fercoq2016}
O.~Fercoq and Z.~Qu.
\newblock {Restarting accelerated gradient methods with a rough strong
  convexity estimate}.
\newblock Research Report 1609.07358, {T{\'e}l{\'e}com ParisTech}, 2016.
\newblock URL \url{https://hal.telecom-paris.fr/hal-02287730}.

\bibitem[Fercoq and Qu(2019)]{Fercoq2019}
O.~Fercoq and Z.~Qu.
\newblock {Adaptive restart of accelerated gradient methods under local
  quadratic growth condition}.
\newblock \emph{{IMA Journal of Numerical Analysis}}, March 2019.
\newblock \doi{10.1093/imanum/drz007}.

\bibitem[Filon(1930)]{Filon1930}
L.~N.~G. Filon.
\newblock On a quadrature formula for trigonometric integrals.
\newblock \emph{Proceedings of the Royal Society of Edinburgh}, 49:\penalty0
  38–47, 1930.
\newblock \doi{10.1017/S0370164600026262}.

\bibitem[Giselsson and Boyd(2015)]{Giselsson2015}
P.~Giselsson and S.~Boyd.
\newblock Monotonicity and restart in fast gradient methods.
\newblock \emph{Proceedings of the IEEE Conference on Decision and Control},
  2015:\penalty0 5058--5063, 02 2015.
\newblock \doi{10.1109/CDC.2014.7040179}.

\bibitem[Gladman et~al.(1991)Gladman, Duncan, and Candy]{GlaDunCan91}
B.~Gladman, M.~Duncan, and J.~Candy.
\newblock Symplectic integrators for long-time integrations in celestial
  mechanics.
\newblock \emph{Celestial Mech. Dynamical Astronomy}, 52:\penalty0 221--240,
  1991.

\bibitem[Greydanus et~al.(2019)Greydanus, Dzamba, and Yosinski]{HNN2019}
S.~Greydanus, M.~Dzamba, and J.~Yosinski.
\newblock {Hamiltonian Neural Networks}.
\newblock In \emph{Advances in Neural Information Processing Systems},
  volume~32. Curran Associates, Inc., 2019.

\bibitem[Hairer(1997)]{Ha1997}
E.~Hairer.
\newblock Variable time step integration with symplectic methods.
\newblock \emph{Applied Numerical Mathematics}, 25\penalty0 (2-3):\penalty0
  219--227, 1997.

\bibitem[Hairer et~al.(2003)Hairer, Lubich, and Wanner]{Hairer2003}
E.~Hairer, C.~Lubich, and G.~Wanner.
\newblock Geometric numerical integration illustrated by the
  {S}t{\"o}rmer–{V}erlet method.
\newblock \emph{Acta Numerica}, 12:\penalty0 399 -- 450, 2003.

\bibitem[Hairer et~al.(2006)Hairer, Lubich, and Wanner]{HaLuWa2006}
E.~Hairer, C.~Lubich, and G.~Wanner.
\newblock \emph{Geometric {N}umerical {I}ntegration}, volume~31 of
  \emph{Springer Series in Computational Mathematics}.
\newblock Springer-Verlag, Berlin, 2nd edition, 2006.

\bibitem[Hall(2015)]{Hall2015}
B.~Hall.
\newblock \emph{Lie Groups, Lie Algebras, and Representations}.
\newblock Graduate Texts in Mathematics. Springer Cham, second edition, 2015.
\newblock \doi{10.1007/978-3-319-13467-3}.

\bibitem[Iserles and Quispel(2018)]{IserlesWhyGNI}
A.~Iserles and G.~R.~W. Quispel.
\newblock Why geometric numerical integration?
\newblock In Kurusch Ebrahimi-Fard and Mar{\'i}a Barbero~Li{\~{n}}{\'a}n,
  editors, \emph{Discrete Mechanics, Geometric Integration and Lie--Butcher
  Series}. Springer International Publishing, 2018.

\bibitem[Jendoubi and May(2010)]{Jendoubi2010}
M.~Jendoubi and R.~May.
\newblock On an asymptotically autonomous system with {T}ikhonov type
  regularizing term.
\newblock \emph{Archiv der Mathematik}, 95:\penalty0 389--399, 10 2010.
\newblock \doi{10.1007/s00013-010-0181-6}.

\bibitem[Jordan()]{Jordan2018}
M.~I. Jordan.
\newblock Dynamical, symplectic and stochastic perspectives on gradient-based
  optimization.
\newblock In \emph{Proceedings of the International Congress of Mathematicians
  (ICM 2018)}, pages 523--549.
\newblock \doi{10.1142/9789813272880_0022}.

\bibitem[Kingma and Ba(2014)]{ADAM}
D.~Kingma and J.~Ba.
\newblock Adam: A method for stochastic optimization.
\newblock \emph{International Conference on Learning Representations}, 12 2014.

\bibitem[Krizhevsky(2009)]{CIFAR}
A.~Krizhevsky.
\newblock Learning multiple layers of features from tiny images.
\newblock Technical report, University of Toronto, 2009.

\bibitem[Lall and West(2006)]{LaWe2006}
S.~Lall and M.~West.
\newblock Discrete variational {H}amiltonian mechanics.
\newblock \emph{J. Phys. A}, 39\penalty0 (19):\penalty0 5509--5519, 2006.

\bibitem[Lecun et~al.(1998)Lecun, Bottou, Bengio, and Haffner]{LeCun1998}
Y.~Lecun, L.~Bottou, Y.~Bengio, and P.~Haffner.
\newblock Gradient-based learning applied to document recognition.
\newblock \emph{Proceedings of the IEEE}, 86\penalty0 (11):\penalty0
  2278--2324, 1998.
\newblock \doi{10.1109/5.726791}.

\bibitem[Lee et~al.(2021)Lee, Tao, and Leok]{Lee2021}
T.~Lee, M.~Tao, and M.~Leok.
\newblock Variational symplectic accelerated optimization on {L}ie groups.
\newblock 2021.

\bibitem[Leimkuhler and Reich(2004)]{LeRe2005}
B.~Leimkuhler and S.~Reich.
\newblock \emph{Simulating {H}amiltonian Dynamics}, volume~14 of
  \emph{Cambridge Monographs on Applied and Computational Mathematics}.
\newblock Cambridge University Press, Cambridge, 2004.

\bibitem[Leok and Shingel(2012{\natexlab{a}})]{LeSh2011}
M.~Leok and T.~Shingel.
\newblock Prolongation-collocation variational integrators.
\newblock \emph{IMA J. Numer. Anal.}, 32\penalty0 (3):\penalty0 1194--1216,
  2012{\natexlab{a}}.

\bibitem[Leok and Shingel(2012{\natexlab{b}})]{LeSh2011_sbvi}
M.~Leok and T.~Shingel.
\newblock General techniques for constructing variational integrators.
\newblock \emph{Front. Math. China}, 7\penalty0 (2):\penalty0 273--303,
  2012{\natexlab{b}}.

\bibitem[Leok and Zhang(2011)]{LeZh2011}
M.~Leok and J.~Zhang.
\newblock Discrete {H}amiltonian variational integrators.
\newblock \emph{IMA Journal of Numerical Analysis}, 31\penalty0 (4):\penalty0
  1497--1532, 2011.

\bibitem[Levin(1982)]{Levin1982}
D.~Levin.
\newblock Procedures for computing one- and two-dimensional integrals of
  functions with rapid irregular oscillations.
\newblock \emph{Mathematics of Computation}, 38\penalty0 (158):\penalty0
  531--538, 1982.

\bibitem[Levin(1996)]{Levin1996}
D.~Levin.
\newblock Fast integration of rapidly oscillatory functions.
\newblock \emph{Journal of Computational and Applied Mathematics}, 67\penalty0
  (1):\penalty0 95--101, 1996.
\newblock ISSN 0377-0427.
\newblock \doi{10.1016/0377-0427(94)00118-9}.

\bibitem[Liu et~al.(2017)Liu, Shang, Cheng, Cheng, and Jiao]{Liu2017}
Y.~Liu, F.~Shang, J.~Cheng, H.~Cheng, and L.~Jiao.
\newblock Accelerated first-order methods for geodesically convex optimization
  on {R}iemannian manifolds.
\newblock In \emph{NeurIPS}, volume~30, pages 4868--4877, 2017.

\bibitem[Marsden and West(2001)]{MaWe2001}
J.~E. Marsden and M.~West.
\newblock Discrete mechanics and variational integrators.
\newblock \emph{Acta Numer.}, 10:\penalty0 357--514, 2001.

\bibitem[Muehlebach and Jordan(2019)]{Muehlebach2019}
M.~Muehlebach and M.~I. Jordan.
\newblock A dynamical systems perspective on {N}esterov acceleration.
\newblock In \emph{Proceedings of the 36th International Conference on Machine
  Learning}, volume~97 of \emph{PMLR}, Long Beach, CA, USA, 2019.

\bibitem[Nemirovsky and Yudin(1983)]{Nem1983}
{A. S.} Nemirovsky and {D. B.} Yudin.
\newblock \emph{Problem Complexity and Method Efficiency in Optimization}.
\newblock Wiley - Interscience series in discrete mathematics. Wiley, 1983.

\bibitem[Nesterov(1983)]{Nes83}
Y.~Nesterov.
\newblock A method of solving a convex programming problem with convergence
  rate $\mathcal{O}(1/k^2)$.
\newblock \emph{Soviet Mathematics Doklady}, 27\penalty0 (2):\penalty0
  372--376, 1983.

\bibitem[Nesterov(2004)]{Nes04}
Y.~Nesterov.
\newblock \emph{Introductory Lectures on Convex Optimization: A Basic Course},
  volume~87 of \emph{Applied Optimization}.
\newblock Kluwer Academic Publishers, Boston, MA, 2004.

\bibitem[Nesterov(2008)]{Nes08}
Y.~Nesterov.
\newblock Accelerating the cubic regularization of {N}ewton's method on convex
  problems.
\newblock \emph{Math. Program.}, 112:\penalty0 159--181, 2008.

\bibitem[O'donoghue and Cand\`{e}s(2015)]{O'donoghue2015}
B.~O'donoghue and E.~Cand\`{e}s.
\newblock Adaptive restart for accelerated gradient schemes.
\newblock \emph{Found. Comput. Math.}, 15\penalty0 (3):\penalty0 715–732, jun
  2015.
\newblock ISSN 1615-3375.
\newblock \doi{10.1007/s10208-013-9150-3}.

\bibitem[Orvieto and Lucchi(2019)]{Orvieto2019}
A.~Orvieto and A.~Lucchi.
\newblock Shadowing properties of optimization algorithms.
\newblock In \emph{Advances in Neural Information Processing Systems},
  volume~32, pages 12692--12703, 2019.

\bibitem[Paul and Rakshit(2020)]{KerasNLP}
S.~Paul and S.~Rakshit.
\newblock Large-scale multi-label text classification.
\newblock \emph{Keras Tutorial}, 2020.
\newblock URL \url{https://keras.io/examples/nlp/multi_label_classification/}.

\bibitem[Phillips(1962)]{Philips1962}
D.~L. Phillips.
\newblock A technique for the numerical solution of certain integral equations
  of the first kind.
\newblock \emph{J. ACM}, 9\penalty0 (1):\penalty0 84–97, jan 1962.
\newblock ISSN 0004-5411.
\newblock \doi{10.1145/321105.321114}.

\bibitem[Powell(1977)]{Powell1977}
M.~J.~D. Powell.
\newblock Restart procedures for the conjugate gradient method.
\newblock \emph{Mathematical Programming}, 12:\penalty0 241--254, 1977.

\bibitem[Renegar and Grimmer(2022)]{Renegar2022}
J.~Renegar and B.~Grimmer.
\newblock A simple nearly optimal restart scheme for speeding up first-order
  methods.
\newblock \emph{Found. Comput. Math.}, 22\penalty0 (1):\penalty0 211–256, feb
  2022.
\newblock ISSN 1615-3375.
\newblock \doi{10.1007/s10208-021-09502-2}.

\bibitem[Roulet and d'Aspremont(2020)]{Roulet2020}
V.~Roulet and A.~d'Aspremont.
\newblock Sharpness, restart, and acceleration.
\newblock \emph{SIAM Journal on Optimization}, 30\penalty0 (1):\penalty0
  262--289, 2020.
\newblock \doi{10.1137/18M1224568}.

\bibitem[Sanders et~al.(2007)Sanders, Verhulst, and Murdock]{Sanders2007}
J.~A. Sanders, F.~Verhulst, and J.~Murdock.
\newblock \emph{Averaging Methods in Nonlinear Dynamical Systems}.
\newblock Applied Mathematical Sciences. Springer New York, 2007.
\newblock ISBN 9780387489186.

\bibitem[Schmitt and Leok(2017)]{ScLe2017}
J.~M. Schmitt and M.~Leok.
\newblock Properties of {H}amiltonian variational integrators.
\newblock \emph{IMA Journal of Numerical Analysis}, 38\penalty0 (1):\penalty0
  377--398, 03 2017.

\bibitem[Schmitt et~al.(2018)Schmitt, Shingel, and Leok]{ScShLe2017}
J.~M. Schmitt, T.~Shingel, and M.~Leok.
\newblock {L}agrangian and {H}amiltonian {T}aylor variational integrators.
\newblock \emph{{BIT} Numerical Mathematics}, 58:\penalty0 457--488, 2018.
\newblock \doi{10.1007/s10543-017-0690-9}.

\bibitem[Smith et~al.(2018)Smith, Kindermans, Ying, and Le]{Smith2018}
S.~Smith, P.~Kindermans, C.~Ying, and Q.~V. Le.
\newblock Don't decay the learning rate, increase the batch size.
\newblock 2018.

\bibitem[Su et~al.(2016)Su, Boyd, and Candes]{SuBoCa16}
W.~Su, S.~Boyd, and E.~Candes.
\newblock A differential equation for modeling {N}esterov's {A}ccelerated
  {G}radient method: theory and insights.
\newblock \emph{Journal of Machine Learning Research}, 17\penalty0
  (153):\penalty0 1--43, 2016.

\bibitem[Tao and Ohsawa(2020)]{Tao2020}
M.~Tao and T.~Ohsawa.
\newblock Variational optimization on {L}ie groups, with examples of leading
  (generalized) eigenvalue problems.
\newblock In \emph{Proceedings of the 23rd International AISTATS Conference},
  volume 108 of \emph{PMLR}, 2020.

\bibitem[Tibshirani(1996)]{Tibshirani1996}
R.~Tibshirani.
\newblock Regression shrinkage and selection via the lasso.
\newblock \emph{Journal of the Royal Statistical Society. Series B
  (Methodological)}, 58\penalty0 (1):\penalty0 267--288, 1996.
\newblock ISSN 00359246.

\bibitem[Tikhonov(1963)]{Tikhonov1963}
A.~N. Tikhonov.
\newblock {Solution of incorrectly formulated problems and the regularization
  method}.
\newblock \emph{{Sov. Math., Dokl.}}, 5:\penalty0 1035--1038, 1963.
\newblock ISSN 0197-6788.

\bibitem[Tikhonov and Arsenin(1977)]{Tikhonov1977}
A.~N. Tikhonov and V.~Y. Arsenin.
\newblock \emph{Solutions of ill-posed problems}.
\newblock V. H. Winston \& Sons, 1977.

\bibitem[Trefethen and Bau(1997)]{Trefethen1997}
L.~N. Trefethen and D.~Bau.
\newblock \emph{Numerical Linear Algebra}.
\newblock Other Titles in Applied Mathematics. SIAM, 1997.
\newblock ISBN 9780898719574.

\bibitem[Verhulst(1996)]{Verhulst1996}
F.~Verhulst.
\newblock \emph{Nonlinear Differential Equations and Dynamic Systems}.
\newblock 1996.
\newblock ISBN 978-3-540-60934-6.
\newblock \doi{10.1007/978-3-642-61453-8}.

\bibitem[Wibisono et~al.(2016)Wibisono, Wilson, and Jordan]{WiWiJo16}
A.~Wibisono, A.~Wilson, and M.~Jordan.
\newblock A variational perspective on accelerated methods in optimization.
\newblock \emph{Proceedings of the National Academy of Sciences}, 113\penalty0
  (47):\penalty0 E7351--E7358, 2016.

\bibitem[Xiao et~al.(2017)Xiao, Rasul, and Vollgraf]{FashionMNIST}
H.~Xiao, K.~Rasul, and R.~Vollgraf.
\newblock {Fashion-MNIST}: a novel image dataset for benchmarking machine
  learning algorithms, 2017.

\bibitem[Yoshida(1990)]{Yoshida1990}
H.~Yoshida.
\newblock Construction of higher order symplectic integrators.
\newblock \emph{Physics Letters A}, 150\penalty0 (5):\penalty0 262--268, 1990.
\newblock ISSN 0375-9601.
\newblock \doi{10.1016/0375-9601(90)90092-3}.

\bibitem[Zare and Szebehely(1975)]{Zare1975}
K.~Zare and V.~G. Szebehely.
\newblock Time transformations in the extended phase-space.
\newblock \emph{Celestial mechanics}, 11:\penalty0 469--482, 1975.

\bibitem[Zhang and Sra(2016)]{Sra2016}
H.~Zhang and S.~Sra.
\newblock First-order methods for geodesically convex optimization.
\newblock In \emph{29th Annual Conference on Learning Theory}, pages
  1617--1638, 2016.

\bibitem[Zhang and Sra(2018)]{Sra2018}
H.~Zhang and S.~Sra.
\newblock An estimate sequence for geodesically convex optimization.
\newblock In \emph{Proceedings of the 31st Conference On Learning Theory},
  volume~75 of \emph{Proceedings of Machine Learning Research}, pages
  1703--1723, 2018.

\bibitem[Zhang et~al.(2018)Zhang, Mokhtari, Sra, and Jadbabaie]{Jadbabaie2018}
J.~Zhang, A.~Mokhtari, S.~Sra, and A.~Jadbabaie.
\newblock Direct {R}unge-{K}utta discretization achieves acceleration.
\newblock In \emph{Advances in Neural Information Processing Systems},
  volume~31. Curran Associates, Inc., 2018.

\bibitem[Zhong et~al.(2019)Zhong, Dey, and Chakraborty]{Zhong2020}
Y.~D. Zhong, B.~Dey, and A.~Chakraborty.
\newblock Symplectic {ODE}-{N}et: learning {H}amiltonian dynamics with control.
\newblock In \emph{International Conference on Learning Representations}, 2019.

\end{thebibliography}
\bibliographystyle{plainnat}

\end{document}